\documentclass{amsart}
\usepackage{etoolbox}
\usepackage[margin=1in,marginparwidth=0.8in]{geometry}
\usepackage{amsmath,amssymb,amsthm,mathtools}
\usepackage[pdfusetitle,unicode,hidelinks]{hyperref}
\usepackage[utf8]{inputenc}
\usepackage[T1]{fontenc}
\IfFileExists{cm-super-t1.map}{}{\usepackage{lmodern}}
\allowdisplaybreaks[2]
\makeatletter
\patchcmd{\@setabstracta}{\skip@20\p@}{\skip@20.708\p@}{}{}
\renewcommand{\@pnumwidth}{2em}
\makeatother

\newcommand{\DOI}[1]{\begingroup
  \urlstyle{same}DOI~\href{https://doi.org/#1}{\nolinkurl{#1}}\endgroup}
\numberwithin{equation}{section}
\newtheorem{theorem}{Theorem}[section]
\newtheorem{lemma}[theorem]{Lemma}
\newtheorem{proposition}[theorem]{Proposition}
\newtheorem{corollary}[theorem]{Corollary}
\theoremstyle{definition}

\theoremstyle{remark}
\newtheorem{remark}[theorem]{Remark}

\newcommand{\A}{\mathbb{A}}
\newcommand{\Pp}{\mathbb{P}}
\newcommand{\Z}{\mathbb{Z}}
\newcommand{\Q}{\mathbb{Q}}
\newcommand{\Gm}{\mathbb{G}_m}
\newcommand{\one}{\mathbf1}
\newcommand{\bx}{\overline{\square}}
\newcommand{\et}{\mathrm{\acute et}}
\newcommand{\letop}{\mathrm{l\acute et}}
\newcommand{\detop}{\mathrm{d\acute et}}
\newcommand{\setop}{\mathrm{s\acute et}}
\newcommand{\ket}{\mathrm{k\acute et}}
\newcommand{\Nis}{\mathrm{Nis}}
\newcommand{\dNis}{\mathrm{dNis}}
\newcommand{\ltr}{\mathrm{ltr}}
\newcommand{\triv}{\mathrm{triv}}
\newcommand{\Spec}{\operatorname{Spec}}
\newcommand{\Cor}{\operatorname{Cor}}
\newcommand{\lCor}{\operatorname{lCor}}
\newcommand{\lSm}{\operatorname{lSm}}
\newcommand{\Shv}{\operatorname{Shv}}
\newcommand{\Hom}{\operatorname{Hom}}
\newcommand{\RHom}{\mathrm R\operatorname{Hom}}
\newcommand{\R}{\mathrm R}
\newcommand{\cofib}{\operatorname{cofib}}
\newcommand{\fib}{\operatorname{fib}}
\newcommand{\DM}{\mathrm{DM}}
\newcommand{\logDM}{\mathrm{logDM}}
\newcommand{\eff}{\mathrm{eff}}
\newcommand{\Loc}{\operatorname{Loc}}
\DeclareMathOperator{\coker}{coker}

\newcommand{\SH}{\mathrm{SH}}
\newcommand{\logSH}{\mathrm{logSH}}
\newcommand{\Sp}{\mathrm{Sp}}
\newcommand{\sNis}{\mathrm{sNis}}
\newcommand{\Sch}{\mathrm{Sch}}
\newcommand{\ParkSch}{\mathrm{Sch}_{\mathrm P}}
\newcommand{\MGL}{\mathrm{MGL}}
\newcommand{\KGL}{\mathrm{KGL}}
\newcommand{\HZ}{\mathrm H\Z}
\newcommand{\Log}{\operatorname{\ell og}}
\newcommand{\op}{\mathrm{op}}
\newcommand{\idem}{\mathrm{idem}}
\newcommand{\PSh}{\operatorname{PSh}}
\newcommand{\Pair}{\mathcal P}
\newcommand{\Ind}{\operatorname{Ind}}
\DeclareMathOperator{\Map}{Map}
\DeclareMathOperator{\Sm}{Sm}
\DeclareMathOperator{\im}{im}
\DeclareMathOperator{\Frac}{Frac}
\DeclareMathOperator{\Th}{Th}
\DeclareMathOperator*{\colim}{colim}

\title[Logarithmic Suslin Rigidity and Comparisons]
  {Logarithmic Suslin Rigidity and Comparisons of logarithmic motives}
\author{Xuyang Chen}
\email{xuyangchen19@fudan.edu.cn}
\date{November 29, 2026}
\hypersetup{
  pdftitle={Logarithmic Suslin Rigidity and Comparisons of logarithmic motives},
  pdfauthor={Xuyang Chen}
}

\begin{document}
\begin{abstract}
Over every field we prove logarithmic Suslin rigidity for sheaves
and unbounded complexes with log transfers and prime-to-characteristic
torsion coefficients. Dividing \'etale and log \'etale effective
motives agree for arbitrary coefficients. For coefficients killed by
an invertible integer, finite \'etale schemes already generate both
categories. Over a perfect field we give an independent proof of
Merici's full faithfulness for restriction to trivial log structures
on the homotopy heart with transfers, and prove it without transfers.
Over a quasi-compact quasi-separated ordinary base, projective-line
projection invariance implies invariance under all projective-space
projections. We obtain a new presentation of effective logarithmic
spectra and an embedding of ordinary projective-line local spectra.
A filtration by marked coordinates proves full faithfulness of
restriction on effective spectra. Ordinary toric modifications
characterize its essential image, and stabilization gives the
fully faithful comparison with Annala--Iwasa spectra.
We prove Park's integral logarithmicity conjecture for algebraic
cobordism. Finite-presentation continuity extends Park's nonconnective
\(K\)-theory comparison on separated finite-dimensional Noetherian
schemes to all quasi-compact quasi-separated ordinary schemes.
\end{abstract}
\maketitle
\par\kern4.487pt
{\setcounter{tocdepth}{1}%
\tableofcontents
}

\section{Introduction}\label{sec:introduction}

Ordinary motivic homotopy theory inverts the projections
\(X\times_S\A^1_S\to X\). Logarithmic motivic homotopy theory
uses the interval \(\bx_S=(\Pp^1_S,\{\infty\})\), with the
point at infinity recorded by the log structure
\cite{BPO,BPOh}. We write \(\underline X\) for the underlying
scheme of a log scheme \(X\) and \(X^\circ\) for the open locus
where its log structure is trivial. In particular,
\(\underline{\bx_S}=\Pp^1_S\) and \(\bx_S^\circ=\A^1_S\).
This distinction gives two ways to compare the theories, by
taking the underlying scheme or by taking the interior.
The logarithmic theory also accommodates cohomology theories
such as Hodge cohomology. For example, for an affine
\(k\)-scheme \(Y\) the degree-zero coherent cohomology of
\(Y\times_k\A^1_k\) is \(\Gamma(Y,\mathcal O_Y)[t]\).
\cite[\S\S1.2--1.3 and Chapter~9]{BPO}

The results on motives with transfers are over a field \(k\).
The comparisons of spectra are over a quasi-compact,
quasi-separated (qcqs) scheme \(S\) with trivial log structure.
Quasi-separatedness ensures that intersections of quasi-compact
open subschemes are quasi-compact. The sites and categories
appearing in the statements below are defined in
Section~\ref{sec:definitions}.

\subsection{Rigidity and change of topology}

Let \(R\) be a commutative coefficient ring. Write
\(\detop\) and \(\letop\) for the dividing \'etale and log
\'etale topologies on \(\lSm/k\), the site of log-smooth
schemes over \(k\). The former adds dividing covers, such as
log blow-ups, to strict \'etale covers. The latter also
admits Kummer \'etale covers, including
\[
  (\A^1_k,\{0\})\longrightarrow(\A^1_k,\{0\}),\qquad
  t\longmapsto t^m,\qquad m>0\ \text{invertible in }k.
\]
For \(t=\detop,\letop\), the category
\(\logDM_t^\eff(k,R)\) is obtained from the unbounded derived
category of \(t\)-sheaves with finite log transfers by localizing
at the maps represented by \(X\times_k\bx_k\to X\).
We denote the motive represented by \(X\) by \(M_t(X)\).
The notation \(\DM_\et^\eff(k,R)\) refers to the corresponding
ordinary category, using ordinary finite transfers and the
\(\A^1\)-projections. Finite log transfers are defined by the
admissible correspondences of Section~\ref{sec:definitions}.

\subsubsection{Sheaf and derived rigidity}

Write \(k_\et\) for the small \'etale site of \(k\), with
its basis of finite \'etale \(k\)-schemes, and \(D(k_\et,R)\)
for the unbounded derived category of its sheaves of
\(R\)-modules. For \(nR=0\) with \(n\) invertible in \(k\),
ordinary \'etale rigidity gives
\(D(k_\et,R)\simeq\DM_\et^\eff(k,R)\).
Its sheaf form says that an \(\A^1\)-invariant \'etale sheaf
with transfers and prime-to-characteristic torsion is pulled
back from \(k_\et\).
\cite[Theorem~4.4]{SV96} \cite[Theorems~7.20 and~9.35]{MVW}
\cite[Theorem~4.5.2]{CD}

In the logarithmic setting, the inclusion
\(\eta:k_\et\hookrightarrow\lSm/k\) treats a finite \'etale
scheme as a scheme with trivial log structure. For \(nR=0\)
with \(n\) invertible in \(k\), Binda, Park, and \O stv\ae r
construct the fully faithful functor
\begin{equation}\label{eq:intro-eta}
 \eta_\sharp:D(k_\et,R)\longrightarrow
            \logDM_{\letop}^\eff(k,R).
\end{equation}
On sheaves it pulls back from \(k_\et\) and equips the
pullback with its canonical log transfers. Their
Conjectures~1.3.3 and~8.5.8 ask whether \eqref{eq:intro-eta}
is essentially surjective.
\cite[Propositions~8.5.4 and~8.5.6--8.5.7 and
Conjectures~1.3.3 and~8.5.8]{BPO}

Their sheaf conjecture asks whether the same pullback
construction recovers every log \'etale sheaf \(F\) with
log transfers and torsion of order invertible in \(k\).
The strict invariance hypothesis is the equivalence
\[
 \R\Gamma_{\letop}(X,F)\xrightarrow{\sim}
 \R\Gamma_{\letop}(X\times_k\bx_k,F)
 \qquad (X\in\lSm/k).
\]
The property concerns derived sections in every degree.
For instance, if \(k\) is algebraically closed and \(n>1\)
is invertible in \(k\), then for the constant sheaf one has
\[
 H^2_{\setop}(\bx_k,\Z/n)\simeq H^2_\et(\Pp^1_k,\Z/n)\ne0,
 \qquad H^2_\et(\Spec k,\Z/n)=0.
\]
The log \'etale pullback from the small site is strictly
interval-invariant.
\cite[Definition~8.5.1, Proposition~8.5.4, Example~8.5.5,
and Conjecture~8.5.9]{BPO}

\medskip
\noindent\textbf{Sheaf rigidity (Theorem~\ref{thm:sheaf}).}
For every such \(F\) the counit
\[
 \eta_\sharp\eta^*F\xrightarrow{\sim}F
\]
is an isomorphism. Here \(\eta^*\) is restriction to \(k_\et\).
It is enough that sections locally be killed by integers
invertible in \(k\), possibly of different orders.
Sections~\ref{sec:perfect} and~\ref{sec:perfection} prove this
using admissible correspondences between marked projective
lines and then passage to a perfect field. The curve argument
is related to affine-line invariance for torsion reciprocity
sheaves, where the relations are imposed using divisors with
multiplicities.
\cite[Definition~2.1.2]{KSY} \cite[Theorem~3.5]{BCKS}
\cite[Theorem~6.3]{Saito}

\medskip
\noindent\textbf{Unbounded derived rigidity (Theorem~\ref{thm:derived}).}
If \(nR=0\) with \(n\) invertible in \(k\), the functor
\eqref{eq:intro-eta} is an equivalence. Moreover, for
\(t=\detop,\letop\) there are compatible equivalences
\[
 I_t:\DM_\et^\eff(k,R)\xrightarrow{\sim}\logDM_t^\eff(k,R),
 \qquad
 J_t\bigl(M_t(X)\bigr)\simeq M_\et^{\mathrm{ord}}(X^\circ).
\]
Here \(I_t\) gives ordinary smooth schemes the trivial
log structure, \(J_t\) is its inverse, and
\(M_\et^{\mathrm{ord}}(X^\circ)\) is the ordinary motive of
the interior. The proof uses a cyclic action on the unmarked
projective line to show \(M_t(\A^1_k)\simeq M_t(\Spec k)\).
The two-chart decomposition of \(\bx_k\) then removes the
boundary. Park's separate rigidity theorem is formulated
after projective-line stabilization in categories without
transfers. \cite[Theorems~9.7--9.8]{ParkEtale}

\subsubsection{Comparison of the two \'etale topologies}

Binda, Park, and \O stv\ae r ask whether the finer topology
gives the same effective motives with arbitrary coefficients
\cite[Conjecture~8.6.1]{BPO}. A direct transfer argument
encounters the admissibility issue for the transpose of a
ramified Kummer map \cite[Remark~8.4.2]{BPO}.

\medskip
\noindent\textbf{Change of topology (Theorem~\ref{thm:topology}).}
For every \(k\) and \(R\), change of topology induces an
equivalence
\[
 T_R:\logDM_{\detop}^\eff(k,R)
       \xrightarrow{\sim}\logDM_{\letop}^\eff(k,R).
\]
The proof first treats a root cover of degree invertible in
\(R\) by using the transpose on the ordinary affine chart of
\(\bx_k\) and Mayer--Vietoris. The derived rigidity theorem
supplies the comparison with \(\mathbb F_\ell\)-coefficients
for every prime \(\ell\) invertible in \(k\). These comparisons
together detect the integral kernel of \(T_{\Z}\).

\subsubsection{Proper generators}

Let \(\logDM_{t,\mathrm{prop}}^\eff(k,R)\) be the localizing
subcategory generated by \(M_t(X)\) for log-smooth \(X\)
with \(\underline X\) proper over \(k\). A localizing
subcategory here is triangulated and closed under arbitrary
coproducts. The earlier Nisnevich comparison assumes
resolution of singularities. The \'etale conjecture is
formulated for \(R=\Z/n\) with \(n\) invertible in \(k\).
\cite[Definition~5.2.8, Theorem~8.2.17, and
Conjecture~8.6.3]{BPO}

\medskip
\noindent\textbf{Proper comparison (Theorem~\ref{thm:proper}).}
For any \(R\) with \(nR=0\) and \(n\) invertible in \(k\),
the two canonical comparisons
\[
 \logDM_t^\eff(k,R)\ \longleftarrow\
 \logDM_{t,\mathrm{prop}}^\eff(k,R)\
 \longrightarrow\ \DM_\et^\eff(k,R),
 \qquad t=\detop,\letop,
\]
are equivalences. The left arrow is inclusion and the right
arrow takes interiors. In fact, the motives of finite
\'etale \(k\)-schemes with trivial log structure already
generate the full category for each of the two topologies.

\subsection{Restriction and representability}

\subsubsection{Restriction on the Nisnevich homotopy heart}

The dividing Nisnevich topology adds dividing covers to
strict Nisnevich covers. Let
\(\mathrm{CI}_{\dNis}^{\log}(k,\Lambda)\) and
\(\mathrm{CI}_{\dNis}^{\ltr}(k,\Lambda)\) denote the categories
of strictly \(\bx\)-invariant dividing Nisnevich sheaves of
\(\Lambda\)-modules, respectively without and with log
transfers. They are the hearts of the homotopy
\(t\)-structures of Binda and Merici. For a perfect field \(k\)
and a commutative ring \(\Lambda\),
Theorem~\ref{thm:restriction-sheaves} proves that
\[
\begin{aligned}
 \omega_\sharp:\mathrm{CI}_{\dNis}^{\log}(k,\Lambda)
    &\longrightarrow \Shv_\Nis(k,\Lambda),\\
 \omega_\sharp:\mathrm{CI}_{\dNis}^{\ltr}(k,\Lambda)
    &\longrightarrow \Shv_\Nis^{\mathrm{tr}}(k,\Lambda)
\end{aligned}
\]
are fully faithful. Both functors evaluate a log sheaf at
ordinary smooth schemes with the trivial log structure.
Merici established the integral case with transfers. His result
extends to coefficient-ring actions by uniqueness of the lift.
Section~\ref{sec:sheaf-restriction} gives a proof that also
applies without transfers. It uses square scaling to
characterize the sections that extend from \({\Gm}_K\) to
\((\Pp^1_K,\{0,\infty\})\).
\cite[Theorems~5.7 and~5.10 and Conjecture~7.3]{BM}
\cite[Proposition~0.1 and Conjecture~0.2]{BMerr}
\cite[Theorem~4.4]{Merici}

\subsubsection{Projective-line locality}

An ordinary Nisnevich sheaf \(F\) on \(\Sm/S\) is
\(\Pp^n\)-local if pullback along
\(X\times_S\Pp^n_S\to X\) induces an equivalence for every
smooth \(X/S\). Ayoub studied \(\Pp^1\)-localization, and
Binda, Park, and \O stv\ae r ask whether it agrees with the
localization at all \(\Pp^n\)-projections after simplicial-circle
stabilization. Theorem~\ref{thm:unstable-projective} proves
that every \(\Pp^1\)-local Nisnevich sheaf of spaces is
\(\Pp^n\)-local for every \(n\geq1\).
Consequently Corollary~\ref{thm:projective} identifies the
canonical functor
\[
 \SH^{\eff,\Pp^1}(S)\xrightarrow{\sim}
 \SH^{\eff,\Pp^\bullet}(S)
\]
as a symmetric monoidal equivalence natural in \(S\).
\cite{Ayoub,Maity} \cite[Definition~4.0.3 and
Question~4.0.4]{BPOh}

Write \(\mathcal C_S=\Shv_{\sNis}(\mathrm{SmlSm}/S,\Sp)\)
for the strict Nisnevich sheaves of spectra on the smooth-pair
site of Section~\ref{sec:spectra-definitions}. In the display,
\((\Pp^\bullet,\Pp^{\bullet-1})\) denotes the pairs
\((\Pp^n_S,\Pp^{n-1}_S)\) with the hyperplane at infinity
marked, for \(n\geq1\). The brackets invert product projections
with the indicated pairs and the composites of admissible
blow-ups along smooth centers when \(\mathrm{SmAdm}\) is listed.
Theorem~\ref{thm:four-log-presentations} gives canonical
equivalences
\[
\begin{aligned}
 \logSH^\eff(S)
 &\simeq \mathcal C_S[\bx^{-1},\mathrm{SmAdm}^{-1}]\\
 &\simeq \mathcal C_S[(\Pp^\bullet,\Pp^{\bullet-1})^{-1}]
 \simeq \mathcal C_S[\bx^{-1}].
\end{aligned}
\]
The last presentation shows that after strict Nisnevich
descent, invariance under the marked line forces invariance
under the other marked projective-space projections and
under the admissible blow-ups. The first three presentations
are due to Binda, Park, and \O stv\ae r.
\cite[Theorem~1.3.1]{BPOh}

Taking interiors gives the fully faithful functor
\(\omega^*:\SH^\eff(S)\hookrightarrow\logSH^\eff(S)\) of
Binda, Park, and \O stv\ae r, with
\((\omega^*F)(Y)=F(Y^\circ)\). They express the embedding as
\[
 \bigl(\A^1\text{-local homotopy theory}\bigr)
 \ \subset\
 \bigl(\bx\text{-local homotopy theory}\bigr).
\]
\cite[Construction~4.0.8]{BPOh}
The underlying-scheme functor gives a fully faithful comparison
on effective spectra
\[
 \rho^*:\SH^{\eff,\Pp^1}(S)\hookrightarrow\logSH^\eff(S),
 \qquad \rho^*F(Y)=F(\underline Y).
\]
Thus we obtain the projective-line version of the slogan
\[
 \bigl(\Pp^1\text{-local homotopy theory}\bigr)
 \ \subset\
 \bigl(\bx\text{-local homotopy theory}\bigr).
\]
Here both categories are stabilized with respect to the
simplicial circle \(S^1\). The left-hand locality refers to
the ordinary projective-line projections. See
Corollary~\ref{cor:projective-embedding} and
\cite[Construction~4.0.6]{BPOh}.

\subsubsection{Restriction of spectra and the essential image}

Giving smooth schemes the trivial log structure induces the
restriction functor in Question~4.0.17. We prove its full
faithfulness, and also that of the right adjoint in the
stable comparison with Annala--Iwasa spectra.
\[
\begin{aligned}
 \lambda^*:\logSH^\eff(S)&\hookrightarrow
                  \Shv_\Nis(\Sm/S,\Sp),\\
 \lambda^*:\logSH(S)&\hookrightarrow \mathrm{MS}_S .
\end{aligned}
\]
The second target \(\mathrm{MS}_S\) is defined in
Section~\ref{sec:ms-log-comparison} by projective-line
stabilization of Zariski sheaves of spectra satisfying
elementary blow-up excision. See
Theorem~\ref{thm:q417-full-scaling} and
Corollary~\ref{cor:ms-log-comparison}.
\cite[Construction~4.0.16, Question~4.0.17, and
Construction~4.0.19]{BPOh}

Corollary~\ref{cor:q418-essential-image} gives an ordinary-scheme
criterion for a Nisnevich sheaf of spectra \(H\) to lie in
the first essential image. It requires the maps
\[
 H(X)\longrightarrow \bigl|C_\bullet^{(r)}(X,H)\bigr|,
 \qquad
 C_m^{(r)}(X,H)=
   \colim_{\Sigma\in I_{r,m}^{\op}} H(X\times_S Y_{\Sigma,S}),
\]
to be equivalences for every smooth \(X/S\) and \(r\geq0\).
The schemes \(Y_{\Sigma,S}\) are smooth toric modifications of
\((\Pp^1_S)^{r+m}\), with index category \(I_{r,m}\) and
cubical realization defined in Section~\ref{sec:toric}.
Over Park's bases this tests whether the logarithm of \(H\)
recovers \(H(X)\) on \(X\times_S\bx_S^r\). We verify the
criterion for restrictions of logarithmic spectra by
filtering the cubical diagram by marked coordinates, and
then use continuity for general qcqs bases.
\cite[Remark~4.0.18]{BPOh}
\cite[Definition~5.26, Proposition~8.4, and Theorem~8.5]{ParkI}

\subsubsection{Algebraic cobordism and \(K\)-theory}

Park's base class \(\ParkSch\) consists of separated
Noetherian schemes of finite Krull dimension.
Conjecture~7.16 asks whether the algebraic cobordism
spectrum \(\MGL_S\) is logarithmic over regular bases in
this class. For an ordinary motivic spectrum \(E\), its
section sheaf in Tate degree \(q\) is
\[
 E^{(q)}(X)=\Map_{\SH(S)}
     \bigl(\Sigma^\infty_{\Pp^1}X_+,\Sigma^{2q,q}E\bigr).
\]
The log extension \(\Log\) and the Tate suspensions are
defined in Sections~\ref{sec:toric} and~\ref{sec:cobordism}.
Theorem~\ref{thm:park-mgl} gives the equivalence
\[
 \MGL_S^{(q)}(X)
 \xrightarrow{\sim}
 \Log\bigl(\MGL_S^{(q)}\bigr)(X\times_S\bx_S^r)
 \qquad(S\in\ParkSch,\ X\in\Sm/S,\ q\in\Z,\ r\geq0).
\]
Its proof uses smooth proper duality, a uniform connectivity
bound on the realized comparison cofiber, and Park's
calculation for integral motivic cohomology.
\cite[Definitions~5.26 and~5.29, Theorem~7.9, and
Conjecture~7.16]{ParkI} \cite[Theorem~1.1]{ParkII}

Let \(\mathrm{logKGL}_{\Z}\) be the logarithmic Bott-periodic
spectrum over \(\Spec\Z\) and let \(f_S:S\to\Spec\Z\) be the
structural morphism. Conjecture~6.5.10 predicts that for
ordinary \(S\) the spectrum
\[
 \log K(S)=
 \Map_{\logSH(S)}(\one_S,f_S^*\mathrm{logKGL}_{\Z})
\]
agrees with the nonconnective \(K\)-theory \(K(S)\) of
perfect complexes on \(S\). Park's representability and
base change give this comparison on \(\ParkSch\).
Theorem~\ref{thm:log-k} extends the natural equivalence
\(\log K(S)\simeq K(S)\) to every qcqs ordinary \(S\),
using finite-presentation continuity and Zariski descent.
\cite[Definition~6.5.9 and Conjecture~6.5.10]{BPOh}
\cite[Theorems~7.10, 8.11, 8.13 and Proposition~8.12]{ParkI}

Sections~\ref{sec:definitions} and~\ref{sec:sheaf-restriction}
set up the categories and prove the comparison for sheaves.
Sections~\ref{sec:perfect}--\ref{sec:consequences} treat
rigidity and change of topology. Section~\ref{sec:projective}
proves the projective-space comparisons.
Sections~\ref{sec:continuity} and~\ref{sec:toric} establish
the continuity and toric constructions used in the
spectral results of Sections~\ref{sec:spectral},
\ref{sec:cobordism}, and~\ref{sec:k-theory}.

\section{Sites, categories, and statements on rigidity}\label{sec:definitions}

\subsection{Logarithmic schemes, transfers, and topologies}

A log structure on a scheme \(\underline X\) is a sheaf of commutative
monoids \(M_X\) on its small \'etale site, with a map
\(M_X\to(\mathcal O_{\underline X},\cdot)\) identifying the inverse
image of \(\mathcal O_{\underline X}^{\times}\) with that sheaf of
units. Its characteristic monoid is
\(\overline M_X=M_X/\mathcal O_{\underline X}^{\times}\).
The trivial log structure is \(\mathcal O_{\underline X}^{\times}\).
We write \(X\) for the log scheme and \(X^\circ\) for the open set
where its log structure is trivial.

Pullback of a log structure pulls back its monoid map and adjoins the
units on the new scheme, identifying elements mapped to units with
those units. A log morphism is a scheme morphism with a compatible
map from the pulled-back target log structure to the source log
structure. It is \emph{strict} when this map is an isomorphism.
A chart is a monoid map \(P\to\mathcal O_{\underline X}(V)\)
that produces the log structure by this construction.

Write \(P^{\mathrm{gp}}\) for the group obtained from a commutative
monoid \(P\) by adjoining additive inverses. The monoid is
\emph{fine and saturated}, abbreviated fs, if it is finitely
generated, embeds in \(P^{\mathrm{gp}}\), and contains every element
of \(P^{\mathrm{gp}}\) having a positive multiple in \(P\).
A log structure is fs when it has fs charts \'etale locally.
Following \cite[\S1.4]{BPO}, we require the charts to exist on
Zariski open neighborhoods. These are \emph{Zariski log structures}. The use of
\'etale log structures in Vidal's theorem is addressed explicitly
in Section~\ref{sec:perfection}.

A log smooth morphism has, \'etale locally, compatible charts
\(P\to Q\), from target to source, such that the kernel and the
torsion in the cokernel of \(P^{\mathrm{gp}}\to Q^{\mathrm{gp}}\)
are finite of invertible orders, and the induced scheme morphism
to the chart fiber product is smooth. For log \'etaleness the
whole cokernel is finite of invertible order and the induced
scheme morphism is \'etale. These are the chart criteria
\cite[Theorem~3.5]{Kato}. Additional chart generators can make the
scheme map \'etale in the log smooth case.

The category \(\lSm/k\) consists of fs log smooth schemes over the
trivial-log field \(k\), separated and of finite type on underlying
schemes. Let \(\mathrm{SmlSm}/k\) be its full subcategory with smooth
underlying schemes. Objects of this subcategory are the pairs
\((Y,D)\), where \(Y\) is smooth and \(D\) is a reduced strict
normal-crossing divisor: locally its components are distinct
coordinate hyperplanes. Its chart is generated by the equations of
these hyperplanes, together with the units. The locus \(X^\circ\) of any \(X\in\lSm/k\)
is smooth and dense. A log blow-up makes its underlying scheme
smooth \cite[Proposition~A.10.2 and Lemma~A.5.10]{BPO}.

An elementary finite log correspondence from \(X\) to \(Y\) is an
integral closed subscheme of \(\underline X\times_k\underline Y\),
finite and surjective over a connected component of \(\underline X\),
with this extra condition: its normalization, given the log
structure pulled back from \(X\), maps to \(Y\) as a log scheme.
This condition is called \emph{admissibility}. The group
\(\lCor_k(X,Y)\) is the free abelian group on these correspondences.
Restriction to the trivial loci injects it into
\(\Cor_k(X^\circ,Y^\circ)\) and respects composition
\cite[Definition~2.1.1 and \S2.3]{BPO}. Composition can therefore
be checked on the interiors using ordinary finite correspondences.
For trivial log structures all finite correspondences are allowed.

For a commutative unital ring \(R\), a presheaf with log transfers is
an additive contravariant functor from this correspondence category
to \(R\)-modules. The represented presheaf is
\[
 R_{\ltr}(X)=\lCor_k(-,X)\otimes_{\Z}R.
\]
Graphs of log morphisms give its usual restriction maps.

A Nisnevich cover is an \'etale cover in which each point has a lift
with the same residue field. A Kummer log morphism is one whose
maps on geometric characteristic monoids are injective and whose
codomain elements have positive multiples in the image. A Kummer
\'etale morphism is log \'etale and Kummer.

A dividing cover is a surjective proper log \'etale monomorphism
in the fs category. Log blow-ups are examples. A log blow-up makes
a finitely generated ideal of the log monoid locally principal,
and is universal for this property. Here a monoid ideal is a subset
stable under adding arbitrary monoid elements. Dividing covers
are isomorphisms over the interior. We follow the terminology of
\cite[Remark~2.2.11]{BPOh}.

We use the strict \'etale topology \(\setop\), the Kummer \'etale
topology \(\ket\), the dividing Nisnevich topology \(\dNis\), the
dividing \'etale topology \(\detop\), and the log \'etale topology
\(\letop\). Strict \'etale covers have the pulled-back log structure.
The last three topologies are generated respectively by strict
Nisnevich, strict \'etale, and Kummer \'etale covers together with
dividing covers. A log \'etale covering family is jointly
surjective after every fs base change.

The topologies are compatible with transfers: the underlying
sheafification has a unique transfer action and is exact
\cite[\S4.2 and Theorem~4.5.7]{BPO}. Denote it by \(a_t\) and
the abelian category of sheaves by \(\Shv_t^{\ltr}(k,R)\).
Forgetting transfers is written \(\gamma^*\).
On a trivial-log object the small sites are the ordinary sites.
For \(t=\dNis,\detop,\letop\), restriction from \(\lSm/k\) to
\(\mathrm{SmlSm}/k\) is an equivalence on sheaves with transfers
\cite[Proposition~4.7.5]{BPO}.

Write \(D(\mathcal A)\) for the unbounded derived category of an
abelian category \(\mathcal A\). We retain its models by complexes
and their chain homotopies when forming diagrams and derived
colimits. A cohomological shift \([1]\) moves a module in degree
zero to degree \(-1\). A cofiber is a mapping cone, and a fiber is
that cone shifted by \([-1]\). The realization of a simplicial
complex is its derived colimit, computed by the direct-sum total
complex of its normalized chains.

The \v{C}ech nerve of a cover \(Y\to X\) has
\(Y^{[a]}=Y^{\times_X(a+1)}\) in degree \(a\), using fs products.
A hypercover is an augmented simplicial covering in which every
degree covers the matching object determined by its lower faces.
The derived categories of sheaves can be presented from complexes
of presheaves by imposing these hypercover relations
\cite[\S4.3 and Appendix~B]{BPO}.

For a complex \(F\), derived sections are denoted
\(\R\Gamma_t(X,F)\). If \(h_{t,R}(X)=a_tR_{\ltr}(X)\), then
\begin{equation}\label{eq:eval}
 \RHom(h_{t,R}(X),F)
 \simeq\R\Gamma_t(X,F).
\end{equation}
The right side is the cohomology on the small site after forgetting
transfers. This natural comparison is induced by the derived
adjunctions for forgetting transfers and restricting to the small
site \cite[Propositions~3.2.4 and~4.3.2]{BPO}. The use of mapping
complexes in \eqref{eq:eval} preserves group actions and diagrams.

\subsection{Effective motives over a field}\label{sec:coeff-fields}

Set
\[
 A=(\A^1,\triv),\quad P=(\Pp^1,\triv),\quad
 \bx=(\Pp^1,\{\infty\}),\quad C=(\A^1,\{0\}).
\]
Products are fs products, the categorical fiber products among fs
log schemes. On charts one takes the integral saturation of the
pushout monoid. With a trivial-log scheme their underlying
schemes are the usual products.
The category \(\logDM_t^\eff(k,R)\) is the localization of
\(D(\Shv_t^{\ltr}(k,R))\) at
\[
 h_{t,R}(X\times\bx)\longrightarrow h_{t,R}(X)
 \qquad(X\in\lSm/k).
\]
The image of \(h_{t,R}(X)\) is \(M_t(X)\), and
\(\one_t=M_t(\Spec k)\). Tensor products are derived, preserve
colimits in each variable, and satisfy
\[
 M_t(X)\otimes M_t(Y)\simeq M_t(X\times_kY).
\]
These are the monoidal constructions of
\cite[\S\S5.2--5.3 and Appendix~B]{BPO}.

The fully faithful right adjoint to localization identifies motives
with the complexes \(F\) for which
\[
 \R\Gamma_t(X,F)\xrightarrow{\sim}
 \R\Gamma_t(X\times\bx,F)
\]
for every \(X\). We call these complexes \emph{strictly
\(\bx\)-invariant}. The same term for a sheaf refers to its
degree-zero complex. Formula \eqref{eq:eval} holds for these
complexes with \(M_t(X)\) as source.

The ordinary category \(\DM_\et^\eff(k,R)\) uses smooth schemes,
ordinary finite transfers, the \'etale topology, and the projections
with \(\A^1\). Its represented motives are
\(M_\et^{\mathrm{ord}}(U)\).
The inclusion \(\lambda:U\mapsto(U,\triv)\) and the interior
functor \(\omega:X\mapsto X^\circ\) act on correspondences.
Their induced motivic functors will be constructed in
Section~\ref{sec:derived}.

The small site \(k_\et\) has finite \'etale \(k\)-schemes as a basis.
Its sheaves have their canonical finite transfers. Restriction to
this site is \(\eta^*\). Its left adjoint with log transfers is
\(\eta_\sharp\), constructed by additive presheaf extension and
sheafification. For torsion annihilated by an invertible integer
it is ordinary pullback equipped with its unique log transfers
\cite[Propositions~8.5.2 and~8.5.6]{BPO}. The canonical adjunction
map \(\eta_\sharp\eta^*F\to F\) is the \emph{counit}.
The same description for elementwise torsion follows from
Section~\ref{sec:perfection}.

\begin{theorem}[Sheaf rigidity {\cite[Conjecture~8.5.9]{BPO}}]\label{thm:sheaf}
Let \(k\) be a field and \(F\) a log \'etale sheaf with log transfers.
Suppose that locally each section is killed by an integer invertible
in \(k\), and that \(F\) is strictly \(\bx\)-invariant.
Then \(\eta_\sharp\eta^*F\to F\) is an isomorphism.
It respects any coefficient-ring action on \(F\).
\end{theorem}

Equivalently, the torsion hypothesis says
\(F=\varinjlim_n F[n]\), where \(F[n]=\ker(n:F\to F)\) and
the directed set uses divisibility among positive integers
invertible in \(k\). A single uniform exponent is not required.
Section~\ref{sec:perfectresult} proves the theorem over perfect fields.
Section~\ref{sec:perfection} descends the counit isomorphism along perfection.

\begin{theorem}[Derived rigidity {\cite[Conjectures~1.3.3 and~8.5.8]{BPO}}]
\label{thm:derived}
Let \(k\) be any field and \(nR=0\), with \(n>0\) invertible in \(k\).
Then the canonical functor
\[
 \eta_\sharp:D(k_\et,R)\xrightarrow{\sim}
 \logDM_{\letop}^\eff(k,R)
\]
is an equivalence. For \(t=\detop,\letop\), inclusion of trivial-log
schemes also induces compatible equivalences
\[
 I_t:\DM_\et^\eff(k,R)\xrightarrow{\sim}
       \logDM_t^\eff(k,R).
\]
Their inverses \(J_t\) send \(M_t(X)\) to
\(M_\et^{\mathrm{ord}}(X^\circ)\).
\end{theorem}

The proof is in Section~\ref{sec:ordinary-comparison}.

\begin{theorem}[Change of topology {\cite[Conjecture~8.6.1]{BPO}}]\label{thm:topology}
For every field \(k\) and every commutative unital ring \(R\), the
canonical change of topology is an equivalence
\[
 T_R:\logDM_{\detop}^\eff(k,R)\xrightarrow{\sim}
            \logDM_{\letop}^\eff(k,R).
\]
\end{theorem}

A full triangulated subcategory is called \emph{localizing} when it
is closed under arbitrary coproducts. Write
\(\logDM_{t,\mathrm{prop}}^\eff(k,R)\) for the smallest such
subcategory containing \(M_t(X)\) for all \(X\) with proper
underlying scheme. This is \cite[Definition~5.2.8]{BPO} and its
\'etale variants.

\begin{theorem}[Proper comparison {\cite[Conjecture~8.6.3]{BPO}}]\label{thm:proper}
If \(nR=0\) with \(n\) invertible in \(k\), then for
\(t=\detop,\letop\) the canonical comparisons
\[
 \logDM_t^\eff(k,R)
 \ \longleftarrow\
 \logDM_{t,\mathrm{prop}}^\eff(k,R)
 \ \longrightarrow\
 \DM_\et^\eff(k,R)
\]
are equivalences. Finite \'etale \(k\)-schemes with trivial log
structure already generate these categories. In particular this
proves the conjecture with its stated coefficients \(R=\Z/n\).
\end{theorem}

An exact functor is \emph{conservative} if it detects isomorphisms,
equivalently if it detects zero objects. A family is jointly
conservative if an object is zero whenever all its images are zero.
Linearity over \(D(\Z)\) means that objects can be tensored with
complexes of abelian groups using the derived tensor product,
compatibly with colimits. A module over a coefficient algebra
means an object with a unital associative action, taken in the
models by complexes and their localizations.

\begin{lemma}[Derived coefficients]\label{lem:coefficients}
The categories of motives above are linear over \(D(\Z)\).
For a ring \(R\), their \(R\)-coefficient versions are the
categories of modules over the coefficient algebra
\(\one\otimes_{\Z}^{\mathbf L}R\) in their integral versions.
In particular, after forgetting the coefficient action, extension
to \(\mathbb{F}_\ell\) sends \(E\) to
\(\cofib(\ell:E\to E)\). Extension to a localization
\(S^{-1}\Z\) inverts the maps \(s:E\to E\), \(s\in S\).
All these constructions commute with change of topology and with
the field extensions used below.
\end{lemma}

\begin{proof}
Before descent, a module over \(R\) is a complex of \(R\)-modules
whose transfer maps are \(R\)-linear. This is precisely a module
over \(R\otimes\one\) in complexes of integral presheaves.
The descent and projection relations are stable under tensoring
with representables, and these generate the presheaf category.
Localization is therefore monoidal. Its category of \(R\)-modules
is obtained by localizing the presheaf module category at the free
\(R\)-module versions of those relations. Indeed a module is local
exactly when its underlying object is local, by the free-module
adjunction. This proves the assertion about coefficients.

The two-term free resolution of \(\mathbb{F}_\ell\) gives the
cofiber formula. A localization of \(\Z\) is the filtered
colimit obtained by multiplication by the elements being inverted,
giving the second formula. Extension of coefficients sends
represented objects to the same objects with new coefficients,
since the integral correspondence groups are free abelian.
The definitions by these relations also prove the stated
compatibilities.
\end{proof}

\begin{lemma}\label{lem:pullback}
Extension of scalars \(b\subset k\) induces a colimit-preserving
functor on effective log motives for \(t=\dNis,\detop,\letop\)
when \(b\) is perfect or the extension is finite separable.
The functor sends \(M_t(X)\) to \(M_t(X_k)\).
\end{lemma}

\begin{proof}
On ordinary finite correspondences the functor is flat
pullback of cycles, which commutes with composition
\cite[\S3.7, especially Theorem~3.7.3 and Corollary~3.7.5]{SV00}.
Let \(Z\) be an elementary log correspondence. Its
normalization \(Z^\nu\) is finite since it is of finite
type over a field. A normal scheme over a perfect field is
geometrically normal, meaning it is normal after every
extension of that field \cite[Tags~038N and 038O]{Stacks}.
For a finite separable extension normality follows because
the base change is \'etale over the original scheme.
Thus \((Z^\nu)_k\) is normal in either case. The map
\((Z^\nu)_k\to Z_k\) is finite and an isomorphism over
the pullback of a dense open of \(Z\). Every generic
point of \(Z_k\) lies over the generic point of \(Z\) by
flatness. Also \(Z_k\) is reduced, by geometric
reducedness over a perfect field or by \'etaleness
in the finite separable case. It follows that the
map is the normalization of \(Z_k\). The log morphism
on \(Z^\nu\) pulls back to give the required log
morphism on each component of that normalization.
The components of the pulled-back cycle are finite and
dominant over components of the pulled-back source, as in
ordinary flat pullback of finite correspondences.
Composition now follows by restriction to the interiors
and the injectivity of that restriction.

This gives a functor on correspondences. Strict base
change of the ground field preserves the covering families
generating the three topologies, commutes with fiber
products over the covered object and preserves \(\bx\).
It therefore preserves hypercovers and the projection
relations. Derived left extension of presheaves now
gives the asserted functor.
\end{proof}

\begin{lemma}\label{lem:descent}
If \(L/k\) is finite separable and \(t=\detop\) or \(\letop\),
then extension of scalars \(q^*\) to \(L\) detects zero objects.
\end{lemma}

\begin{proof}
Viewing \(L\)-schemes as \(k\)-schemes gives a functor on log
correspondences. The cycle of a correspondence from \(U\)
to \(V\) over \(L\) is included from the open and closed
subscheme \(\underline U\times_L\underline V\) of
\(\underline U\times_k\underline V\), using the diagonal
of the finite \'etale map \(\Spec L\to\Spec k\).
The normalization and the log morphism are unchanged.
Composition is preserved because the intersections and
pushforwards are supported on the locus where the
\(L\)-structures agree. This is also the correspondence
functor given by the adjunction in
\cite[Remark~2.1.2(iv)]{BPO}. Forgetting the field preserves
covers and their matching objects in each slice over an
\(L\)-object. It sends \(U\times_L\bx_L\) to
\(U\times_k\bx_k\). Thus derived left extension induces
a colimit-preserving functor \(q_\sharp\). On representable
presheaves, its composite with \(q^*\) is product with
\(\Spec L\), naturally for correspondences. Passing to the
derived categories and localizing gives
\[
 q_\sharp q^* E\simeq E\otimes M_t(\Spec L).
\]

Suppose \(q^*E=0\). For \(a\geq0\), put
\(T_a=(\Spec L)^{\times_k(a+1)}\), the \(a\)-th term of
the \v{C}ech nerve of \(q\). A projection to one factor
provides a section of \(T_a\times_k\Spec L\to T_a\).
It follows that \(E\otimes M_t(T_a)\) is a retract of
\(E\otimes M_t(T_a)\otimes M_t(\Spec L)=0\).
By \v{C}ech descent the realization of \(M_t(T_\bullet)\)
is the unit. Tensoring with \(E\) and taking this realization
gives \(E=0\).
\end{proof}

\subsection{Spectra over an ordinary base}\label{sec:spectra-definitions}

Write \(\Sch\) for qcqs ordinary schemes, as in
\cite[\S1.7]{BPOh}. For \(S\in\Sch\), let \(\Sm/S\) be smooth
finitely presented schemes over \(S\) and let
\(\Pair_S=\mathrm{SmlSm}/S\) be the site of smooth pairs.
In a smooth pair \((Y,D)\), the divisor \(D\) is a relative
strict normal crossings divisor. Locally its equations
are distinct coordinates in a relative \'etale chart over \(S\).
The associated log structure is given by the chart
\(\mathbb N^a\to\mathcal O_Y,\ e_i\mapsto t_i\) for
the local equations \(t_i\). The chart definition applies
over nonreduced bases. For example, on
\(\Spec((k[\epsilon]/(\epsilon^2))[t])\), the function
\(t+\epsilon\) is invertible away from \(t=0\), but at
\((t,\epsilon)\) it is not a unit times a power of \(t\).
The rule using all functions invertible away from the divisor
would give a different log structure there.

Write \(\Sp\) for spectra, and use \(\Map\) for a mapping
spectrum in a stable category. Our categories and their
diagrams are understood as \(\infty\)-categories. A spectral
sheaf satisfies descent with limits in \(\Sp\).
We use ordinary sheaf descent in the strict Nisnevich
presentations. Write
\[
 \mathcal N_S=\Shv_{\Nis}(\Sm/S,\Sp),\qquad
 \mathcal C_S=\Shv_{\sNis}(\Pair_S,\Sp),\qquad
 \mathcal L_S=\logSH^{\eff}(S).
\]
The strict Nisnevich topology has Nisnevich covers with
the log structures pulled back. The effective logarithmic
category \(\mathcal L_S\) is the \(S^1\)-stable category of
dividing Nisnevich sheaves of spectra invariant under product
with \(\bx\). Its strict Nisnevich presentation is
\[
 \mathcal L_S\simeq
 \mathcal C_S[(\Pp^\bullet,\Pp^{\bullet-1})^{-1}].
\]
The brackets mean localization at the maps represented
by \(Y\times_S P_n\to Y\) for every pair \(Y\), where
\(P_n=(\Pp^n_S,\Pp^{n-1}_S)\), the hyperplane \(x_0=0\)
being marked, and \(n\geq1\)
\cite[Corollary~3.4.7]{BPOh}.
We write \(M(Y,D)\), or \(M(Y)\) when the log structure is
understood, for the represented spectrum in this category.
Write \(M_S(Y,D)\) for its image in \(\logSH(S)\), where
the pointed projective-line suspension
\(T_S=\cofib(\one_S\xrightarrow{1}M(\Pp^1_S))\) is made
tensor-invertible. In this definition the log structure on
\(\Pp^1_S\) is trivial, the arrow uses the section \(1\),
and \(\one_S=M(S)\) is the tensor unit. A cofiber in a
stable category is the pushout along the map to zero.
Tensor products of represented spectra are represented by
products of pairs.

We will use smooth-center blow-up invariance in \(\mathcal L_S\).
An admissible blow-up along a smooth center is one whose center
is contained in the boundary and is given by coordinate
equations in charts adapted to the boundary. Its new boundary
is the reduced total transform. Such blow-ups induce
equivalences on represented spectra
\cite[Definition~2.5.1 and Theorem~3.3.3]{BPOh}.
The invariance used here is a property of the localized
category. It is stronger than invariance under the dividing
covers in the definition of the topology.

\begin{lemma}[Square scaling on the unmarked line]\label{lem:unmarked-scaling}
For an ordinary qcqs base \(B\) and a unit \(a\) on \(B\), the
automorphism \(z\mapsto a^2z\) acts as the identity on
\(M(\Pp^1_B)\) in the homotopy category of \(\mathcal L_B\).
\end{lemma}

\begin{proof}
First consider a projective line over \(B\). For a unit \(c\), the
family \(z\mapsto z+c\tau\) defines a rational map
\(\Pp^1_z\times\bx_\tau\dashrightarrow\Pp^1_z\) with trivial log
structure on the target. In inverse coordinates \(u=1/z\) and
\(v=1/\tau\) it is \([v+cu:uv]\), with boundary \(v=0\). Blow up
\(u=v=0\). In the chart \(v=u w\) the map is
\([w+c:u w]\) with remaining base locus \(u=0,\ w=-c\).
The chart \(u=v h\) gives the same center, since \(c\) is a unit.
Near this center the base ideal is \((u,w+c)\). A second blow-up
resolves the map. Both centers are smooth over \(B\) and lie in
the boundary with normal crossings. Use as boundary the strict
transforms and exceptional divisors, each with multiplicity one. Their projections, followed by the interval
projection, induce equivalences on represented effective spectra
by \cite[Theorem~3.3.3]{BPOh} and interval invariance. The sections
at \(\tau=0,1\) lift to the resolved family and induce the same
inverse to the projection. This proves that \(z\mapsto z+c\) acts
as the identity on \(M(\Pp^1_B)\).

Conjugating by inversion gives the same result for the lower
elementary matrix. Write \(E_{ij}(c)\) for the corresponding
elementary matrix and put
\[
 w(b)=E_{12}(b)E_{21}(-b^{-1})E_{12}(b).
\]
Then \(w(b)w(-1)=\operatorname{diag}(b,b^{-1})\), which acts by
\(z\mapsto b^2z\). Hence square scaling acts as the identity on
the projective line.

\end{proof}

\section{Restriction on strictly invariant sheaves}\label{sec:sheaf-restriction}

Fix a perfect field \(k\) and a commutative unital coefficient ring
\(\Lambda\). Without transfers, strictly interval-invariant sheaves
are defined by the same invariance of derived sections on the
dividing Nisnevich site, using sheaves of \(\Lambda\)-modules on
ordinary log morphisms. Denote their category by
\(\mathrm{CI}_{\dNis}^{\log}(k,\Lambda)\), and denote the category
with log transfers by \(\mathrm{CI}_{\dNis}^{\mathrm{ltr}}(k,\Lambda)\).
These are the degree-zero objects for the homotopy \(t\)-structures
of \cite[Theorem~5.7]{BM}. Write \(\omega_\sharp F\) for the
ordinary Nisnevich sheaf \(V\mapsto F(V,\triv)\), keeping
the usual finite transfers when present. On trivial log
schemes the derived dividing Nisnevich sections are the
ordinary derived Nisnevich sections of this restriction.
With log transfers they are computed on the underlying sheaf
by the comparison in \eqref{eq:eval}.

\begin{theorem}[Restriction for sheaves {\cite[Conjecture~7.3]{BM}}]\label{thm:restriction-sheaves}
For every perfect \(k\) and every commutative unital \(\Lambda\),
the restriction functors
\[
\begin{aligned}
 \mathrm{CI}_{\dNis}^{\log}(k,\Lambda)
   &\xrightarrow{\ \omega_\sharp\ }\Shv_{\Nis}(k,\Lambda),\\
 \mathrm{CI}_{\dNis}^{\mathrm{ltr}}(k,\Lambda)
   &\xrightarrow{\ \omega_\sharp\ }\Shv_{\Nis}^{\mathrm{tr}}(k,\Lambda)
\end{aligned}
\]
are fully faithful. Thus every morphism of restricted sheaves, in the
indicated target category, extends uniquely to a morphism of the
logarithmic sheaves.
\end{theorem}

\subsection{Cohomological inputs and local schemes}\label{sec:inputs}

Throughout the proof \(F\) is a sheaf in
\(\mathrm{CI}_{\dNis}^{\log}(k,\Lambda)\) and \(H=\omega_\sharp F\).
We write \(C_F(Y)=\R\Gamma_{\dNis}(Y,F)\), the complex of derived
sections, and \(C_H(V)=\R\Gamma_{\Nis}(V,H)\) for an ordinary scheme
\(V\). Complexes use cohomological indexing, so
\(H^i(C[1])=H^{i+1}(C)\). The cofiber of a map of complexes
is its mapping cone and its fiber is the cofiber shifted by \([-1]\).
A square of complexes is homotopy cartesian if its upper left term is
the fiber of the difference of the two maps from adjacent terms into
the lower right term.

We use four established facts.
\begin{enumerate}
\item \emph{Purity.} If \(V\subseteq Y\) is a dense open with the
  pulled back log structure and \(Y\) has smooth underlying scheme,
  \(F(Y)\to F(V)\) is injective \cite[Theorem~5.10]{BM}.
\item \emph{Truncation.} A complex of dividing Nisnevich sheaves is
  called interval-local when its derived sections are invariant under
  product with \(\bx\). Its usual derived truncations and its
  individual cohomology sheaves are again interval-local
  \cite[Theorem~5.7]{BM}. One can describe these complexes as
  presheaves of complexes satisfying cohomological descent for the
  site and interval invariance \cite[\S2.3]{BM}.
\item \emph{Blow-ups at the boundary.} Blow up a smooth center contained
  in the boundary and in normal crossings position with it. The new
  boundary is the reduced total transform of the old boundary. The
  resulting map induces an isomorphism on \(C_F\). This is
  \cite[Proposition~7.0.2 and Theorem~7.2.10]{BPO} applied in the
  effective category without transfers. Here normal crossings position
  means that the ideal of the center and the boundary components are
  given locally by subsets of a coordinate system.
\item \emph{Adding a smooth divisor.} Let \(J\) be a smooth divisor
  meeting the existing boundary of \(Y\) transversely, with no common
  component. Give \(J\) the boundary induced by \(Y\), and let \(Y^+\)
  be \(Y\) with \(J\) added to the boundary. If the normal line
  bundle of \(J\) in \(Y\) has a chosen trivialization, the
  logarithmic Gysin triangle identifies
  \[
     \cofib\bigl(C_F(Y)\longrightarrow C_F(Y^+)\bigr)
      \;\simeq\;
     \fib\bigl(C_F(J\times\Pp^1)\xrightarrow{\infty^*}C_F(J)\bigr)[1].
     \tag{*}
  \]
  The projective-line factor on the right has trivial log structure.
  In the motivic formulation, the cofiber of the map from \(Y^+\)
  to \(Y\) is the Thom motive of the normal line. The Thom motive
  for a trivial line on \(J\) is the cofiber of the infinity section
  \(J\to J\times\Pp^1\). Taking derived sections with coefficients
  in \(F\) gives \((*)\). Use \cite[Propositions~7.0.2 and 7.4.5
  and Theorem~7.5.4]{BPO}. The equivalence is compatible with
  restricting \(Y\) to an open subscheme and restricting the
  chosen trivialization there. See also
  \cite[Proposition~3.2.13, Theorem~3.2.19 and Remark~3.2.20]{BPOh}.
\end{enumerate}

A \emph{strict Nisnevich excision square} consists of an open inclusion
and a strict \'etale map inducing an isomorphism over its reduced
closed complement, together with their intersection. Derived sections
send such a square to a homotopy cartesian square. This remains true for
the limit squares below. We evaluate on function fields and henselian
local schemes by filtered colimits over smooth models and pointed
\'etale neighborhoods. Divisors, morphisms of finite presentation,
and the excision squares used here can be represented after shrinking
those models. We may apply derived excision directly on smooth models and then
take filtered colimits, which are exact on modules and commute
with finite limits. For the field value \(C_H(\Spec K)\) this
gives \(H(K)[0]\). In fact a positive-degree class on a smooth
model vanishes after a Nisnevich cover locally. A member of the
cover over the generic point with residue field \(K\) admits a
section after shrinking to a dense open, so the class becomes
zero in the colimit. The passage to the limit is the same in
strict Nisnevich excision, where one gives each neighborhood
the pulled back boundary. See \cite[Lemmas~2.6 and 4.3 and
the proof of Claim~4.5]{BM} for these excision diagrams. For brevity, we sometimes
write \(F(K)\) for \(H(K)\).

\begin{lemma}[Coordinates at a point]\label{lem:coordinates}
Let \(x\) be any point of a smooth \(k\)-scheme \(X\), with residue
field \(K\). Set \(n=\dim\mathcal O_{X,x}\). A regular system of
parameters \(t_1,\ldots,t_n\) can be used to identify the
henselization \(\mathcal O_{X,x}^h\) with
\(K\{t_1,\ldots,t_n\}\), the henselization of the local ring of
\(\A^n_K\) at the origin. If \(X\) has a strict normal crossings
boundary, we may use local equations of its components through
\(x\) as the first \(r\) parameters, where \(r\) is their number.
\end{lemma}

\begin{proof}
The field \(K\) is finitely generated over \(k\). Choose a separating
transcendence basis \(b_1,\ldots,b_m\) of \(K/k\), so that \(K\) is
finite separable over \(K_0=k(b_1,\ldots,b_m)\). This is possible
because \(k\) is perfect. Lift the \(b_i\) to regular functions near
\(x\). Together with the \(t_j\) they define a map to
\(\A^{m+n}_k\) that is \'etale at \(x\). Indeed, the
\(dt_j\) give the conormal directions at \(x\), and the \(db_i\)
give a basis of the differentials of \(K/k\).
Pass to the generic fiber over \(\A^m_k\), giving an \'etale
neighborhood over \(\A^n_{K_0}\) of its origin. The local ring at
\(x\) is unchanged since the nonzero rational functions of the
first \(m\) variables are units there. In the base change to \(K\),
the point picked by the identity embedding of the residue field
has residue \(K\). Its henselization agrees with the henselization
at \(x\), since this base change is \'etale with unchanged residue
field at the selected point. It also agrees with the henselization
of \(\A^n_K\) at the origin by the same property. The normal
crossings hypothesis permits the specified choice of the \(t_j\).
Cf.\ \cite[Lemma~4.3 and its proof]{BM}.
\end{proof}

\subsection{Sections on the marked line and at a boundary divisor}
\label{sec:line}\label{sec:heightone}

For a base \(S\) and \(I\subseteq\{0,\infty\}\), write
\(P_{I,S}=S\times(\Pp^1,I)\). If \(S\) is already logarithmic,
this notation keeps its boundary on all four products. Write
\(L_S=P_{\{0,\infty\},S}\). Its interior over an ordinary \(S\)
is \({\Gm}_S\).

\begin{lemma}[The two affine charts]\label{lem:charts}
The square
\[
\begin{array}{ccc}
 C_F(P_{\varnothing,S}) &\longrightarrow& C_F(P_{\{0\},S})\\
 \big\downarrow && \big\downarrow\\
 C_F(P_{\{\infty\},S}) &\longrightarrow& C_F(L_S)
\end{array}
\]
is homotopy cartesian. In particular, \(F(P_{\varnothing,S})\simeq
F(S)\) via projection. For an ordinary \(S\), there is an exact
sequence
\[
 0\longrightarrow H(S)\longrightarrow F(L_S)
    \xrightarrow{\delta_S} H^1_{\Nis}(\Pp^1_S,H).
\]
The map \(\delta_S\) is the ordinary Mayer--Vietoris connecting
map applied to the restriction to \({\Gm}_S\). It is surjective
when \(S=\Spec K\) for a finitely generated field \(K/k\).
\end{lemma}

\begin{proof}
On the affine chart at zero the maps in one direction of the square
are identities, because the infinity mark is absent. On the affine
chart at infinity the maps in the other direction are identities.
On their intersection both directions are identities. Zariski
descent proves the assertion for the square, since finite limits
of complexes commute with each other. Both singly marked lines
are copies of the interval over \(S\). The difference of their
restriction maps on \(H^0\) takes \((a,b)\) to the pullback of
\(a-b\) from \(S\). Pullback \(F(S)\to F(L_S)\) has a left
inverse given by the section at \(1\). The long exact
cohomology sequence of the square gives the claimed exactness
and the identification of \(F(P_{\varnothing,S})\).
Restriction of the square to the ordinary affine open cover
of \(\Pp^1_S\) gives a map to the ordinary Mayer--Vietoris
square. The map on the upper left term is the cohomology
identification for the trivial log structure. Naturality
of connecting maps gives the formula for \(\delta_S\).
When \(S=\Spec K\), the two singly marked lines have vanishing
first cohomology by strict interval invariance and the field
case of Section~\ref{sec:inputs}.
\end{proof}

\begin{lemma}[Action of a square scaling]\label{lem:scaling}
Let \(S\) be ordinary and smooth over \(k\) and \(a\) a unit on \(S\).
The automorphism \(t\mapsto a^2 t\) of the projective line
with trivial log structure over \(S\) induces the identity on
\(C_H(\Pp^1_S)\) in the derived category.
\end{lemma}

\begin{proof}
On the smooth pairs over \(k\), \(C_F\) is a strict Nisnevich
sheaf of spectra by taking the Eilenberg--Mac Lane spectra of its
derived section complexes. It satisfies dividing descent and
interval invariance. It therefore belongs to \(\mathcal L_k\).
Lemma~\ref{lem:unmarked-scaling} over \(S\), followed by the
smooth left adjoint from \(S\) to \(k\), gives the claimed
action on the represented projective line over \(S\).
Evaluating the derived sections gives the assertion.
\end{proof}

\begin{proposition}[An ordinary identity for extension to \(L_K\)]
\label{prop:criterion}
Let \(K/k\) be a finitely generated field extension and write
\(L_K=L_{\Spec K}\). The restriction \(F(L_K)\to H({\Gm}_K)\)
identifies \(F(L_K)\) with the subgroup of sections satisfying
\begin{equation}\label{eq:criterion}
 f(s^2t)-f(t)=f(s^2)-f(1)
          \qquad\text{in }H({\Gm}_{s,K}\times{\Gm}_{t,K}).
\end{equation}
For example, \(f(s^2t)\) denotes pullback along
\((s,t)\mapsto s^2t\). Thus membership in this subgroup is
preserved by every morphism of ordinary Nisnevich sheaves.
\end{proposition}

\begin{proof}
The restriction is injective by purity. Take
\(\ell\in F(L_K)\), put \(S={\Gm}_{s,K}\), and pull back \(\ell\) to
\(L_S\). The two pullbacks along \((s,t)\mapsto(s,s^2t)\)
and \((s,t)\mapsto(s,t)\) have the same image under \(\delta_S\).
Indeed, the scaling fixes the two marked divisors, so the
connecting map is natural for it. Lemma~\ref{lem:scaling} gives
its trivial action on \(H^1_{\Nis}(\Pp^1_S,H)\). By
Lemma~\ref{lem:charts}, their difference comes from \(H(S)\).
Evaluation at \(t=1\) identifies that difference with the
right hand side of \eqref{eq:criterion}. The calculation
over \(K\) here and below can be carried out over smooth
models and passed to the colimit described in
Section~\ref{sec:inputs}.

Conversely, let \(f\) satisfy \eqref{eq:criterion}. The surjectivity of
\(\delta_K\) and ordinary Mayer--Vietoris give a decomposition
\[
 f(t)=\ell(t)+b(t)+d(t^{-1}),
 \qquad \ell\in F(L_K),\quad b,d\in H(\A^1_K).
\]
Absorb \(b(0)+d(0)\) into \(\ell\) so that \(b(0)=d(0)=0\).
The difference \(f-\ell\) satisfies the same identity. On
\({\Gm}_{s,K}\times{\Gm}_{t,K}\) it follows that
\[
 \bigl(b(s^2t)-b(t)\bigr)
 +\bigl(d(s^{-2}t^{-1})-d(t^{-1})\bigr)
\]
is pulled back from \({\Gm}_{s,K}\). Each summand extends to one
affine chart of the \(t\)-line. By subtracting the other summand
from the section pulled back from the base, it also extends to
the other chart. The sections that extend to both charts
are precisely the image of \(H(\Pp^1_S)\) by the sheaf axiom.
By Lemma~\ref{lem:charts} these are the pullbacks from \(H(S)\).
Evaluating at \(t=0\) in the first chart and at \(t^{-1}=0\) in
the second shows that both summands vanish. The equality
\(b(s^2t)=b(t)\) extends to all of \(\A^2_{s,t}\) by purity.
Restriction to the line \(s=0\) gives \(b=0\).
For \(d\) use \(s^{-1}\) and \(t^{-1}\) as the coordinates of
another affine plane to get \(d=0\) by the same argument.
Thus \(f\) is the restriction of \(\ell\).
\end{proof}

\begin{proposition}\label{prop:local}
Let \(R=K\{t\}\), let \(E=\Frac R\), and let \((R,0)\) denote
\(\Spec R\) with its closed point marked. Inside \(H(E)\) one has
\begin{equation}\label{eq:local}
 F(R,0)=\im\bigl(H(R)\to H(E)\bigr)
         +\im\bigl(F(L_K)\to H(E)\bigr).
\end{equation}
The second map is restriction first to \({\Gm}_K\) and then along
the punctured henselian neighborhood of \(0\).
\end{proposition}

\begin{proof}
Put \(P=P_{\{\infty\},\Spec K}\) and \(L=L_K\).
For \(P\), use the open complement of \(0\) and the henselian
neighborhood \(R\) with trivial log structure. For \(L\), use the same
open complement with the same induced log structure, and use
\((R,0)\) as the henselian neighborhood. The intersection in
both excision squares is \(\Spec E\) with trivial log structure.
Taking the cofiber of the map between the two squares of derived
sections yields the isomorphism
\[
 \cofib\bigl(C_F(P)\to C_F(L)\bigr)
   \xrightarrow{\sim}
 \cofib\bigl(C_H(\Spec R)\to C_F(R,0)\bigr).
\]
For instance, the square whose rows are
\(C_F(P)\to C_F(L)\) and \(C_H(\Spec R)\to C_F(R,0)\)
is homotopy cartesian by the two excision squares. The
isomorphism of horizontal cofibers follows. In constructing the two squares, use pointed \'etale neighborhoods
in the chart at \(0\). Removing the other points over \(0\) ensures
that the fiber is just the selected point with residue field \(K\).
The complement of \(0\) glues with each such neighborhood to
give a strict excision square. Its intersection with the
neighborhood gives \(\Spec E\) in the filtered limit.

Since \(H^1(C_F(P))=0\), the left cofiber receives a surjection
on \(H^0\) from \(F(L)\). The right cofiber has the image of
\(H(R)\) as the kernel of the map into its \(H^0\) from
\(F(R,0)\). The comparison gives \eqref{eq:local} inside
\(F(R,0)\), which injects into \(H(E)\) by purity and passage
to the filtered colimit.
\end{proof}

For two strictly interval-invariant sheaves \(F,G\), let
\(\varphi:\omega_\sharp F\to\omega_\sharp G\) be an ordinary
Nisnevich sheaf morphism. It preserves the first summand of
\eqref{eq:local} by naturality. It preserves the second summand
by Proposition~\ref{prop:criterion}. Thus it preserves the
images of sections at henselian local schemes along boundary
divisors, using Lemma~\ref{lem:coordinates} to reduce a
codimension-one henselian local ring to this form.

\subsection{Extension of sections and full faithfulness}
\label{sec:hartogs}\label{sec:finish}

\begin{lemma}\label{lem:reduced}
For a log scheme \(Y\) in the smooth basis, define
\[
 B_F(Y)=
 \ker\!\left(H^1_{\dNis}(Y\times\Pp^1,F)
        \xrightarrow{\infty^*} H^1_{\dNis}(Y,F)\right),
\]
using the trivial log structure on the projective-line factor.
Then \(B_F\) is a strictly interval-invariant sheaf. Its restrictions
to dense strict opens are injective.
\end{lemma}

\begin{proof}
On the full log site take a complex of presheaves \(I\) computing
the derived sections of \(F\). Choose it fibrant in the dividing
Nisnevich local model structure, so \(I(Y)\) computes \(C_F(Y)\)
for every \(Y\). Binda and Merici show that
\(Y\mapsto I(Y\times V)\) also computes derived sections of a
complex of dividing Nisnevich sheaves for any log smooth \(V\)
\cite[Lemma~2.10]{BM}.
This is the derived internal Hom out of the free sheaf on \(V\),
whose value on \(Y\) is given by product with \(V\).
Apply this with the trivial-log projective
line. Taking the derived fiber of restriction along \(\infty\)
and shifting gives a complex \(J\) whose derived sections are
\[
 \R\Gamma_{\dNis}(Y,J)
   \simeq \fib\bigl(C_F(Y\times\Pp^1)\to C_F(Y)\bigr)[1].
\]
This is interval-local by the invariance of \(C_F\) and commutation
of the two products. In the displayed fiber, the map
onto \(C_F(Y)\) has a section induced by projection. Hence
its \(H^0\) is the kernel on sections, which vanishes by
Lemma~\ref{lem:charts}. Its negative cohomology also vanishes.
The shifted complex of derived sections has zero cohomology in
negative degrees on every object. Sheafifying its cohomology
presheaves gives the cohomology sheaves of \(J\), so \(J\) is
cohomologically nonnegative. For such a complex the degree-zero
derived sections equal the sections of its zeroth cohomology
sheaf. This follows by applying derived sections to the
truncation triangle in degrees \(0\) and \(\geq1\).
Thus \(B_F\) is the zeroth cohomology sheaf of \(J\). The
truncation theorem makes \(B_F\) strictly interval-invariant.
Purity for \(B_F\) gives the claimed injectivity. These
properties also hold on the local schemes evaluated by
filtered colimit.
\end{proof}

\begin{lemma}\label{lem:puncture}
Let \(T\) be the spectrum of a henselization of a local ring
on a smooth \(k\)-scheme, with the log structure of a strict
normal crossings boundary. Let \(x\) be its closed point.
If \(\dim T=n\geq2\), restriction induces
\(F(T)\simeq F(T-\{x\})\).
\end{lemma}

\begin{proof}
Let \(r\) be the number of boundary components through \(x\).
We prove the statement by descending induction on \(r\) for
fixed \(n\), with all choices of the smooth model allowed.

If \(r=n\), use Lemma~\ref{lem:coordinates} to view \(T\) as
the henselian neighborhood of \(\infty^n\) in
\(\bx_K^n\), where \(K\) is the residue field. Let
\(W=\bx_K^n-\{\infty^n\}\). A section on \(W\) restricts
to the two opens where, respectively, the first and the
second coordinate avoid infinity. By interval invariance the
restrictions are pulled back from the first and second
\(\A^1_K\), respectively. Their restrictions to the common
open agree. Restrict further to \(\A^2_K\) by setting the
remaining coordinates to zero. The resulting equality is
\(p_1^*b=p_2^*d\), where \(p_i\) are the two projections and
\(b,d\in H(\A^1_K)\). Setting each coordinate to zero in turn
shows that \(b,d\) come from the same section of \(H(K)\).
One of the two opens is dense in \(W\), so purity shows that
the original section equals the pullback of that constant
section on \(W\). Consequently \(F(\bx_K^n)\to F(W)\) is surjective.
For the passage to \(T\), use excision on derived sections.
With \(C=C_F\), the map
\[
 \cofib\bigl(C(\bx_K^n)\to C(W)\bigr)
   \longrightarrow \cofib\bigl(C(T)\to C(T-\{x\})\bigr)
\]
is an isomorphism. The \(H^0\) of the left side is zero,
because \(F(\bx_K^n)\to F(W)\) is surjective and
\(H^1(C(\bx_K^n))=H^1(C(\Spec K))=0\) by interval
invariance. Thus \(F(T)\to F(T-\{x\})\) is surjective by
the long exact sequence of the right side. Purity gives
injectivity.

If \(r<n\), choose the parameters as in
Lemma~\ref{lem:coordinates} and add the divisor \(J\)
defined by \(t_{r+1}=0\), obtaining \(T^+\). Give \(J\)
the boundary induced by the original \(r\) components.
Its normal line is trivialized by \(t_{r+1}\). The
exact sequence from \((*)\) in Section~\ref{sec:inputs}
gives
\[
 F(T)\longrightarrow F(T^+)\longrightarrow B_F(J)
\]
and the corresponding sequence after removal of \(x\).
The maps commute with restriction, with the same
trivialization restricted to the smaller normal line.
A section of \(F(T-\{x\})\) extends to \(T^+\) by induction,
applied to its restriction to \(T^+-\{x\}\). The
obstruction to this extension coming from \(F(T)\) is its
image in \(B_F(J)\). The image vanishes on \(J-\{x\}\)
by exactness on the punctured pair. That open is dense
in \(J\) since \(\dim J=n-1\geq1\).
Lemma~\ref{lem:reduced} gives vanishing on \(J\), and
exactness gives a lift to \(F(T)\). Its restriction to
\(T-\{x\}\) is the prescribed section, since the map
from \(F(T-\{x\})\) to \(F(T^+-\{x\})\) is injective
by purity on the common interior. This completes the induction.
\end{proof}

\begin{proposition}\label{prop:codim}
Let \(Y\) have smooth underlying scheme and strict normal crossings
boundary. If a closed subset \(Z\) has codimension at least two,
then \(F(Y)\to F(Y-Z)\) is an isomorphism.
\end{proposition}

\begin{proof}
Purity gives injectivity. For surjectivity, fix a section on \(Y-Z\).
At a generic point \(x\) of a component of \(Z\), the punctured local
scheme is contained in \(Y-Z\). Pull back the section to the
punctured henselian neighborhood of \(x\). Lemma~\ref{lem:puncture}
gives a section on the henselization. The henselian section
comes from a strict \'etale neighborhood with a point over \(x\)
having the same residue field. The \'etale map becomes an
isomorphism on an open neighborhood of \(x\) in the reduced
component of \(Z\), after restriction of the source. Indeed
the selected point over the generic point gives a rational
section of the boundary restriction of the \'etale map.
It spreads to a section over an open of that component. Its
image is open in the restriction of the \'etale neighborhood
to the component, because a section of an \'etale morphism is
an open immersion. Shrink the base away from the other
components near \(x\), then remove the closed complement of
the image in the preimage of the reduced closed set. After
restriction, the \'etale map induces an isomorphism over that
closed set. On the punctured overlap the two sections agree
after shrinking by the colimit description. One can also
take the source neighborhood irreducible and use purity
there once the sections agree on a dense open of the overlap. The strict Nisnevich square then
glues the sections to one defined near \(x\).
Equivalently, one can first use the Nisnevich square for
the local ring at \(x\) and its henselization and then
spread the section back to a neighborhood in \(Y\).

Perform this extension at each generic point of \(Z\).
Purity makes the sections agree on all overlaps because their
restrictions to the original dense open agree. The remaining
closed subset has strictly smaller dimension. Repeating
the procedure gives a section on all of \(Y\).
\end{proof}

\begin{proof}[Proof of Theorem~\ref{thm:restriction-sheaves}]
First consider \(\varphi:\omega_\sharp F\to\omega_\sharp G\)
without transfers. Given \(a\in F(X)\) on an object of the smooth
basis, apply \(\varphi\) to its restriction on \(X^\circ\).
At the generic point of any component of the boundary, the
preservation following Proposition~\ref{prop:local} gives an
extension of the resulting section to the henselian local
scheme with its closed point marked. Strict Nisnevich
excision descends the extension to a neighborhood of that
generic point, exactly as in the gluing argument in
Proposition~\ref{prop:codim}. By purity the extensions agree
where they are both defined. Their complement in \(X\) has
codimension at least two, so Proposition~\ref{prop:codim}
for \(G\) extends the section to \(X\).

The section on \(X\) is unique by the injection
\(G(X)\to(\omega_\sharp G)(X^\circ)\). Uniqueness gives
\(\Lambda\)-linearity and compatibility with all log
morphisms because each induces a map of interiors. It also
shows that restriction on morphisms is injective. The
smooth basis computes the full category of dividing
Nisnevich sheaves, so the constructed map is a morphism
on the full site.

For the version with transfers, apply the result to the
underlying sheaves. These are still strictly interval-invariant
by the cohomology comparison in Section~\ref{sec:sheaf-restriction}.
For a finite log correspondence, its action restricts to the
action of the corresponding ordinary finite correspondence
between interiors, by the definition of \(\lCor(k)\).
Thus the two possible composites with the lifted map
agree on the source interior. Purity gives equality on the
source on the smooth basis. Compatibility on the full
site follows either by the equivalence with the smooth
basis for sheaves with log transfers
\cite[Proposition~4.7.5]{BPO}, or by testing on a dividing
cover by smooth pairs and then on their interiors.
Hence the lifted map has log transfers and is unique.
\end{proof}

\begin{remark}[Scope of the proof]
The scaling identity takes place over \({\Gm}_s\), with \(s\)
a variable. In characteristic two, the calculation of the base
locus, the matrix identity and the specialization \(s=0\)
remain valid. The henselian coordinate lemma works over finite
fields because the residue field is separably generated over \(k\).
The two blow-ups in Lemma~\ref{lem:unmarked-scaling} have explicit smooth
centers, so smooth-center blow-up invariance applies over every
perfect field. The truncation theorem is used for the single
cohomology sheaf in Lemma~\ref{lem:reduced}.
\end{remark}

\section{A direct proof of sheaf rigidity over a perfect field}
\label{sec:perfect}
Throughout this section \(k\) is perfect. For a log \'etale sheaf \(F\)
with log transfers, denote by \(H\) the restriction
\(H(U)=F(U,\triv)\) to smooth schemes with trivial log structure.
It is an ordinary \'etale sheaf with finite transfers. For a function
field \(K/k\), by which we mean a finitely generated field extension,
sections over \(K\) and over the schemes over \(K\) occurring below
mean the filtered colimit of sections over smooth \(k\)-models.
The same convention is used for Nisnevich cohomology. Maps, finite
covers, and the finite open covers used below descend after shrinking
a smooth integral model of \(K\). A Nisnevich cover of such a model
has a lift of the generic point with residue field \(K\), which gives
a section over a dense open. Positive-degree Nisnevich cohomology
classes vanish locally for the Nisnevich topology, so their images
vanish after such a shrinking. In particular, the filtered
Nisnevich cohomology of \(\Spec K\) is concentrated in degree zero.
Filtered colimits of abelian groups are exact.

Let \(\epsilon_*\) denote restriction from the log \'etale topology
to the dividing Nisnevich topology. The \emph{homotopy heart} consists
of strictly \(\bx\)-invariant dividing Nisnevich sheaves with log
transfers \cite[Theorem~5.7]{BM}. It is an abelian subcategory in
which kernels and cokernels agree with those in the category of
dividing Nisnevich sheaves with log transfers. When \(\epsilon_*F\)
belongs to the heart, we also write \(F\) for this restriction.
The injectivity theorem (called purity)
\cite[Theorem~5.10]{BM} gives injectivity of \(F(X)\to F(V)\) for a
dense open log subscheme \(V\subset X\) with \(X\) smooth as a scheme.
In particular, \(H(U)\to H(k(U))\) is injective for \(U\) smooth
integral.

\subsection{The marked projective line and sections with prescribed poles}
\label{sec:twopoints}

For the boundary calculations, \(F\) is a log \'etale sheaf whose
dividing Nisnevich restriction is strictly \(\bx\)-invariant, and
\(nF=0\) for some integer \(n>0\) invertible in \(k\).
For a function field \(K/k\), the map \(H(K)\to H(K')\) is
injective for any finite separable extension \(K'/K\). Indeed, the
extension spreads out to a finite \'etale surjection of smooth
models after shrinking, and \(H\) is an \'etale sheaf.

\begin{lemma}\label{lem:two}
For every function field \(K/k\), pullback along the projection
induces an isomorphism
\[
 H(K)\xrightarrow{\ \sim\ }F(\Pp^1_K,\{0,\infty\}).
\]
\end{lemma}

\begin{proof}
Take \(v\in F(\Pp^1_K,\{0,\infty\})\) and write \(v(a)\in H(K)\)
for its evaluation at \(a\in K^\times\). We first prove that
\begin{equation}\label{eq:multiplicative}
 v(a)+v(b)=v(ab)+v(1)\qquad (a,b\in K^\times).
\end{equation}
Put \(c=ab\), \(A=a+b\), and \(B=c+1-A\). If \(B=0\), the
two monic quadratics with roots \(a,b\) and \(c,1\) coincide,
so their cycles of roots are equal.

Suppose \(B\ne0\). The rational function
\[
 q(y)=\frac{y^2-Ay+c}{By}
\]
defines a finite morphism \(\Pp^1_y\to\Pp^1_t\) of degree two.
Its poles at \(y=0\) and \(y=\infty\) are simple since \(B,c\ne0\). Thus
\[
 q^*\{\infty\}=\{0\}+\{\infty\}.
\]
Its graph, viewed as a correspondence in the reverse projection
from the \(t\)-line to the \(y\)-line, has normalization \(\Pp^1_y\).
The pulled-back source log structure is exactly the compactifying
log structure at \(0,\infty\). This verifies admissibility in the
literal sense of \cite[Definition~2.1.1]{BPO}. The finite-chart
equation is
\[
 y^2-(A+Bt)y+c=0.
\]
This is finite flat over the \(t\)-chart, and the same holds at
infinity by the two simple poles. The construction is valid in
characteristic two as well. It spreads to a smooth model of \(K\)
after making \(a,b,B,c\) regular and \(B,c\) invertible.

The fiber cycles at \(t=0\) and \(t=1\) are \([a]+[b]\) and
\([c]+[1]\), respectively, with multiplicities if roots coincide.
The two evaluations of a section of \(F(\bx_K)=H(K)\) agree.
Apply this fact to the transfer of \(v\) along the correspondence
to obtain \eqref{eq:multiplicative}.

Thus \(a\mapsto v(a)-v(1)\) is a group homomorphism on \(K^\times\).
For each \(a\), its value becomes zero after adjoining an \(n\)-th
root of \(a\), since the homomorphism then takes \(a\) to \(n\)
times the value on a root. This extension can be chosen finite
separable, and \(H\) is \'etale separated. Hence \(v(a)=v(1)\).
The same argument over \(K(y)\), with \(a=y\), gives the equality
at the generic point of the projective line. Purity on smooth
models now gives \(v=p^*v(1)\), for \(p:\Pp^1_K\to\Spec K\).
\end{proof}

\begin{lemma}\label{lem:p1}
For every function field \(K/k\),
\[
 H(\Pp^1_K)=H(K),\qquad H^1_{\Nis}(\Pp^1_K,H)=0.
\]
\end{lemma}

\begin{proof}
Apply Lemma~\ref{lem:charts} to the dividing Nisnevich
restriction of \(F\), writing \(L_K=(\Pp^1_K,\{0,\infty\})\).
It gives the equality on sections and identifies
\(H^1_{\Nis}(\Pp^1_K,H)\) with the quotient of \(F(L_K)\)
by the constants \(H(K)\). The comparison to the dividing
Nisnevich site and the calculation on function fields were
justified above. Lemma~\ref{lem:two} makes the quotient zero.
\end{proof}

\begin{lemma}\label{lem:split}
Let \(K/k\) be a function field and \(S\subset \Pp^1(K)\) a finite
nonempty set. Write \(U_S=\Pp^1_K-S\). The sum of the restrictions
induces an isomorphism
\[
 \bigoplus_{a\in S}
   H(\Pp^1_K-\{a\})/H(K)
      \ \xrightarrow{\ \sim\ }\ H(U_S)/H(K).
\]
This isomorphism is natural for automorphisms of \(\Pp^1_K\)
preserving \(S\).
\end{lemma}

\begin{proof}
Constants are embedded by pullback because each of the opens
is dense and restriction from \(\Pp^1_K\) is injective.
For disjoint nonempty subsets \(S_1,S_2\) with union \(S\),
\(U_{S_1}\) and \(U_{S_2}\) cover \(\Pp^1_K\), and their
intersection is \(U_S\). Lemma~\ref{lem:p1} and the
Mayer--Vietoris sequence show that the sum map
\(H(U_{S_1})\oplus H(U_{S_2})\to H(U_S)\) is surjective and has
kernel the antidiagonal copy of \(H(K)\). Here using the sum map
instead of the difference map changes the sign on the second
summand. Modulo constants this gives a direct sum. Induction gives
the claimed isomorphism, which is canonical because the entire
map is the sum of the given restrictions.
\end{proof}

\begin{proposition}\label{prop:evaluation}
For every function field \(K/k\), \(s\in H(\A^1_K)\), and \(a\in K\),
one has \(s(a)=s(0)\) in \(H(K)\).
\end{proposition}

\begin{proof}
The case \(n=1\) is immediate. After the finite separable extension
\(K(\mu_n)/K\) we may assume that \(K\) contains all \(n\)-th roots of
unity. Equality descends since \(H\) is \'etale separated.
Use coordinates \(z\) and \(w\) on two copies of \(\Pp^1_K\), and
let \(f:\Pp^1_K\to\Pp^1_K\) be \(w=z^n\). Put
\[
 G=\mu_n(K),\qquad S=\mu_n(K)\subset \Pp^1_z(K),\qquad
 V=\Pp^1_w-\{1\}.
\]
Then \(U_S=f^{-1}(V)\). The group \(G\) acts on the source by
\(z\mapsto\zeta z\), freely and transitively on \(S\). The points
\(0,\infty\) are fixed and both lie in \(U_S\). The coordinate
\(x(w)=a/(1-w)\) identifies \(V\) with \(\A^1_K\) when \(a\ne0\),
taking \(0\) to \(a\) and \(\infty\) to \(0\).
For \(a=0\) the desired equality is tautological.
The section \(b=f^*x^*s\in H(U_S)\) is \(G\)-invariant, since
\(f\circ g=f\) for all \(g\in G\), and satisfies
\(b(0)-b(\infty)=s(a)-s(0)\).

By Lemma~\ref{lem:split}, the invariant class of \(b\) in
\(H(U_S)/H(K)\) is the sum of the translates of its component
at \(1\), since \(G\) acts freely and transitively on \(S\).
Choose a representative \(c\in H(\Pp^1_K-\{1\})\) of that
component. Then
\[
 b-\sum_{g\in G}g^*(c|_{U_S})\in H(K).
\]
The two points \(0,\infty\) are fixed by \(G\) and lie in \(U_S\).
Their difference of evaluations kills constants and takes every
summand to \(c(0)-c(\infty)\). Hence
\[
 s(a)-s(0)=n\bigl(c(0)-c(\infty)\bigr)=0 .
\]
\end{proof}

\begin{corollary}\label{cor:affine}
The ordinary restriction \(H\) is \(\A^1\)-invariant.
\end{corollary}
\begin{proof}
For \(s\in H(\A^1_K)\), apply Proposition~\ref{prop:evaluation}
over \(K(T)\) at \(a=T\), with \(T\) transcendental. Purity
identifies this equality at the generic point with
\(s=p^*s(0)\) on \(\A^1_K\), where \(p\) is the projection.
For a smooth integral \(U/k\) and \(s\in H(U\times\A^1)\),
the difference \(s-p^*i_0^*s\) restricts to zero on
\(\A^1_{k(U)}\). It therefore vanishes on \(U'\times\A^1\)
for a dense open \(U'\subset U\), and vanishes everywhere by
purity. The argument on connected components proves the result
for every smooth \(U/k\).
\end{proof}

\subsection{Reconstruction and elementwise torsion}
\label{sec:perfectresult}

\begin{proof}[Proof of Theorem~\ref{thm:sheaf} over a perfect field]
Let \(F\) satisfy the hypotheses of the theorem. With integral
coefficients, the derived restriction with log transfers
\cite[\S5.3]{BPO} satisfies
\[
 \R\Gamma_{\dNis}(X,\R\epsilon_*F)
   \simeq \R\Gamma_{\letop}(X,F).
\]
Thus \(\R\epsilon_*F\) is strictly \(\bx\)-invariant in the
dividing Nisnevich topology. Binda and Merici show that the
usual cohomological truncations preserve such complexes
\cite[Theorem~5.7]{BM}. It follows that
\(H^0(\R\epsilon_*F)=\epsilon_*F\) belongs to the homotopy heart.

For the bounded torsion argument, let \(E\) be a log \'etale sheaf
with log transfers such that \(\epsilon_*E\) belongs to the heart
and \(nE=0\) for some positive integer \(n\) invertible in \(k\).
Corollary~\ref{cor:affine} gives \(\A^1\)-invariance of the
ordinary restriction \(U\mapsto E(U,\triv)\). By classical Suslin
rigidity \cite[Theorem~7.20]{MVW}, this ordinary restriction is the
\'etale pullback of \(\eta^*E\). Set \(E_0=\eta_\sharp\eta^*E\).
It has its canonical log transfers and is strictly \(\bx\)-invariant
by \cite[Propositions~8.5.4 and 8.5.6]{BPO}, so the truncation
argument also puts \(\epsilon_*E_0\) in the heart. The counit
\(E_0\to E\) restricts to the classical isomorphism on ordinary
smooth schemes.
Restriction to objects with trivial log structure is exact as
a functor to ordinary Nisnevich sheaves with transfers, since
dividing Nisnevich covers on these objects are ordinary Nisnevich covers.
By purity, a sheaf in this heart that vanishes on these objects
vanishes on each log-smooth object with smooth underlying scheme.
Such objects detect all sheaves by \cite[Proposition~4.7.5]{BPO}.
Thus the kernel and cokernel of \(E_0\to E\) vanish, and \(E_0\to E\)
is an isomorphism of log \'etale sheaves with transfers.

For the original sheaf \(F\) and each integer \(n>0\) invertible in \(k\),
let \(F[n]=\ker(n\colon F\to F)\) in log \'etale sheaves with
transfers. Restriction to the dividing Nisnevich topology preserves
kernels, which are taken on sections. Since \(\epsilon_*F\) belongs
to the homotopy heart, so does \(\epsilon_*F[n]\), as the kernel of
multiplication by \(n\) in that heart. The bounded torsion argument
applies to \(F[n]\). These isomorphisms
are canonical and hence compatible with the inclusions for
divisibility among positive integers invertible in \(k\).

By the torsion hypothesis, \(F=\varinjlim_n F[n]\) as a log
\'etale sheaf with transfers. Indeed, the union is all of \(F\)
after sheafification. Likewise
\(\eta^*F=\varinjlim_n\eta^*(F[n])\) on the small site, since
restriction to the small site preserves colimits and kernels
\cite[proof of Proposition~8.5.2]{BPO}. Pullback
\(\eta_\sharp\) preserves colimits, so taking the colimit of the
counits on \(F[n]\) gives the required isomorphism.
Naturality of the counit gives compatibility with the coefficient action.
\end{proof}

\section{Passage to a perfect field}\label{sec:perfection}

We give the base-change argument in characteristic \(p>0\).
The perfection \(k^{\mathrm{perf}}\) is the purely inseparable
extension formed by adjoining all \(p^r\)-th roots of elements
of \(k\) for \(r\geq 0\). We take correspondences with coefficients
in \(R=\Z[1/p]\), since any abelian group with torsion of order
prime to \(p\) carries a canonical \(R\)-module structure.

\begin{lemma}\label{lem:perfection}
Let \(L/k\) be a purely inseparable algebraic extension in
characteristic \(p>0\) and put \(R=\Z[1/p]\).
Base change on finite log correspondences with \(R\)-coefficients
is fully faithful. For each log-smooth \(X/k\), strict base change
induces an equivalence on the small log \'etale categories with
their covering families and the Zariski log-structure convention.
\end{lemma}

\begin{proof}
Use the cycle base-change functor on ordinary finite correspondences
\cite[\S1]{Suslin}. The projection
\(X_L\times_L Y_L\to X\times_k Y\) is an integral universal
homeomorphism. It identifies integral supports by taking the reduced
inverse image. Finiteness and surjectivity of a support over a
component of the source are preserved and reflected. For reflection,
the reduced inverse image maps integrally and surjectively to the
original support. If the reduced inverse image is finite over the
new source, its composite map to the original source is integral.
Universal closedness and the condition of having finite fibers then
descend along the integral surjection onto the original support.
The original support is separated and of finite type over the
original source, so it is proper and quasi-finite, hence finite.

Admissibility may be tested on the normalizations. The map from the
normalization of the reduced inverse image to the original
normalization is integral and radicial. In affine charts this follows
by taking integral closure in the purely inseparable extension of
the function field. Every element of that integral closure has a
\(p\)-power in the original normal ring. At corresponding local
rings \(A\subset B\), a fraction in \(\operatorname{Frac}(A)\) is
a unit of \(B\) exactly when it is a unit of \(A\), by integrality
and normality applied to the fraction and its inverse.

For clarity, admissibility implies that the generic point of a
support lies in the product of the trivial loci. This condition is
preserved and reflected. On the normalization, a section of the
log structure pulled back from the source is locally represented
by a chart element and a unit. At corresponding points of the two
normalizations the characteristic monoids of this log structure
are identical, and the map on units is injective. Thus the map
on pulled-back log monoids is injective. If an upstairs lift of a
target chart generator \(q\) is represented by a source chart
element \(b\) and a unit, that unit must be the fraction \(q/b\),
where \(q,b\) here denote their structure functions. They are
generically nonzero. The unit test shows that this lift descends
to the original normal local ring. The chart relations and the
agreement of the lifts on overlaps descend by injectivity of the
map on log monoids. This descends the log morphism itself.
Thus admissible supports correspond bijectively. The coefficient of a
scalar-extended support is the length of its local ring at the generic
point. For finite \(L/k\), if \(E\) is the function field of the support,
\(E\otimes_k L\) is an Artinian local \(E\)-algebra of dimension
\([L:k]\), and its residue field is purely inseparable over \(E\).
The length is the quotient of these two degrees, hence a power of
\(p\), including \(1\). For arbitrary \(L\), the ideal of the reduced
inverse image descends to a finite stage. After enlarging that stage,
its quotient is the reduction of the inverse image at that stage.
Indeed, on each of finitely many affine charts the generators of the
ideal are nilpotent after extension to \(L\), hence already nilpotent
at the finite stage by faithful flatness. The quotient is reduced
by the same faithful flatness. Its extension to \(L\) is reduced by
construction, so further extension does not change the generic
length. The resulting coefficient is a unit in \(R\).
Ordinary scalar extension respects composition,
and composition of finite log correspondences is the underlying
cycle composition on the trivial loci. This gives full faithfulness.

A universal homeomorphism of schemes is a map that is a
homeomorphism on spaces after every scheme base change. The
projection \(X_L\to X\) is strict and is an ordinary universal
homeomorphism on underlying schemes. By \cite[Remark~2.2(1), p.~240]{Vidal}
it is an fs
Kummer universal homeomorphism. Here the latter condition means
that the characteristic monoid maps are Kummer and the map
remains a homeomorphism on spaces after every fs log base change.
Vidal proves that fs base change gives equivalences of the small
log \'etale sites and the small Kummer \'etale sites
\cite[Theorem~0.1]{Vidal}.

Vidal allows log structures given by charts on \'etale neighborhoods.
We verify that his equivalence restricts to the Zariski log structures
used in the book. Strict base change preserves such structures.
For reflection, suppose an object \(V\to X\) in Vidal's category
has Zariski log structure after base change. Ordinary \'etale
topological invariance \cite[Tag~04DZ]{Stacks} identifies the
small \'etale sites of \(\underline V\) and \(\underline V_L\), and
their Zariski opens correspond. Strictness identifies their
characteristic monoid sheaves under this equivalence. Around a
point of \(V_L\), choose a Zariski chart \(P\) that is sharp
(its only invertible element is its identity), as in
\cite[Remark~A.5.8]{BPO}. Its group \(P^{\mathrm{gp}}\) is free abelian
of finite rank. The chart induces a homomorphism from \(P\) to
the sections of the characteristic monoid on the corresponding
open of \(V\).

We can lift this homomorphism to a chart after shrinking that
open. For each basis element of \(P^{\mathrm{gp}}\), the lifts to
\(M_V^{\mathrm{gp}}\) of its image in
\(\overline M_V^{\mathrm{gp}}\) form an
\(\mathcal O_{\underline V}^\times\)-torsor, a sheaf of local lifts
on which units act freely and transitively. Such a torsor is the
sheaf of bases of a line bundle and has sections Zariski locally.
Choose lifts of all basis elements on a common neighborhood.
They define a group homomorphism \(P^{\mathrm{gp}}\to
M_V^{\mathrm{gp}}\) whose restriction to \(P\) lands in \(M_V\).
Indeed, for each element of \(P\), a monoid lift exists \'etale
locally and differs from the group lift by a unit. The resulting
prelog map is a chart, since the induced map on characteristic
monoids is the chart isomorphism after identifying \'etale
stalks with those upstairs, and adjoining units supplies the
isomorphism on units. Thus the Zariski chart condition descends.

The separatedness and finite-type conditions on the underlying
schemes descend by faithful flatness. In the full sites, covering
families are preserved and reflected by Vidal's theorem. For the
sites with Zariski charts this can also be seen on the generators.
A Kummer \'etale family is a cover if it is jointly surjective.
The equivalence for Kummer \'etale objects and the homeomorphisms
on their underlying spaces preserve and reflect this condition.
For a dividing cover, the monomorphism condition says that its
diagonal into its fs self-product is an isomorphism, so it is
preserved and reflected by the equivalence. Properness descends
by faithful flatness on schemes, and surjectivity is unchanged.
Here strict base change has the usual scheme fiber product as
its underlying scheme \cite[Lemma~A.5.6]{BPO}.
We obtain the asserted equivalence for the small sites in
\cite[\S3.2]{BPO}, with the book's convention.
\end{proof}

\begin{proof}[Completion of the proof of Theorem~\ref{thm:sheaf}]
Put \(L=k^{\mathrm{perf}}\). Write \(b\) for scalar extension on
\(R\)-linear log correspondences. Restriction of presheaves with
transfers along \(b\) has a left adjoint, given by additive left
Kan extension along the opposite of \(b\). Denote by \(\mathcal P\) the
image of \(F\) by this left adjoint, and by \(F_L\) its log \'etale
sheafification with transfers. Full faithfulness in
Lemma~\ref{lem:perfection} gives
\(\mathcal P(X_L)=F(X)\) on every object over \(k\). This identification is
natural for log transfers.

The assignment \(F\mapsto F_L\) is left adjoint on sheaves to
restriction along \(b\), because scalar extension preserves
covering families and their fiber products. Therefore it
preserves colimits. The construction is additive, so the image
of \(F[n]\) is killed by \(n\). It follows that \(F_L\), the
filtered colimit of these images, has elementwise torsion of
order prime to \(p\).

The map \(\mathcal P(X_L)\to F_L(X_L)\) is an isomorphism for every
\(X/k\). In fact, on the small log \'etale site over \(X_L\)
the presheaf \(\mathcal P\) identifies with the restriction of \(F\) on
the small site over \(X\), by Lemma~\ref{lem:perfection}. It is
already a sheaf on this site. Sheafification with log transfers
agrees with ordinary log \'etale sheafification on the underlying
presheaves \cite[Theorem~4.5.7]{BPO}, and the small-site
restriction commutes with sheafification \cite[\S3.2]{BPO}.
Consequently there are also natural identifications
\[
 \R\Gamma_{\letop}(X_L,F_L)\simeq\R\Gamma_{\letop}(X,F)
\]
from the equivalence of the small sites and
\cite[Proposition~3.2.4]{BPO}. They apply to \(X\times\bx\) as
well, so strict \(\bx\)-invariance holds on scalar-extended
objects over \(L\).

To verify strict invariance on every \(Y/L\), there is a finite
Zariski open cover of \(Y\) by opens that are scalar extensions
of objects over \(k\). Indeed, by the log-smooth chart
description \cite[Lemma~A.5.9]{BPO}, Zariski locally \(Y\)
admits a strict smooth map to a toric chart
\(\Spec L[P]\) with its standard log structure. Further
Zariski localization of the source gives a strict \'etale
map to \(\Spec L[P]\times\A^d_L\) for some \(d\), by choosing
relative \'etale coordinates for a smooth morphism. This
target is scalar-extended from \(k\). The strict \'etale
object descends by ordinary topological invariance of the
\'etale site \cite[Tag~04DZ]{Stacks}, carrying the pulled-back
log structure. Every finite intersection of the resulting opens
is itself an open in one of those charts and descends, since
scalar extension by \(L/k\) induces a homeomorphism on underlying
spaces. The functor \(Y\mapsto Y\times\bx_L\) preserves these
covers and their intersections. Zariski descent for the two
cohomology complexes thus proves strict \(\bx\)-invariance of
\(F_L\).

The perfect-field result in Section~\ref{sec:perfectresult} applies
to \(F_L\). To descend its
conclusion, write \(b_\sharp\) for the construction
\(F\mapsto F_L\) and \(b^*\) for restriction along \(b\).
The calculation on the small log \'etale site of \(X_L\)
shows that the adjunction map \(E\to b^*b_\sharp E\) is an
isomorphism for every log \'etale sheaf \(E\) with
\(R\)-linear log transfers over \(k\). On the small sites over
the fields, scalar extension is an equivalence by ordinary
\'etale topological invariance. Write \(d_\sharp\) for the
induced extension of sheaves and \(d^*\) for its inverse.
The restriction functors satisfy
\(\eta_k^*b^*=d^*\eta_L^*\), where the subscripts distinguish
the ground fields. Taking left adjoints gives
\[
 b_\sharp\eta_{k,\sharp}\simeq \eta_{L,\sharp}d_\sharp .
\]
Moreover, \(d_\sharp\eta_k^*F\simeq\eta_L^*F_L\), by the
isomorphism \(F\simeq b^*F_L\) on finite \'etale objects.
These identifications take the extension of the counit over
\(k\) to the counit over \(L\). The latter is an isomorphism,
and applying \(b^*\) proves the result over \(k\).

For completeness, ordinary pullback from the small sites is
compatible with these identifications on the underlying sheaves.
For a small-site sheaf \(G\), ordinary pullback is the
sheafification of the presheaf extension formed from the
groups \(G(T)\) indexed by maps \(X\to T\) to finite \'etale
objects \(T\) with trivial log structure. Such log morphisms
are exactly the ordinary maps \(\underline X\to T\).
Ordinary topological invariance on \(\underline X\) and on the
ground field identifies this indexing category and coefficient
diagram with the ones over \(X_L\) and \(L\). Sheafification on
the corresponding small log \'etale sites then agrees. Over
\(L\), the pullback with transfers agrees with ordinary pullback
by \cite[Proposition~8.5.6]{BPO} for bounded torsion, and
by filtered colimits for elementwise torsion. The preceding
compatibilities give the same identification over \(k\).
The counit over \(k\) commutes with the coefficient action by naturality.
\end{proof}

\section{Rigidity for unbounded torsion motives}\label{sec:derived}

\subsection{Boundary neighborhoods and the cyclic action}
\label{sec:normal}\label{sec:cyclic}

The following elementary construction will be used twice.
Here \(M=M_t\), for \(t=\dNis,\detop,\letop\), with arbitrary coefficients.
It only uses strict Nisnevich distinguished squares: an open
immersion, a strict \'etale map that is an isomorphism on the
reduced closed complement, and their fiber product. Motives
send these squares to derived pushouts
\cite[Proposition~5.2.10]{BPO}.

\begin{lemma}\label{lem:normal}
Let \(X=(\underline X,Z+B)\) have smooth underlying scheme and
strict normal-crossing boundary. On a coordinate neighborhood
where \(Z=\{t=0\}\), the relative object
\[
 N(X,Z)=\cofib(M(X-Z)\to M(X))
\]
is identified by strict Nisnevich excision with
\[
 M(Z,B|_Z)\otimes
 \cofib\bigl(M(\Gm)\to M(C)\bigr).
\]
The comparison is natural for the power cover \(u^m=t\).
on the relative objects of that cover it is equivariant for
the action \(u\mapsto\zeta u\).
\end{lemma}

\begin{proof}
Choose an \'etale coordinate map \(\underline X\to\A^d\)
whose first coordinate is \(t\) and whose boundary is a union
of coordinate hyperplanes. Set \(X_1=(Z,B|_Z)\times C\).
It too maps strictly \'etale to the same coordinate log scheme.
In the scheme fiber product of \(\underline X\) and
\(\underline X_1\), the inverse image of the first divisor is
\(Z\times_{\A^{d-1}}Z\). The diagonal is open and closed there,
because \(Z\to\A^{d-1}\) is \'etale and separated.
Remove its complement in this inverse image, a closed subset
of the full fiber product. The resulting open \(W\) maps
strictly \'etale to \(X\) and \(X_1\) and is an isomorphism on
the designated divisor. The two distinguished squares identify
their relative cofibers. The product formula gives the display.

For the power cover pull back \(W\) along \(u^m=t\).
Its maps to the two covered objects are strict \'etale and
are isomorphisms on the reduced designated divisors. The
deck action is the action on \(u\). Every comparison is
equivariant, and the relative map is the identity on
\(M(Z,B|_Z)\) tensored with the power map on the line.
\end{proof}

Assume \(n\geq2\), \(n{R}=0\), \(n\) is invertible in \(k\),
and \(k\) contains a primitive \(n\)-th root of unity \(\zeta\).
Let \(G\) be the cyclic group of order \(n\) with generator
\(g\), acting on \(P\) by \(t\mapsto\zeta t\).
It fixes \(0,\infty\). Let \(D_0=\{t^n=1\}\), a disjoint
union of \(n\) rational points forming one free \(G\)-orbit.
Write \(j:U=P-D_0\hookrightarrow P\) and \(i:D_0\hookrightarrow P\).
On \(U\) the function \(f=(1-t^n)^{-1}\) is invariant under
\(G\) and maps \(0,\infty\) to \(1,0\), respectively.
At infinity its expression in the coordinate \(v=1/t\)
is \(v^n/(v^n-1)\), so it is regular there.

For an object \(B\) with \(G\)-action, the homotopy orbits
\(B_{hG}\) are defined as the realization of the action bar
construction. Its term in simplicial degree \(a\) is
\(\bigoplus_{G^a} B\). One end face acts by the
corresponding group element, the other end face deletes
one entry, interior faces multiply adjacent entries, and
degeneracies insert the identity. We use the convention
in degree one with faces acting by the group element and
by the identity, in that order. For the unit with trivial
action the generator \([g]\) of the normalized
degree-one chains is a cycle because its two faces agree.
Thus there is a map
\(\sigma:\one[1]\to\one_{hG}\) defined with integral
coefficients and preserved by the changes in
Section~\ref{sec:coeff-fields}.
The fixed points induce \(s_0,s_\infty:\one_{hG}\to M_t(P)_{hG}\).
Set
\[
 d_t=(s_0-s_\infty)\sigma:\one_t[1]\longrightarrow M_t(P)_{hG}.
\]
The integral dividing Nisnevich version is \(d_{\dNis,\Z}\).
For \(x\in A(k)\) the graph of the point defines a map
\([x]:\one_t\to M_t(A)\).

\begin{proposition}\label{prop:pairing}
For \(t=\detop\) or \(\letop\) there is a map
\(\alpha_t:M_t(P)_{hG}\to M_t(A)[1]\) satisfying
\[
 \alpha_t d_t=([1]-[0])[1].
\]
\end{proposition}

\begin{proof}
First work in the derived category of sheaves before
imposing the projection relations. Let \(F_P\) be the
restriction to the ordinary small \'etale site of \(P\) of
\(\gamma^*a_t{R}_{\ltr}(A)\). Its \(G\)-action is
induced by the action on \(P\).
For an equivariant sheaf \(F\), write
\[
 H_G^r(P,F)=H^r\bigl(\R\Gamma(P_{\et},F)^{hG}\bigr).
\]
Here \({}^{hG}\) is derived group invariants, computed
by the product total complex of the group cochain
construction applied to an equivariant derived
sections complex. In particular \(H_G^r\) gives long
exact sequences for equivariant short exact sequences.
These equivariant derived sections can be formed by
resolving in equivariant sheaves. Forgetting the group
action preserves injectives since it has an exact left
adjoint given by the sum of the translates.

For each \(a\), formula \eqref{eq:eval} identifies the
derived mapping complex from
\(a_t{R}_{\ltr}(P\times G^a)\) to
\(a_t{R}_{\ltr}(A)\) with the product indexed by
\(G^a\) of the derived sections of \(F_P\).
By naturality, the action and all simplicial maps give
the group cochain construction. Mapping out of a
realization takes the derived limit of this
cosimplicial diagram of mapping complexes. Thus
\begin{equation}\label{eq:bar}
 H_G^r(P,F_P)\simeq
 \Hom_{D(\Shv_t^{\ltr}(k,{R}))}
 \bigl((a_t{R}_{\ltr}(P))_{hG},a_t{R}_{\ltr}(A)[r]\bigr).
\end{equation}
The comparison also identifies the maps induced by
the two fixed points with restriction to them.

Let \(\chi:G\to({R},+)\) send \(g\) to \(1\).
This character is defined because \(n{R}=0\).
It gives a degree-one cohomology class with constant
coefficients by using constant sections as group
cochains. Equivalently it is the class of the extension
of equivariant constant sheaves
\(0\to{R}\to{R}^2\to{R}\to0\)
in which \(g\) acts on the middle term by
\((a,b)\mapsto(a+b,b)\).
Its pullback to \(D_0\) is zero. In fact, at
\(g^a x_0\), for any chosen \(x_0\in D_0(k)\),
the splitting \(b\mapsto(ab,b)\) of this extension
defines an equivariant splitting on the whole orbit.
The product \(ab\) means \(a\)-fold addition in
\({R}\), which is independent of the representative
of \(a\bmod n\).

On the small \'etale site there is an exact sequence of
equivariant sheaves
\[
 0\longrightarrow j_!{R}_U\longrightarrow{R}_P
   \longrightarrow i_*{R}_{D_0}\longrightarrow0 .
\]
Here \(j_!\) is extension by zero, left adjoint to
restriction to \(U\). Exactness can be checked at
geometric stalks in \(U\) and \(D_0\).
The cohomology of \(i_*{R}_{D_0}\) equals
the following calculation of the cohomology on \(D_0\) with its \(G\)-action.
direct image for a closed immersion is exact on these
sites, and it is right adjoint to exact inverse image.
The long exact sequence therefore gives a class
\(\beta\in H_G^1(P,j_!{R}_U)\) lifting the class
of \(\chi\) on \(P\).
Both fixed points lie in \(U\). At either point the
inclusion \(j_!{R}_U\to{R}_P\) restricts to the
identity of the constant coefficient sheaf, so the
restriction of any such \(\beta\) is the class of \(\chi\).

The graph of \(f\) is an invariant section of \(F_P|_U\).
By adjunction it gives an equivariant map
\(j_!{R}_U\to F_P\) of small-site sheaves. The
image of \(\beta\) is a derived map with transfers
by \eqref{eq:bar}, and applying localization gives
\(\alpha_t\).
The two restrictions of the image of \(\beta\) are the
group one-cocycles \(h\mapsto\chi(h)[1]\) and
\(h\mapsto\chi(h)[0]\) with values in sections over
\(\Spec k\), placed in the derived sections complex.
Here \([1]\) and \([0]\) inside the expressions denote
the sections given by the graphs of the two affine
points. Evaluating on the degree-one chain \([g]\)
gives those two sections because \(\chi(g)=1\).
By \eqref{eq:bar} these evaluations are exactly the
composites of \(\alpha_t\) with \(s_0\sigma\) and
\(s_\infty\sigma\), with the cohomological shift.
Subtracting proves the formula.
\end{proof}

\subsection{Affine-line invariance and removal of the boundary}
\label{sec:integral-line}

Continue with a perfect field containing \(\zeta\).
Write \(M_{\Z}(X)\) for the integral dividing
Nisnevich motive and \(\one_{\Z}=M_{\Z}(\Spec k)\).
The reduced projective line motive here is
\[
 Q=\fib(M_{\Z}(P)\longrightarrow\one_{\Z}).
\]
The fixed point \(\infty\) splits the structural map
equivariantly. Taking bar realizations is exact on
fiber sequences of complexes, so there is a splitting
\(M_{\Z}(P)_{hG}\simeq Q_{hG}\oplus(\one_{\Z})_{hG}\).
The two fixed-point maps have the same composite to
\((\one_{\Z})_{hG}\). Hence \(d_{\dNis,\Z}\) comes from
a map \(\widetilde d:\one_{\Z}[1]\to Q_{hG}\).

Write \(\pi_b\) for the homology sheaf in degree \(b\),
equivalently the cohomology sheaf in degree \(-b\),
of the \(\bx\)-local representative of an integral
dividing Nisnevich motive. Binda--Merici's homotopy
\(t\)-structure means that the ordinary cohomological
truncations preserve those local complexes
\cite[Theorem~5.7]{BM}. In particular a represented
motive is connective, meaning that its negative
homology sheaves vanish. Indeed the inclusion of local
complexes preserves both truncations and the derived
localization, as its left adjoint, preserves connective
objects. One may also use the connectivity theorem
\cite[Theorem~3.2]{BM}. The underlying-smooth site used
there gives the same category by
\cite[Proposition~5.4.2]{BPO}.

Merici computes \(\pi_0 M_\Z(P)=\Z\) via the structural
map and \(\pi_1M_\Z(P)\simeq\omega^*\Gm\) via the
tautological line bundle
\cite[Proposition~3.6 and the preceding calculation]{Merici}.
Here \(\omega^*\Gm\) is the sheaf that evaluates on \(X\)
as the group of invertible regular functions on \(X^\circ\),
with ordinary transfers on the interior. The integral
unit is concentrated in degree zero in the local
presentation. The splitting and connectivity give
\begin{equation}\label{eq:Q}
 \pi_bQ=0\quad(b\leq0),\qquad \pi_1Q\simeq\omega^*\Gm.
\end{equation}
The identification used in degree one is compatible
on field sections with extension of perfect fields.
Indeed the map is the class of the tautological bundle,
using \(H^1_{\mathrm{Nis}}(\Pp^1,\Gm)=\operatorname{Pic}(\Pp^1)\).
It can equivalently be constructed by sending the log
motive to the ordinary motive of its interior and then
taking the ordinary class of that bundle. This is the
adjunction description of the map in \cite[(3.5.2)]{Merici}.
Taking interiors and pulling back the bundle commute with
extension of fields. The induced map on \(\Gm\) sends a
unit to the same unit after scalar extension.

\begin{lemma}\label{lem:low}
The map
\[
 k^\times=\Hom(\one_\Z[1],Q)
       \longrightarrow\Hom(\one_\Z[1],Q_{hG})
\]
induced by the degree-zero term in the bar construction
is surjective.
\end{lemma}

\begin{proof}
Evaluation at \(\Spec k\) is exact on dividing Nisnevich
sheaves: any small Nisnevich cover of a field has a
section and dividing covers there are isomorphisms.
Thus \(V=\RHom(\one_\Z,Q)\) has \(H_b(V)=0\) for
\(b\leq0\) and \(H_1(V)=k^\times\), by \eqref{eq:eval}
and \eqref{eq:Q}. Here \(H_b\) is ordinary complex
homology in degree \(b\).
The unit is compact by \cite[Proposition~5.2.5]{BPO},
meaning that maps from it and its shifts commute with
direct sums. Hence its derived mapping functor
commutes with direct sums. It also preserves finite
derived colimits, since it is exact, and consequently
commutes with simplicial realizations. One can see
the last assertion by building a realization from its
finite-degree skeleta and taking their sequential
derived colimit, which is a cone of a map of direct sums.
So \(\RHom(\one_\Z,Q_{hG})\simeq V_{hG}\).

Filter the normalized bar complex by simplicial degree.
Its homology spectral sequence has
\(E^2_{a,b}=H_a(G,H_b(V))\), where \(H_a(G,-)\) is the
group homology computed by the action bar complex of a
module. It is enough to take \(V\) concentrated
homologically in degrees at least one, by truncating
below that degree before applying the bar construction.
The filtration is then finite in each total degree.
In total degree one the single term is \(E^2_{0,1}\).
It is the quotient of \(H_1(V)\) by all action
differences, and no higher differential enters or
leaves it. The map from the degree-zero bar term
surjects onto this quotient, as required.
\end{proof}

\begin{proposition}\label{prop:affine}
For any field and coefficients as in
Theorem~\ref{thm:derived}, the structure map
\(M_t(A)\longrightarrow\one_t\) is an isomorphism
for \(t=\detop\) and \(t=\letop\).
\end{proposition}

\begin{proof}
The case \(n=1\) has zero coefficients. For \(n\geq2\)
and perfect \(k\) containing \(\zeta\), choose
\(u\in k^\times\) lifting \(\widetilde d\) in
Lemma~\ref{lem:low}. Choose a finite separable extension
\(L/k\) containing an \(n\)-th root of \(u\). This is
possible because \(u\ne0\) and \(n\) is invertible in \(k\).
Under the identification of field sections in \eqref{eq:Q},
the pullback of the lift is \(u\in L^\times\).
Since multiplication by \(n\) in that group is
the \(n\)-th power map, the pullback of the lift is an
\(n\)-fold multiple. The same holds for \(d_{\dNis,\Z}\)
over \(L\). Changing coefficients and topology sends it
to \(d_t=0\) over \(L\). Proposition~\ref{prop:pairing}
then gives \([1]=[0]\) on the affine line over \(L\).

If \(p:M_t(A_L)\to\one_t\) is the structure map,
affine-coordinate multiplication \(m:A_L\times A_L\to A_L\)
gives
\[
 \operatorname{id}-[0]p
   =M_t(m)\bigl(\operatorname{id}\otimes([1]-[0])\bigr)=0.
\]
Since \(p[0]\) is the identity, \(p\) is invertible.
The cofiber of the corresponding structure map
over \(k\) therefore vanishes by Lemma~\ref{lem:descent}.
For a perfect field that does not yet contain \(\zeta\),
the finite separable extension obtained by adjoining it
and the same descent argument give the result.

For arbitrary \(k\), its prime field \(b\) is perfect and
\(n\) is invertible in \(b\). Apply the already proved case
to \(b\) with the same coefficient ring. The map over \(k\)
is its pullback by Lemma~\ref{lem:pullback}.
\end{proof}

\begin{lemma}\label{lem:boundary}
Under the coefficient hypothesis of Theorem~\ref{thm:derived},
the inclusion induces \(M_t(X^\circ)\simeq M_t(X)\) for
\(t=\detop,\letop\).
\end{lemma}

\begin{proof}
The open charts \(A,C\) of \(\bx\) meet in \(\Gm\).
By Proposition~\ref{prop:affine} the map \(M_t(A)\to M_t(\bx)\)
is an equivalence: both structural maps to the unit are
equivalences. Mayer--Vietoris gives
\begin{equation}\label{eq:boundary-line}
 M_t(\Gm)\xrightarrow{\sim}M_t(C).
\end{equation}
Make a log blow-up \(Y\to X\) with smooth underlying scheme,
using \cite[Proposition~A.10.2]{BPO}. It preserves interiors
and motives. Work on its normal-crossing coordinate
neighborhoods. For a boundary component \(Z\), Lemma~\ref{lem:normal}
and \eqref{eq:boundary-line} give \(N(Y,Z)=0\).
Remove the other components successively. Intersections of
the neighborhoods inherit the coordinates, so Zariski
\v{C}ech descent proves the assertion on \(Y\), hence on \(X\).
\end{proof}

\subsection{Comparison with ordinary motives}
\label{sec:ordinary-comparison}

\begin{proof}[Proof of Theorem~\ref{thm:derived}]
Work under the coefficient hypothesis of the theorem and take
\(t=\detop\) or \(\letop\).
The inclusion \(\lambda\) of smooth schemes with
trivial log structure and the functor
\(\omega:X\mapsto X^\circ\) induce functors on
correspondences by \cite[Chapter~2]{BPO}.
The inclusion takes ordinary \'etale covers to
strict \'etale covers. Taking interiors sends
covers of either of the two log sites to ordinary
\'etale covers. Both preserve the fiber products
used to form hypercovers. In the case of interiors,
this follows also because a map from a trivial-log
scheme factors uniquely through the interior.
The derived left extensions on sheaves are the
functors of \cite[\S4.3]{BPO}.

Taking interiors takes the \(\bx\)-projection relations
to the \(\A^1\)-projection relations. For the inclusion
of trivial-log schemes, the \(\A^1\)-projection
relations hold in the target by
Proposition~\ref{prop:affine} and the product formula.
Consequently the functors descend to
\[
 I_t:\DM^\eff_{\et}(k,{R})\longrightarrow
        \logDM^\eff_t(k,{R}),\qquad
 J_t:\logDM^\eff_t(k,{R})\longrightarrow
        \DM^\eff_{\et}(k,{R}).
\]
The equality \(\omega\lambda=\operatorname{id}\)
gives \(J_tI_t\simeq\operatorname{id}\).
The open immersions \(\lambda\omega(X)\to X\)
are natural for log correspondences. For a
correspondence from \(X\) to \(Y\) the two
composites both restrict its cycle to the source
interior, since the log condition then forces it
to land in \(Y^\circ\). This gives a natural
transformation \(I_tJ_t\to\operatorname{id}\).
It is an isomorphism on every \(M_t(X)\) by
Lemma~\ref{lem:boundary}.

These motives generate under shifts, cofibers and
arbitrary coproducts. Indeed the corresponding
representable presheaves with transfers are projective
generators of the presheaf category, and hence
generate its unbounded derived category. For the
last implication, bounded-above complexes admit
bounded-above projective resolutions, whose finite
subcomplexes obtained by cutting off the low degrees
recover the resolution by a sequential colimit.
Each projective is a summand of a sum of
representables. A subcategory closed under cofibers
and coproducts is closed under such summands,
since a summand is the sequential derived colimit
of its idempotent projection. Every complex is the
sequential derived colimit of its bounded-above
cohomological truncations. The localization functor
carries these generators to \(M_t(X)\). Thus
\(I_tJ_t\to\operatorname{id}\) is an isomorphism on
the whole category. No cohomological-dimension
hypothesis on \(k\) enters this generation argument.
The functors \(I_t\) for the two topologies agree
under change of topology because they are induced
by the same inclusion.

Ordinary unbounded \'etale rigidity
\cite[Theorem~4.5.2]{CD} gives the canonical equivalence
\(\rho_!:D(k_{\et},{R})\to\DM^\eff_{\et}(k,{R})\).
Its right adjoint is restriction of the
\(\A^1\)-local complexes to the small site after
forgetting transfers. The right adjoint to \(I_t\)
is restriction of local complexes to trivial-log
smooth schemes, by the derived adjunction in
\cite[(4.3.5)]{BPO}. That restriction lands
in \(\A^1\)-local complexes by
Proposition~\ref{prop:affine} and derived evaluation.
Hence the right adjoint of \(I_{\letop}\rho_!\)
is the small-site restriction defining the
right adjoint of \(\eta_\sharp\). Uniqueness of
left adjoints identifies the two functors, so \(\eta_\sharp\)
is an equivalence.
\end{proof}

\section{Unrestricted change of topology}\label{sec:topology}

We now remove the torsion hypothesis from the comparison of
the two logarithmic topologies. The proof of
Theorem~\ref{thm:derived} supplies \(T_{\mathbb{F}_\ell}\)
for every prime \(\ell\) invertible in \(k\). The new geometric
input is descent for a root cover once its degree is a unit
in the coefficient ring.
Write \(\mathcal C_R=\logDM_{\detop}^{\eff}(k,R)\) and
\(\mathcal D_R=\logDM_{\letop}^{\eff}(k,R)\) in this section.

\subsection{Descent for root covers}

\begin{lemma}\label{lem:root-line}
Let \(m>0\) be invertible in \(k\) and \(R\).
For \(q_m:C\to C\), \(u\mapsto u^m\), the augmented fs
\v{C}ech nerve has
\[
 |M_{\detop}(C^{[\bullet]})|\xrightarrow{\sim}M_{\detop}(C)
\]
in \(\mathcal C_R\).
\end{lemma}

\begin{proof}
Finite separable extension detects zero objects by
Lemma~\ref{lem:descent}, so assume \(\mu_m\subset k\) and
write \(G=\mu_m(k)\). The fs self-products of the root cover
are \(C\times G^a\). Their realization is the group quotient \(M(C)_{hG}\).
Since \(m\) is invertible in \(R\), this object is the direct
summand given by the idempotent \(m^{-1}\sum_{g\in G}g\).
This description is valid for arbitrary complexes: the
averaging operator contracts the positive-degree part of the
group bar resolution.

Use the power map on the two charts of \(\bx\), obtaining
a diagram of Mayer--Vietoris triangles
\[
\begin{matrix}
 M(\Gm)_{hG}&\longrightarrow&M(A)_{hG}\oplus M(C)_{hG}
       &\longrightarrow&M(\bx)_{hG}&\longrightarrow\\
 \downarrow&&\downarrow&&\downarrow\\
 M(\Gm)&\longrightarrow&M(A)\oplus M(C)
       &\longrightarrow&M(\bx)&\longrightarrow .
\end{matrix}
\]
The first vertical map is an equivalence because the power
map on \(\Gm\) is a strict finite \'etale \(G\)-torsor.
On \(A\), the transpose \(\tau\) of the power map is an
admissible finite correspondence, since both log structures
are trivial. The cycle equalities
\[
 M(q_m)\tau=m\,\mathrm{id},\qquad
 \tau M(q_m)=\sum_{g\in G}M(g)
\]
show that \(M(A)_{hG}\to M(A)\) is an equivalence.
On \(\bx\), the structural equivalence to the unit is
equivariant, and \(\one_{hG}\simeq\one\) by averaging.
The remaining vertical map, on \(C\), is therefore an
equivalence. The coordinate at infinity replaces \(\zeta\)
by \(\zeta^{-1}\), which gives the same group quotient.
\end{proof}

\begin{lemma}\label{lem:root-chart}
Let \(X=(\underline X,D)\) have smooth underlying scheme and
strict normal-crossing boundary. Let \(Y\to X\) be obtained
by adjoining \(m\)-th roots of boundary coordinates on a
chart, with \(m\) invertible in \(k\) and \(R\).
Its augmented fs \v{C}ech nerve is a colimit diagram of
dividing \'etale motives with coefficients in \(R\).
\end{lemma}

\begin{proof}
Start with one coordinate \(t\), defining \(Z\subset D\),
and put \(B=D-Z\). After a finite separable extension the
cover is Galois of group \(G=\mu_m(k)\).
On \(X-Z\) it is strict finite \'etale, so the comparison
on the open complements is already an equivalence.
Lemma~\ref{lem:normal} identifies the comparison on relative
cofibers with
\[
 M(Z,B|_Z)\otimes
 \left(N(C)_{hG}\longrightarrow N(C)\right),
 \qquad N(C)=\cofib(M(\Gm)\to M(C)).
\]
It is an equivalence by Lemma~\ref{lem:root-line} and
\'etale descent for the power map on \(\Gm\).
The cofiber triangles give the result on \(X\).
The construction is local on coordinate neighborhoods. Their inverse images are \(G\)-stable and their intersections
have the same description, so Zariski descent glues it.

For several coordinates adjoin the roots successively.
The intermediate schemes remain smooth with strict
normal-crossing boundary. The commuting deck groups give
the product group. Iterated group quotients are its single
bar realization. The one-coordinate equivalences therefore
give the assertion. Descend the finite separable extension
by Lemma~\ref{lem:descent}.
\end{proof}

For the following descent argument, suppose that every integer invertible in
\(k\) is invertible in \(R\). Let
\[
 \mathcal A=\Shv_{\detop}^{\ltr}(k,R),\qquad
 \mathcal B=\Shv_{\letop}^{\ltr}(k,R),\qquad
 a:\mathcal A\rightleftarrows\mathcal B:j
\]
be sheafification and the fully faithful inclusion.
A \emph{full root cover} on a normal-crossing chart with
boundary coordinates \(t_1,\ldots,t_r\) means
\[
 \Spec\bigl(\mathcal O_X[u_1,\ldots,u_r]/
               (u_1^m-t_1,\ldots,u_r^m-t_r)\bigr)\longrightarrow X,
\]
with the boundary given by the vanishing of the \(u_i\).
Here \(m\) is invertible in \(k\). A common multiple absorbs
different root orders. A strictly henselian local ring is a
henselian local ring with separably closed residue field.
The henselian condition means that coprime monic factorizations
modulo the maximal ideal lift to the ring.
The extra covers can be tested, dividing \'etale locally,
on the full root covers of normal-crossing coordinate
charts. To see this, first use
\cite[Proposition~4.4.1]{BPO} to make a log \'etale cover
Kummer \'etale after a target dividing cover, then make
the target underlying scheme smooth. The Kummer chart
criterion supplies a full root cover dominating the cover
strict \'etale locally. Equivalently, over a strictly
henselian neighborhood the standard root covers form a
cofinal family and have strictly henselian local underlying schemes.
See \cite[Corollary~3.6 and \S4.7(d)]{Illusie}.
After a further strict \'etale cover their deck groups
are constant. We use these root covers in the following lemmas.

\begin{lemma}\label{lem:exact-inclusion}
The functor \(j\) is exact and preserves arbitrary coproducts.
\end{lemma}

\begin{proof}
Kernels are computed on sections. Consider an epimorphism
\(F\to H\) in \(\mathcal B\) and \(s\in H(X)\).
Choose a log \'etale cover on which it lifts.
The preceding refinement, after dividing \'etale localization
of \(X\), gives a finite Galois root cover \(Y\to X\)
of order \(d\) on which a lift \(b\in F(Y)\) exists.
A standard root cover over a strict henselization is finite
with strictly henselian local source, so a strict \'etale cover containing a point over its closed
point has a section. The resulting chart morphism descends
to a strict \'etale neighborhood of finite presentation.
Thus this refinement is by actual dividing \'etale covers
and finite root covers.

The section \(d^{-1}\sum_{g\in G}g^*b\) is invariant and
still maps to \(s|_Y\). Since \(F\) is a log \'etale sheaf
it descends to \(F(X)\). This proves local surjectivity
for the dividing \'etale topology, hence exactness of \(j\).

Both topologies have finite covering refinements on their
finite-type objects, and the requisite finite intersections are quasi-compact
\cite[Proposition~3.2.5]{BPO}. Thus filtered colimits are
computed on sections for both categories and \(j\) preserves
them. Exactness gives finite coproducts, so \(j\) preserves
arbitrary coproducts.
\end{proof}

For a full root cover \(q:Y\to X\) with constant group \(G\),
write \(h_X=a_{\detop}R_{\ltr}(X)\) and set
\[
 K_q=\coker\bigl((h_Y)_G\longrightarrow h_X\bigr).
\]
Here \((h_Y)_G\) is the quotient by the action differences.
The degree \(|G|\) is invertible in \(R\).

\begin{lemma}\label{lem:root-kernel}
The arrow \((h_Y)_G\to h_X\) is injective, and \(K_q[0]\)
is the cofiber of the augmented root-cover \v{C}ech
realization in \(D(\mathcal A)\). Moreover,
\[
 \ker\bigl(D(a):D(\mathcal A)\to D(\mathcal B)\bigr)
      =\Loc\{K_q:q\text{ a full root cover as above}\}.
\]
\end{lemma}

\begin{proof}
Check injectivity first on correspondence presheaves.
Above an integral support \(W\) in \(T\times X\), the
supports in \(T\times Y\) form a single \(G\)-orbit if
any is admissible. Indeed their generic points are the
components of a \(G\)-torsor over the function field of
\(W\), since every correspondence sends the generic
interior into \(X^\circ\). Admissibility is preserved by
the deck action, so either the whole orbit is present or
none of it is. Pushforward on its coinvariants multiplies
by the residue-field degree of a component, a divisor of
\(|G|\). This is a unit of \(R\). Different orbits have
different image supports. The presheaf map is therefore
injective, and exact dividing \'etale sheafification
preserves it.

The fs \v{C}ech terms are \(Y\times G^a\).
Averaging identifies their realization with \((h_Y)_G[0]\).
The injectivity just proved identifies the augmented
cofiber with \(K_q[0]\). Since \(q\) is a log \'etale
cover, \(aK_q=0\).

Conversely, if \(aH=0\) for \(H\in\mathcal A\), each
section becomes zero on a log \'etale cover. Refine that
cover as above. Dividing \'etale locally a full root
cover kills the section, so its representing map
\(h_X\to H\) factors through \(K_q\). These factorizations
give an epimorphism from a sum of the \(K_q\) to \(H\).
Its kernel is again killed by \(a\). Repeating gives a
resolution by sums of \(K_q\). Thus every such sheaf
belongs to \(\Loc\{K_q\}\).

Both \(a\) and \(j\) are exact and preserve coproducts.
They induce the same adjunction on unbounded derived
categories, with \(D(j)\) fully faithful. Let
\[
 Q(E)=\fib(E\to jaE).
\]
It preserves coproducts and cofibers and has image
\(\ker D(a)\). The objects \(h_X\) generate \(D(\mathcal A)\). This follows either by sheafification from presheaf
generators or from \cite[Proposition~B.1.5]{BPO}.
For \(h_X\), the only possible cohomology of \(Q(h_X)\)
is in degrees zero and one, and both sheaves are killed
by \(a\). The preceding resolution argument places
\(Q(h_X)\) in \(\Loc\{K_q\}\). Applying the
colimit-preserving functor \(Q\) to the generators proves
the asserted equality for all unbounded objects.
\end{proof}

\begin{proposition}\label{prop:inverted-tame}
If every integer invertible in \(k\) is invertible in \(R\),
then \(T_R:\mathcal C_R\to\mathcal D_R\) is an equivalence.
\end{proposition}

\begin{proof}
Before imposing the projection relations, Lemma~\ref{lem:root-kernel}
describes change of topology as localization at the objects
\(K_q\). Sheafification sends each projection relation to
the same relation in the finer topology. We may therefore
impose the two sets of relations in either order.
After imposing \(\bx\)-invariance, every \(K_q\) is zero
by Lemma~\ref{lem:root-chart}. The additional localization
does nothing. The kernel description includes unbounded
complexes, so this proves the proposition in their full scope.
\end{proof}

\subsection{Integral and arbitrary coefficients}

\begin{proof}[Proof of Theorem~\ref{thm:topology}]
Put
\[
 R_0=
 \begin{cases}
 \Q,&\operatorname{char}k=0,\\
 \Z_{(p)},&\operatorname{char}k=p>0.
 \end{cases}
\]
The second ring inverts all integers prime to \(p\).
Proposition~\ref{prop:inverted-tame} proves that \(T_{R_0}\)
is an equivalence. Theorem~\ref{thm:derived} and its
compatibility with change of topology prove that
\(T_{\mathbb{F}_\ell}\) is an equivalence for every prime
\(\ell\) invertible in \(k\).

Let \(E\in\ker T_{\Z}\). Compatibility of coefficient
extension with \(T\), and Lemma~\ref{lem:coefficients}, give
\[
 E\otimes_{\Z}^{\mathbf L}R_0=0,\qquad
 E/\ell:=\cofib(\ell:E\to E)=0
       \quad(\ell\in k^\times).
\]
Thus every such \(\ell\) acts invertibly on \(E\).
In characteristic zero this means
\(E\simeq E\otimes_{\Z}^{\mathbf L}\Q=0\).
In characteristic \(p\) it means
\(E\simeq E\otimes_{\Z}^{\mathbf L}\Z_{(p)}=0\).
These identifications are the internal derived
localizations of Lemma~\ref{lem:coefficients}.

Change of topology \(T_{\Z}\) is a localization: in
presheaf models it adds the finer hypercover relations.
Since it has zero kernel it is an equivalence. It is
symmetric monoidal. Lemma~\ref{lem:coefficients} now
extends this equivalence to module categories for every
commutative ring \(R\). The resulting functor is \(T_R\),
because it is induced by the same identity on log
correspondences. This proves the theorem.
\end{proof}

\section{Proper generators and cohomological consequences}
\label{sec:consequences}

\begin{proof}[Proof of Theorem~\ref{thm:proper}]
By Theorem~\ref{thm:derived}, either \'etale motivic category
with the stated torsion coefficients is equivalent to
\(D(k_\et,R)\). The free sheaves represented by finite
\'etale \(k\)-schemes generate that unbounded derived
category. Indeed the corresponding presheaves are projective
generators, and exact sheafification preserves their
generation under shifts, cones, and coproducts.
Under \(\eta_\sharp\) and \(I_t\) they give the motives
of those same finite \'etale schemes with trivial log
structure. These schemes are proper.

Consequently the proper part is the entire logarithmic
category, so the left comparison in the theorem is an
equivalence. The right comparison is the restriction
of the interior equivalence \(J_t\), and is therefore
an equivalence as well. The canonical compatibility
with change of topology was proved with
Theorem~\ref{thm:derived}.
\end{proof}

\begin{corollary}[The groups in {\cite[Remark~8.6.2]{BPO}}]\label{cor:constant}
Let \(n\) be invertible in \(k\), \(X\in\lSm/k\), and \(r\in\Z\).
For \(t=\detop,\letop\) there are canonical identifications
\[
 \Hom_{\logDM_t^\eff(k,\Z/n)}
       (M_t(X),\one_t[r])
       \simeq H^r_\et(X^\circ,\Z/n).
\]
The change-of-topology map is the identity under these
identifications.
\end{corollary}

\begin{proof}
The equivalence \(J_t\) takes \(M_t(X)\) to the ordinary
motive of \(X^\circ\) and the unit to the ordinary constant
coefficient object. With invertible-order torsion that
constant sheaf is already \(\A^1\)-local
\cite[Proposition~3.1.7]{CD}. Derived evaluation gives
the displayed cohomology. The functors are compatible
with change of topology, so the map is the identity.
In the dividing \'etale category \(\one_t\) is the
image of the constant sheaf after projection localization.
\end{proof}

\begin{corollary}[Sections and boundary restriction]\label{cor:boundary-sections}
For \(F\) as in Theorem~\ref{thm:sheaf}, put \(G=\eta^*F\).
Then
\[
 F(X)\simeq H^0_\et(X^\circ,G|_{X^\circ}).
\]
For a union \(D'\) of components of a strict normal-crossing
boundary \(D\), the restriction
\[
 F(Y,D)\longrightarrow
 F(Y-D',D|_{Y-D'})
\]
is an isomorphism.
\end{corollary}

\begin{proof}
For bounded torsion over a perfect field this is the
section calculation for the pulled-back sheaf
\cite[Proposition~8.5.6 and its proof]{BPO}.
Section~\ref{sec:perfection} transports it across
perfection. For elementwise torsion take the filtered
colimit over the bounded torsion subgroups. Sections
commute with these colimits on the small logarithmic site,
whose finite-refinement property is \cite[Proposition~3.2.5]{BPO}, and on
the quasi-compact ordinary interior
\cite[Tag~03Q5]{Stacks}.
Both terms in the restriction formula are sections
on the same interior \(Y-D\), and their comparison
is the identity there.
\end{proof}

For a separably closed field, Theorem~\ref{thm:sheaf}
therefore gives the constant sheaf with group \(F(k)\).
For an arbitrary field it retains the Galois action
encoded by the small-site sheaf \(G\).

\begin{remark}[Scope of the coefficient hypotheses]
Theorem~\ref{thm:topology} is a comparison of the two
logarithmic topologies with arbitrary coefficients.
The comparison with ordinary \'etale motives and
the proper-generation conclusion use the torsion
hypothesis in Theorems~\ref{thm:derived} and
\ref{thm:proper}. Strict invariance of a sheaf means
invariance in every cohomological degree
\cite[Conjecture~8.5.9]{BPO}. The corresponding question for
invariance of sections appears in \cite[Remark~8.5.10]{BPO}.
\end{remark}

\section{Projective-line localization vs.\ projective-spaces localization}
\label{sec:projective}

\subsection{Projection invariance for sheaves of spaces}

Take \(S\in\Sch\) and use products over \(S\). We give the proof for
sheaves of spaces first and then pass to spectra.

By a space we mean a homotopy type. A Nisnevich sheaf of spaces is
a presheaf satisfying descent, with the limit of the diagram of a
cover and its intersections taken in spaces. In particular, if
\(Y=Y_1\cup Y_2\) is a Zariski cover by two opens, then
\[
  F(Y)\simeq F(Y_1)\times^h_{F(Y_1\cap Y_2)} F(Y_2).
\]
Here \(\times^h\) denotes the homotopy pullback. A point of this
pullback consists of points on the two opens and a path between
their restrictions. A sheaf is \(\Pp^n\)-local if the projection
induces an equivalence
\(F(X)\to F(X\times\Pp^n_S)\) for every \(X\in\Sm/S\).

\begin{theorem}\label{thm:unstable-projective}
Every \(\Pp^1\)-local Nisnevich sheaf of spaces on \(\Sm/S\) is
\(\Pp^n\)-local for all \(n\geq 0\).
\end{theorem}

A pointed sheaf \(H\) is a sheaf with a chosen point of \(H(S)\),
which gives compatible basepoints by pullback. Write \(G=\Omega H\)
for its sheaf of based loop spaces. We use concatenation for
multiplication of loops and reversal for inversion. These operations
are natural, and translation by any loop is an equivalence with
inverse given up to homotopy by translation by its reverse. Write
\(\fib\) for the homotopy fiber, taken at the specified basepoint.

\begin{lemma}\label{lem:projective-line-loops}
Suppose \(H\) is pointed and \(\Pp^1\)-local. For any smooth
\(U/S\), put
\[
\begin{aligned}
  D(U)&=G(U), & A(U)&=G(U\times\A^1_S),\qquad\\
  T(U)&=G(U\times\Gm), & B(U)&=\fib(A(U)\xrightarrow{\mathrm{ev}_1}D(U)),
\end{aligned}
\]
where \(\Gm=\A^1_S-\{0\}\) and \(\mathrm{ev}_1\) is restriction
along the section with coordinate \(1\). Let
\(\alpha_+,\alpha_-:A(U)\to T(U)\) be the pullbacks along
\(t\mapsto t\) and \(t\mapsto t^{-1}\), respectively. There is
a natural equivalence of spaces
\[
  \Phi_U:A(U)\times B(U)\longrightarrow T(U),\qquad
  (a,b)\longmapsto \alpha_+(a)\alpha_-(b).
\]
It carries left multiplication on the first factor by \(a_0\)
to left multiplication by \(\alpha_+(a_0)\).
\end{lemma}

\begin{proof}
For this proof abbreviate the four spaces by \(D,A,T,B\). On the
two affine charts of \(\Pp^1_S\), the coordinates are reciprocal
on their overlap. Descent and \(\Pp^1\)-locality therefore identify
\(H(U)\) with
\[
 H(U\times\A^1_S)\times^h_{H(U\times\Gm)}H(U\times\A^1_S),
\]
by pullback from \(U\). A loop in \(H(U\times\Gm)\) determines
a point of this pullback with basepoints in both charts. Its
corresponding point of \(H(U)\) is in the basepoint component,
because evaluation at \(1\) on the first chart is a left inverse
to pullback from \(U\). A path in the pullback to its basepoint
gives chart loops whose restrictions multiply, with one reversed,
to the class of the given loop. Thus every component of \(T\)
has a representative \(\alpha_+(a)\alpha_-(a')\). Let
\(d=\mathrm{ev}_1(a')\), and let \(\bar d\in A\) be its pullback
from \(U\). Since \(\alpha_-(\bar d)=\alpha_+(\bar d)\), we may
replace the two factors by \(a\bar d\) and \(\bar d^{-1}a'\).
The second has evaluation homotopic to the unit and so defines
a point of \(B\). This proves surjectivity on components.

Inverting the second factor identifies the fiber of \(\Phi_U\)
over the unit with
\[
  (A\times^h_T A)\times^h_D *,
\]
where the map to \(D\) evaluates the second chart at \(1\).
Indeed a path from \(\alpha_+(a)\alpha_-(b)\) to the unit is
equivalent to a path from \(\alpha_+(a)\) to \(\alpha_-(b^{-1})\).
Looping the descent equivalence identifies \(A\times^h_T A\)
with \(D\), compatibly with this evaluation. The displayed fiber
is therefore contractible. Left translation in \(A\) by \(a_0\)
and right translation in \(B\) by \(b_0\) correspond under
\(\Phi_U\) to left translation by \(\alpha_+(a_0)\) and right
translation by \(\alpha_-(b_0)\) in \(T\). These equivalences
identify the unit fiber with the fiber at \(\Phi_U(a_0,b_0)\).
Since there is such a point in every component of \(T\), all
homotopy fibers are contractible. The asserted naturality and
compatibility with left multiplication follow from the formula.
\end{proof}

In \(\Pp^2_S\), let \(U_i=(x_i\ne 0)\) and \(W=U_1\cup U_2\).
Write \(t_i=x_i/x_0\) on \(U_0=\A^2_S\), and define
\[
  V=U_0\cap W,\qquad V_i=U_0\cap U_i=(t_i\ne 0),\qquad
  V_{12}=V_1\cap V_2,\qquad W_{12}=U_1\cap U_2.
\]
We write \(\A^1_i\) and \({\Gm}_i\) for the affine line and its
open complement of zero in coordinate \(t_i\).

\begin{lemma}\label{lem:projective-plane-loops}
Let \(H\) be a pointed \(\Pp^1\)-local Nisnevich sheaf of spaces and
\(G=\Omega H\). Every class in \(\pi_0G(V)\) is the product of a
class from \(\pi_0G(U_0)\) followed by a class from \(\pi_0G(W)\).
\end{lemma}

\begin{proof}
Let
\[
 B_1(Y)=\fib\bigl(G(\A^1_1\times Y)\xrightarrow{\mathrm{ev}_1}G(Y)\bigr).
\]
This is itself the loop space of a pointed \(\Pp^1\)-local
Nisnevich sheaf, namely the analogous fiber formed with \(H\)
in place of \(G\). Descent and locality are preserved under
limits, which also commute with taking loops. Use
Lemma~\ref{lem:projective-line-loops} in coordinate \(1\) to
identify the diagram for the two opens of \(V\) with
\[
 G(U_0)\times B_1(\A^1_2)
   \longrightarrow G(V_2)\times B_1({\Gm}_2)
   \longleftarrow G(V_2).
\]
The left arrow is restriction on both factors by naturality.
The right arrow is inclusion of the first factor since
\(V_2=\A^1_1\times{\Gm}_2\) and the restriction to \(V_{12}\)
uses the coordinate \(t_1\). Its homotopy pullback is
\[
  G(V)\simeq G(U_0)\times
     \fib\bigl(B_1(\A^1_2)\to B_1({\Gm}_2)\bigr).
\]
In forming this pullback the path in the first component can
be eliminated because the right arrow is the identity there.
Apply Lemma~\ref{lem:projective-line-loops} to \(B_1\) in
coordinate \(2\). With
\[
  B_{12}=\fib\bigl(B_1(\A^1_2)\xrightarrow{\mathrm{ev}_1}B_1(S)\bigr),
\]
it gives \(B_1({\Gm}_2)\simeq B_1(\A^1_2)\times B_{12}\), with
restriction along \(t_2\) given by inclusion of the first factor.
We obtain
\begin{equation}\label{eq:projective-plane-split}
  G(V)\simeq G(U_0)\times\Omega B_{12}.
\end{equation}
Restriction from \(U_0\) is inclusion of the first factor.
Moreover, multiplication by a loop from \(G(U_0)\) on the left
acts by left multiplication on the first factor and leaves the
second factor fixed. In the pullback diagram this follows from
the last assertion of Lemma~\ref{lem:projective-line-loops}: on
each chart this loop multiplies the first factor by its restriction,
and on the overlap the path in the second component is unchanged.

A loop at the unit in \(G(V_{12})\) gives a point in \(G(V)\)
by taking units on both opens and using that loop as the gluing
path. Under \eqref{eq:projective-plane-split} its first component
is the unit and its second component is the loop obtained by
projecting first to \(B_1({\Gm}_2)\) and then to \(B_{12}\).
Here we take gluing paths from the restriction of the first
chart to that of the second, and identify the fiber of an
inclusion of the first factor with loops in the second by
using that same path orientation.

The regular map
\[
  W_{12}\longrightarrow\A^2_S,\qquad
  [x_0:x_1:x_2]\longmapsto(x_0/x_1,x_0/x_2)
\]
restricts on \(V_{12}\) to \((t_1^{-1},t_2^{-1})\).
Pullback along this map sends the composite
\(B_{12}\to B_1(\A^1_2)\to G(U_0)\) to \(G(W_{12})\).
On restriction to \(V_{12}\), the inversions in both
coordinates place this map in the second factor of each
application of Lemma~\ref{lem:projective-line-loops}. More
explicitly, the coordinate \(2\) pullback of
\(B_{12}\to B_1(\A^1_2)\) along \(t_2^{-1}\) is its inclusion
as the second factor of \(B_1({\Gm}_2)\). Then the pullback
along \(t_1^{-1}\) is the inclusion of \(B_1({\Gm}_2)\) as
the second factor of \(G(V_{12})\). Consequently the composite
to \(B_{12}\) under the two projections is the identity.
Loop the map \(B_{12}\to G(W_{12})\) and use the resulting
loops as gluing paths for the cover \(W=U_1\cup U_2\), with
units on both opens. Restricting these points of \(G(W)\) to
\(V\) gives the points \((1,k)\) of \eqref{eq:projective-plane-split},
for \(k\in\Omega B_{12}\). Multiplying on the left by points
from \(G(U_0)\) gives all components of \(G(V)\).
\end{proof}

\begin{proof}[Proof of Theorem~\ref{thm:unstable-projective}]
Let \(F\) be a \(\Pp^1\)-local Nisnevich sheaf of spaces.
For each \(n\), it suffices to prove the claim at \(X=S\) for
every such sheaf. Indeed, for any smooth \(X/S\) the sheaf
\(F_X(Y)=F(X\times_S Y)\) is \(\Pp^1\)-local by the hypothesis
on \(F\). It satisfies Nisnevich descent because taking the
product with \(X\) preserves Nisnevich covers and their
\v{C}ech nerves. The cases \(n=0,1\) hold since \(\Pp^0_S\cong S\)
and \(F\) is \(\Pp^1\)-local.

If \(F\) is \(\Pp^d\)-local and \(E\) is a vector bundle of rank
\(d+1\) over a smooth \(Y/S\), then the bundle of lines
\(\Pp_Y(E)\to Y\) induces an equivalence
\(F(Y)\to F(\Pp_Y(E))\). On every open of a trivializing
cover and on its intersections this is the pullback along a
product projection with \(\Pp^d_S\). Descent gives the
equivalence on \(Y\).

For \(n=2\), first consider components. Use the cover by
\(U_0\) and \(W\) from Lemma~\ref{lem:projective-plane-loops}.
Choose a point \(s\in F(\Pp^2_S)\) and
let \(c\in F(S)\) be its restriction along \([1:1:1]\). The open
\(W\) carries the projection \([x_0:x_1:x_2]\mapsto[x_1:x_2]\).
Let \(L\subset\mathcal O^2\) be the tautological line bundle on
\(\Pp^1_S\), whose fiber at a line is that line. The projective
bundle \(Q=\Pp_{\Pp^1_S}(\mathcal O\oplus L)\) maps to
\(\Pp^2_S\) by the inclusion \(\mathcal O\oplus L\subset
\mathcal O^3\). Over \(W\) it is an isomorphism, since the
projection of a line representing a point of \(W\) to the
last two coordinates determines the unique line in \(\mathcal O^2\)
above which it lies. By \(\Pp^1\)-locality and the projective
bundle observation, pullback is an equivalence \(F(S)\simeq F(Q)\).
Evaluating at the lift of \([1:1:1]\) is its inverse, so
\(s|_W\) is in the component of the pullback of \(c\).

On \(U_0\) the morphism
\[
  U_0\times\Pp^1_S\longrightarrow\Pp^2_S,\qquad
  (t_1,t_2,[u:v])\longmapsto[u:ut_1+v:ut_2]
\]
has restrictions at \([1:0]\) and \([0:1]\) equal to the chart
inclusion and the constant map with value \([0:1:0]\). The triple
on the right defines a morphism over every base: it generates
the same ideal locally as \(u,v\). Both endpoint evaluations
are inverse to pullback along \(U_0\times\Pp^1_S\to U_0\)
on \(F\), so the two restrictions of \(s\) are in the same
component. Since \([0:1:0]\) lies in \(W\), this proves that
\(s|_{U_0}\) is also in the component of the pullback of \(c\).

Use \(c\) as a basepoint for \(F\). After choosing paths from the
restrictions of \(s\) to the pullbacks of \(c\) on the two opens,
the remaining gluing datum is a loop in \(F(V)\) at \(c|_V\).
Changing the two paths by loops on their respective opens
changes the gluing class by left multiplication by the image
from \(\pi_1F(U_0,c)\) and right multiplication by the image from
\(\pi_1F(W,c)\), up to inversion in the second group. This is
the path description of the homotopy pullback for the cover.
By Lemma~\ref{lem:projective-plane-loops}, these two changes can
make the gluing class the unit. The resulting point lies in the
component of the pullback of \(c\) on \(\Pp^2_S\).

The section \([1:1:1]\) gives injectivity on components.
For each \(c\in F(S)\) and \(r\geq1\), apply the component
bijection just proved to the Nisnevich sheaf of \(r\)-fold based
loops of \(F\) at \(c\). This sheaf is \(\Pp^1\)-local since
loop spaces are defined by limits. The resulting component
bijections are the isomorphisms on \(\pi_r\) at every basepoint.
The component bijection and these isomorphisms give an
equivalence \(F(S)\to F(\Pp^2_S)\). Applying the same result to
\(F_X\) proves \(\Pp^2\)-locality.

For the induction step, let \(n\geq3\) and assume that every
\(\Pp^1\)-local Nisnevich sheaf of spaces is \(\Pp^d\)-local
for \(0\leq d<n\). Write
\(\mathcal O_S^{n+1}=A\oplus B\) with \(A\) and \(B\) free of
ranks \(2\) and \(n-1\). Let \(L_A\) and \(L_B\) be
the tautological lines on \(\Pp(A)\) and \(\Pp(B)\).
Form the schemes
\[
\begin{aligned}
  Q_A&=\Pp_{\Pp(B)}(A\oplus L_B),&
  Q_B&=\Pp_{\Pp(A)}(L_A\oplus B),\\
  Q_{AB}&=\Pp_{\Pp(A)\times\Pp(B)}(L_A\oplus L_B).
\end{aligned}
\]
Bundles in this display are pulled back to the indicated bases.
The first two schemes map to \(\Pp^n_S=\Pp(A\oplus B)\) by
inclusion of the indicated bundles, and forgetting one of the
chosen lines maps \(Q_{AB}\) to each of them. This gives a square
\[
\begin{matrix}
 Q_{AB}&\longrightarrow&Q_B\\
 \downarrow&&\downarrow\\
 Q_A&\longrightarrow&\Pp^n_S.
\end{matrix}
\]
The opens complementary to \(\Pp(A)\) and to \(\Pp(B)\) cover
\(\Pp^n_S\). Over the complement of \(\Pp(A)\), projection to
\(B\) determines the chosen line in \(B\) uniquely. On this
open both horizontal arrows of the square are isomorphisms.
Indeed the projection of the line in \(A\oplus B\) to \(B\)
has some coordinate invertible locally, so its image is a
line subbundle there and defines the inverse morphisms.
Over the complement of \(\Pp(B)\), the same reasoning with
projection to \(A\) shows that both vertical arrows are
isomorphisms. Applying \(F\) to the square restricted to
either open or their intersection therefore gives a homotopy
pullback square. The inverse images of these opens cover each
scheme in the display. Descent for those covers and the
commutation of limits with limits give
\[
  F(\Pp^n_S)\simeq F(Q_A)\times^h_{F(Q_{AB})}F(Q_B).
\]
The bases of the three projective bundles are \(\Pp^{n-2}_S\),
\(\Pp^1_S\), and \(\Pp^1_S\times\Pp^{n-2}_S\).
Their respective fibers are projective spaces of dimensions
\(2\), \(n-1\), and \(1\). All these dimensions are less
than \(n\). By induction and Zariski descent for projective
bundles, the structure maps identify \(F(Q_A)\), \(F(Q_B)\),
and \(F(Q_{AB})\) with \(F(S)\). These identifications commute
with the arrows since all arrows are over \(S\). The pullback
formula therefore makes \(F(S)\to F(\Pp^n_S)\) an equivalence.
Applying this conclusion to \(F_X\) completes the induction and
gives the required equivalences
\[
  F(X)\xrightarrow{\sim}F(X\times_S\Pp^n_S)
  \qquad (X\in\Sm/S,\ n\geq0).\qedhere
\]
\end{proof}

\subsection{The stable comparison}

\begin{corollary}[Answer to {\cite[Question~4.0.4]{BPOh}}]\label{thm:projective}
The canonical comparison
\[
  \SH^{\mathrm{eff},\Pp^1}(S)\longrightarrow
  \SH^{\mathrm{eff},\Pp^\bullet}(S)
\]
of \cite[Question~4.0.4]{BPOh} is a symmetric monoidal equivalence
for every \(S\) in \(\mathrm{Sch}\) of \cite{BPOh}, natural with
respect to pullback in \(S\).
\end{corollary}

\begin{proof}
In \cite[Definition~4.0.3]{BPOh}, the two unstable localizations
are the full subcategories of Nisnevich sheaves of spaces
local for \(\Pp^1\) and for all \(\Pp^n\), respectively.
Theorem~\ref{thm:unstable-projective} identifies these full
subcategories. Pointing and then forming \(S^1\)-spectra
therefore gives the same category on both sides. Explicitly,
an \(S^1\)-spectrum is described by a sequence of pointed
objects with each object identified with the based loops of
the next. The condition on projective spaces is the same at
each level of the sequence by the theorem.

The two unstable localizations use products of the indicated
projections with every smooth scheme. Hence the localizations
are compatible with the cartesian monoidal structure, since
smooth representable presheaves generate sheaves under
colimits and the localized product preserves colimits in
each variable. The stable comparison is the symmetric monoidal
functor induced by the comparison of these localizations,
and an equivalence on underlying categories is a symmetric
monoidal equivalence in this situation. Base change sends each
smooth scheme, Nisnevich cover, and projective-space projection
to the corresponding one over the new base. The localization
comparisons thus commute with pullback by their universal
properties.
\end{proof}

\subsection{The marked projective line}

Use the site \(\Pair_S\) and sheaf category \(\mathcal C_S\)
of Section~\ref{sec:spectra-definitions}. Write
\(\mathrm{SmAdm}\) for composites of admissible blow-ups
along smooth centers. The notation \([\bx^{-1}]\) imposes
all product projections with \(\bx\).

Binda, Park, and \O stv\ae r establish the symmetric monoidal
equivalences \cite[Theorem~1.3.1]{BPOh}
\[
 \log\SH^{\mathrm{eff}}(S)
 \simeq \mathcal C_S[\bx^{-1},\mathrm{SmAdm}^{-1}]
 \simeq \mathcal C_S[(\Pp^\bullet,\Pp^{\bullet-1})^{-1}].
\]
The first term is their effective logarithmic category, stabilized
with respect to \(S^1\). The second term uses the marked line and
admissible blow-ups \cite[Corollary~3.4.3]{BPOh}. The third uses
all marked projective spaces \cite[Corollary~3.4.7]{BPOh}.

\begin{theorem}[Four presentations of effective logarithmic spectra]
\label{thm:four-log-presentations}
For every quasi-compact, quasi-separated scheme \(S\) with trivial
log structure, there are canonical symmetric monoidal equivalences
\[
\begin{aligned}
 \log\SH^{\mathrm{eff}}(S)
 &\simeq \mathcal C_S[\bx^{-1},\mathrm{SmAdm}^{-1}]\\
 &\simeq \mathcal C_S[(\Pp^\bullet,\Pp^{\bullet-1})^{-1}]\\
 &\simeq \mathcal C_S[\bx^{-1}].
\end{aligned}
\]
In particular, a strict Nisnevich sheaf of spectra on these pairs
that is invariant under product with \(\bx\) is invariant under
product with every \(P_n\) and under admissible blow-ups along smooth
centers.
\end{theorem}

\begin{proof}
Work in \(\mathcal C_S[\bx^{-1}]\). For a smooth pair \(Y\),
write \(M(Y)\) for the image of its represented suspension spectrum
with a disjoint basepoint, and write \(\one=M(S)\) for the tensor
unit. The localization is symmetric monoidal. Indeed, the
spectra represented by smooth pairs generate under colimits and
shifts, and tensoring a projection in the defining class with
the spectrum of a smooth pair gives another such projection.
The tensor product preserves colimits in each variable, and
\(M(Y\times_S Y')\simeq M(Y)\otimes M(Y')\). Strict Zariski descent
expresses \(M(Y)\) as the colimit of the represented spectra of
the nonempty intersections of any finite strict open cover of \(Y\).
For two opens this is the pushout for the open union, and the
general finite case follows by taking successive open unions.
Here and below all colimits and fibers are formed in the
\(\infty\)-categorical sense. We use \(\Sigma\) for suspension,
and the cofiber of a map \(C\to D\) is the pushout \(D\amalg_C 0\).

\smallskip
\noindent\emph{The two charts of the marked line.}
Set
\[
 A=M(\A^1_S,\varnothing),\qquad
 B=M(\A^1_S,\{0\}),\qquad
 T=M(\Gm,\varnothing).
\]
Here \(\Gm=\A^1_S-\{0\}\). Let \(f:T\to A\) and \(g:T\to B\)
be the chart maps given by
\(t\mapsto t\) and \(t\mapsto t^{-1}\), respectively. Since
the pushout for the two charts is \(M(\bx)\simeq\one\), we have
a cofiber sequence
\[
 T\xrightarrow{(f,-g)}A\oplus B
   \xrightarrow{(\epsilon_A,\epsilon_B)}\one,
\]
where \(\epsilon_A,\epsilon_B\) are the structure maps. The section
\(1\) lies in the unmarked open of the second chart and splits
\(\epsilon_B\). Write \(K=\fib(\epsilon_B)\) and \(p_K:B\to K\)
for the projection in the resulting splitting \(B\simeq\one\oplus K\).
The projection from the fiber of \((\epsilon_A,\epsilon_B)\) to
\(A\oplus K\) is an equivalence, since the \(\one\)-component of
\(B\) is determined by the \(A\)-component. Taking the negative of
the \(K\)-coordinate gives the equivalence
\((f,p_K g):T\simeq A\oplus K\). Thus the arrow \(f\) is identified
with projection onto \(A\), and the inclusion \(s:K\to T\) of the
second summand satisfies \(p_K g s\simeq\mathrm{id}_K\).

\smallskip
\noindent\emph{The intersections of the projective charts.}
Fix \(n\geq 2\). In \(P_n\) let \(U_i\) be the open pair
\(x_i\ne 0\) with inherited boundary. Write
\(I=\{1,\ldots,n\}\), \(W=\bigcup_{i\in I}U_i\), and
\(V=U_0\cap W\). The pair \(U_0\) is the unmarked affine \(n\)-space
with coordinates \(t_i=x_i/x_0\). For nonempty \(J\subseteq I\),
put \(W_J=\bigcap_{i\in J}U_i\) and \(V_J=U_0\cap W_J\).
The spectra \(M(V_J)\) for nonempty \(J\) compute \(M(V)\)
by descent. The diagram has arrows from larger to smaller subsets.
Adjoining the empty subset with value \(M(U_0)=A^{\otimes n}\)
gives the \(n\)-dimensional cube formed by tensoring \(n\) copies
of the arrow \(f:T\to A\). In the \(i\)-th coordinate the value
is \(T\) for \(i\in J\) and \(A\) otherwise.

For any diagram \(D\) on the nonempty subsets, there is a
natural map
\[
 \colim_{\varnothing\ne J\subseteq I} D(J)
 \longrightarrow \Sigma^{n-1}D(I).
\]
It is induced by projecting to the diagram supported on \(I\)
alone, where the value is \(D(I)\). This projection is natural
because that diagram is the right Kan extension from \(I\).
The value at a proper subset is the limit over the empty category
of arrows from that subset to \(I\), hence is zero. In the
simplicial formula for the colimit of the projected diagram,
the chains beginning at \(I\) give the cone on the subdivision of
the boundary of the simplex on \(I\). The face omitting \(I\) has
value zero. Collapsing the boundary therefore gives the
\((n-1)\)-fold suspension of \(D(I)\).

The equivalence of arrows \(f\simeq(A\oplus K\to A)\) splits the
augmented cube for \(V\) as a direct sum indexed by subsets
\(H\subseteq I\). The summand for \(H\) has \(K\) in the
coordinates in \(H\) and \(A\) in all others, at exactly those
vertices \(J\) containing \(H\). Between its nonzero vertices the
maps are identities. For \(H=\varnothing\), the diagram is
constant and its indexing category has the initial object \(I\).
Thus its colimit is \(A^{\otimes n}\) and its augmentation is the
identity. For \(H=I\), the diagram is supported on the deepest
intersection and gives \(\Sigma^{n-1}K^{\otimes n}\). If \(H\) is
nonempty and proper, choose a direction \(i\notin H\). The maps
in that direction are equivalences, including between zero
vertices. The cofiber of the augmentation from the colimit
over nonempty subsets to the empty vertex is computed by taking
cofibers one direction at a time. This follows by induction from
the definition of cofiber as a pushout. Taking direction \(i\) first
shows that this cofiber is zero. The empty vertex is also zero
because \(H\ne\varnothing\), so this summand contributes zero.
Consequently there is a splitting
\[
 M(V)\simeq A^{\otimes n}\oplus\Sigma^{n-1}K^{\otimes n}
\]
with first projection the map to \(M(U_0)\). Under the natural map
from \(M(V)\) to \(\Sigma^{n-1}M(V_I)\) constructed above, the
inclusion of the second summand becomes
\(\Sigma^{n-1}s^{\otimes n}\).

Write \(a:\Sigma^{n-1}K^{\otimes n}\to M(W)\) for the second
summand followed by \(M(V)\to M(W)\). On \(W_I\) the morphism
\[
 h:W_I\longrightarrow(\A^1_S,\{0\})^n,\qquad
 [x_0:\cdots:x_n]\longmapsto (x_0/x_1,\ldots,x_0/x_n)
\]
respects the log structures. Indeed, the boundary on \(W_I\) has
equation \(u=x_0/x_1\), and each target coordinate pulls back to
\(u\) times the unit \(x_1/x_i\). Thus each target log generator
pulls back to the source log generator multiplied by the
corresponding unit. On \(V_I\) the map is
\((t_1,\ldots,t_n)\mapsto(t_1^{-1},\ldots,t_n^{-1})\).
Apply the natural map to the deepest intersection for \(W\) and
then apply \(h\) and \(p_K^{\otimes n}\), with the same suspension
throughout. This gives a left inverse to \(a\). In fact, the map
\(V_J\to W_J\) is natural in \(J\), so the composite with \(a\)
is \(\Sigma^{n-1}(p_K g s)^{\otimes n}\), the identity.

The two-open cover \(P_n=U_0\cup W\) exhibits \(M(P_n)\)
as the cofiber of the map from \(M(V)\) to
\(M(U_0)\oplus M(W)\) with one component negated. Choose the
negative sign on the \(M(U_0)\)-component. In the displayed
splitting of \(M(V)\), the first component is minus the
identity on \(A^{\otimes n}\) and zero on the second summand.
An automorphism of the target fixing \(M(W)\) removes the map
from the first summand to \(M(W)\). Thus
\(M(P_n)\simeq\cofib(a)\), compatibly with the map from
\(M(W)\). Since \(a\) has a left inverse, it is the inclusion of
a direct summand in the stable category, and \(M(W)\to M(P_n)\)
has a right inverse.

\smallskip
\noindent\emph{A projective bundle with one marked section.}
Let \(L\subset\mathcal O^n\) be the tautological line bundle on
the unmarked \(\Pp^{n-1}_S\) with coordinates \([x_1:\cdots:x_n]\).
Over this base take the bundle of lines
\(Q=\Pp_{\Pp^{n-1}_S}(\mathcal O\oplus L)\) and mark the
section \(\Pp(L)\). The inclusion
\(\mathcal O\oplus L\subset\mathcal O^{n+1}\) gives a map
\(Q\to P_n\). Concretely, the first coordinate defines the
boundary on \(P_n\), and its pullback vanishes precisely along
\(\Pp(L)\) as a Cartier divisor. The map is therefore strict.
It is an isomorphism above \(W\). On \(W\), projection to the
last \(n\) coordinates gives the morphism to \(\Pp^{n-1}_S\)
that recovers the point of the bundle base.
The section \(\Pp(\mathcal O)\) is disjoint from the marked section
and maps to \(o=[1:0:\cdots:0]\).

A Zariski trivialization of \(L\) identifies the projection of
the marked bundle \(Q\) with a product projection with factor
\(\bx\). By descent, \(M(Q)\to M(\Pp^{n-1}_S)\) is an
equivalence, with inverse induced by the unmarked section.
It follows that \(M(W)\to M(Q)\to M(P_n)\) factors through
the map \(M(o):\one\to M(P_n)\). Since the composite has a
right inverse, \(M(o)\) has a right inverse. The structure
map \(M(P_n)\to\one\) is already its left inverse. This proves
that the structure map is an equivalence.

For \(n=1\) the equivalence is the defining relation for
\(\bx\). For every smooth pair \(X\) and every \(n\),
tensoring with \(M(X)\) now proves that the projection from
\(X\times_S P_n\) is inverted in \(\mathcal C_S[\bx^{-1}]\).
The canonical further localization at the marked projective
spaces is therefore an equivalence. Combining it with the
three original presentations proves the theorem.
\end{proof}

\subsection{The embedding of projective-line local spectra}

Let \(\rho:\Pair_S\to\Sm/S\) send \(Y\) to its underlying scheme
\(\underline Y\). The right adjoint comparison of
\cite[Construction~4.0.6]{BPOh} is given by precomposition with
\(\rho\) in the marked projective-space presentation. We can now
express it using locality for the ordinary projective line.

\begin{corollary}\label{cor:projective-embedding}
For every qcqs ordinary scheme \(S\), the functor
\[
 \rho^*:\SH^{\eff,\Pp^1}(S)\longrightarrow\logSH^\eff(S),
 \qquad (\rho^*F)(Y)=F(\underline Y),
\]
is fully faithful. Under Corollary~\ref{thm:projective}, it identifies
with the right adjoint in \cite[Construction~4.0.6]{BPOh}.
\end{corollary}

\begin{proof}
On the sites, \(\rho\) is left adjoint to
\(\lambda(X)=(X,\triv)\). A morphism from \(Y\) to
\(\lambda(X)\) is uniquely determined by an ordinary map
\(\underline Y\to X\), and \(\rho\lambda=\operatorname{id}\).
On presheaves, precomposition with \(\rho\) is therefore the
left Kan extension along \(\lambda^\op\), with right adjoint
restriction along \(\lambda\). The unit of this adjunction is
an equivalence, so \(\rho^*\) is fully faithful on presheaves.

Taking underlying schemes carries strict Nisnevich covers and
their \v{C}ech nerves to ordinary Nisnevich covers and their
\v{C}ech nerves. Moreover,
\(\underline{Y\times_S\bx_S}=\underline Y\times_S\Pp^1_S\).
Thus \(\rho^*F\) satisfies strict Nisnevich descent and marked-line
invariance for every \(\Pp^1\)-local Nisnevich sheaf of spectra
\(F\). Theorem~\ref{thm:four-log-presentations} identifies this
full subcategory of presheaves with \(\logSH^\eff(S)\).
Full faithfulness follows from the presheaf comparison. The right
adjoint in \cite[Construction~4.0.6]{BPOh} has the same evaluation
on represented objects, which identifies the two functors.
\end{proof}

The underlying-scheme functor thus gives the projective-line
analogue of the embedding in \cite[Construction~4.0.8]{BPOh}.
At the \(S^1\)-stable level it reads
\[
 \bigl(\Pp^1\text{-local homotopy theory}\bigr)
 \ \subset\
 \bigl(\bx\text{-local homotopy theory}\bigr).
\]

\section{Continuity of the categories and of mapping spectra}\label{sec:continuity}

A strict Nisnevich
distinguished square is a cartesian square
\[
\begin{matrix}
 V\times_X U&\longrightarrow&V\\
 \downarrow&&\downarrow p\\
 U&\longrightarrow&X
\end{matrix}
\]
in which \(U\to X\) is a quasi-compact open immersion,
\(p\) is strict \'etale and finitely presented, and \(p\)
is an isomorphism over the reduced closed complement of \(U\).
Use the strict Nisnevich topology of Section~\ref{sec:spectra-definitions}.

For \(X\in\Pair_S\), let \(h_S(X)\) be its representable
presheaf of spaces. Begin with presheaves of spectra on \(\Pair_S\).
Impose the empty-object relation and the excision maps
\begin{equation}\label{eq:nisrelation}
 \Sigma^\infty h_S(U)_+
 \amalg_{\Sigma^\infty h_S(U\times_XV)_+}
 \Sigma^\infty h_S(V)_+
 \longrightarrow\Sigma^\infty h_S(X)_+.
\end{equation}
The pushout in this formula is the pushout in spectra-valued
presheaves. Impose also the projection maps
\begin{equation}\label{eq:projective-relations}
 \Sigma^\infty h_S(X\times(\Pp^n,\Pp^{n-1}))_+
   \longrightarrow\Sigma^\infty h_S(X)_+,
 \qquad n\geq1.
\end{equation}
Finally make the pointed projective-line suspension invertible.

Its object before inversion is \(T_S\) from
Section~\ref{sec:spectra-definitions}. This describes
\(\logSH(S)\) by \cite[Corollary~3.4.7]{BPOh}.

\begin{lemma}\label{lem:nis}
For a qcqs base, strict Nisnevich sheaves of spectra on \(\Pair_S\)
are exactly the presheaves sending the empty object to zero and
every strict Nisnevich distinguished square to a pullback.
\end{lemma}
\begin{proof}
The finite-presentation Nisnevich excision criterion is
\cite[Proposition~C.5(4)]{HoyoisDuality}. The same argument applies
to strict covers here: a strict \'etale object over a pair \(X\)
is an ordinary \'etale scheme over \(\underline X\) with the
pulled-back log structure. Every Nisnevich refinement and
distinguished-square construction in that proof remains in
\(\Pair_S\). The construction can therefore be performed on these
strict \'etale slices, for each \(X\). Applying the space-valued
criterion to \(\Omega^\infty\Sigma^j F\) for every integer
\(j\) gives the assertion for spectra.
\end{proof}

This is descent for covers. Hypercompletion imposes descent for all
hypercovers as well. The finite presentation above uses ordinary
sheaf descent, as in the strict Nisnevich presentation of
\cite[Corollary~3.4.7]{BPOh}.

An object \(C\) of a presentable stable category is
\emph{compact} if \(\Map(C,-)\) preserves filtered colimits.
The category is compactly generated when its compact objects
generate under colimits. A finite cell object is built from a
specified family using finitely many shifts and finite colimits.
A retract is a direct summand. The idempotent completion formally
adjoins all retracts.
An exact functor between stable categories is one preserving finite
colimits, equivalently finite limits.

\begin{lemma}\label{lem:compact}
For every qcqs \(S\), the objects
\(\Sigma^j T_S^{\otimes m}\otimes M_S(X)\), with
\(X\in\Pair_S\) and \(j,m\in\Z\), form a family of compact
generators of \(\logSH(S)\). Pullback preserves compact objects.
\end{lemma}
\begin{proof}
Representable spectral presheaves are compact, by evaluation.
Both sides of \eqref{eq:nisrelation} and
\eqref{eq:projective-relations} are finite colimits of compact
representables. By Lemma~\ref{lem:nis}, all the required
locality conditions use only finite limits of evaluations.
Filtered colimits of spectra commute with finite limits, so
local objects are closed under filtered colimits computed
objectwise. The localized representables are therefore compact
and generate the localized category.
Only strict Nisnevich excision was used in this argument. In particular,
it applies to spectra themselves, including sphere-spectrum
coefficients. It does not use an eventual coconnectivity hypothesis
from the discussion of other topologies and coefficient categories in
\cite[Remark~2.6.10]{BPOh}.

Tensoring a representable with \(T_S\) is a cofiber of two
representables, since products of pairs remain pairs. Thus this
tensor functor preserves compact objects. Formal inversion of
\(T_S\) gives the stated compact generators. The cyclic permutation
on \(T_S^{\otimes 3}\) is homotopic to the identity
\cite[Remark~3.3.5]{BPOh}, so the inversion is computed by
the colimit of the sequence of compactly generated categories with
transition functor \(-\otimes T_S\).
Every compact object is a retract of a finite cell object in these
generators. Pullback sends each displayed generator to the same
generator over the new base, and preserves finite colimits and
retracts. Hence it preserves compact objects.
\end{proof}

Let \(S=\varprojlim_i S_i\) be a cofiltered inverse limit,
with affine transition maps and qcqs stages. Fix a stage \(S_0\)
and restrict to the cofinal subsystem over that stage.
Write \(p_i:S\to S_i\) for the projections and
\(p_{i0}:S_i\to S_0\) for the structure maps.
Only systems of affine schemes are needed for the final argument.

\begin{lemma}\label{lem:finite-data}
Pairs in \(\Pair_S\), their morphisms, and finite diagrams of
such morphisms descend to a stage \(S_i\), with equalities
holding at some later stage. The same is true of the relations
\eqref{eq:nisrelation} and \eqref{eq:projective-relations}.
\end{lemma}
\begin{proof}
The underlying schemes and smooth morphisms descend by finite
presentation. Use a finite chart cover of a pair. Its monoids
have finitely many generators and relations, and its chart
functions, units, log morphisms, and gluing data descend after
increasing the stage. Relative normal-crossing coordinates also
descend. These are the approximation statements used in
\cite[\S3.1, especially Propositions~3.1.4 and~3.1.6]{BPOh}.

For a distinguished square, its quasi-compact open immersion and
\'etale map descend. Choose a finitely presented closed subscheme
with the complement as its support. The \'etale map is an
isomorphism over its reduction by the distinguished-square
condition. The reduction induces an equivalence on small
\'etale sites, so the map is already an isomorphism over the
chosen closed subscheme \cite[Tag~04DZ]{Stacks}. The closed
subscheme, the isomorphism, and strictness descend as well, so
the square is distinguished at a later stage. The projective pairs in
\eqref{eq:projective-relations} and their projections are defined
universally over \(\Z\), so they also descend.
\end{proof}

\begin{proposition}\label{prop:continuity}
For a compact \(C_0\in\logSH(S_0)\) and an arbitrary
\(F_0\in\logSH(S_0)\), put
\(C_i=p_{i0}^*C_0\), \(F_i=p_{i0}^*F_0\),
\(C=p_0^*C_0\), and \(F=p_0^*F_0\). Then
\begin{equation}\label{eq:continuity}
 \operatorname*{colim}_i
 \Map_{\logSH(S_i)}(C_i,F_i)
 \xrightarrow{\ \sim\ }
 \Map_{\logSH(S)}(C,F).
\end{equation}
In particular, this holds when \(C_0=\one_{S_0}\).
\end{proposition}
\begin{proof}
We include the categorical step, because \(F_0\) need not
be compact. Write \(\mathcal A^\omega\) for the full category
of compact objects in a compactly generated stable category
\(\mathcal A\). Let
\(\mathcal P_i=\PSh(\Pair_{S_i},\Sp)\) and
\(\mathcal P=\PSh(\Pair_S,\Sp)\).
Lemma~\ref{lem:finite-data} identifies \(\Pair_S\) with the
filtered colimit of the categories \(\Pair_{S_i}\), including
their morphisms and finite diagrams. Consequently
\begin{equation}\label{eq:presheaf-compacts}
 \mathcal P^\omega\simeq
 \left(\operatorname*{colim}_i
                    \mathcal P_i^\omega\right)^{\idem}.
\end{equation}
Here and below the colimit on small categories uses the induced exact
functors, and the superscript \(\idem\) denotes idempotent
completion. To verify \eqref{eq:presheaf-compacts}, recall that
compact presheaves are retracts of finite cell objects in the spectral
representables. Maps between representables are the suspension
spectra of the sets of morphisms of pairs. These sets are filtered
colimits of their stagewise sets. Finite cell constructions use finite
limits and colimits of mapping spectra, so the same assertion holds
for their mapping spectra. Finally adjoining retracts gives
\eqref{eq:presheaf-compacts}.

At each stage, first quotient this small stable category by the
cofibers of the empty-object, excision, and projection relations,
and then adjoin retracts. A quotient here has the following universal
meaning: exact functors out of it are precisely exact functors out of
the original category that send the specified cofibers to zero.
Localization by maps between compact objects has this compact-object
description \cite[\S\S5.1--5.2]{BGT}. Every relation over
\(S\) descends by Lemma~\ref{lem:finite-data}. The universal
description of the quotient therefore shows that forming these
quotients commutes with the filtered colimit. Idempotent completion
does so as a construction in idempotent-complete stable categories.
Finally invert \(T\). By Lemma~\ref{lem:compact}, tensoring
with \(T\) preserves compact objects. Its formal inversion is a
colimit of these small categories, and colimits over the suspension
index commute with the filtered colimit over the bases. We obtain
\begin{equation}\label{eq:compact-continuity}
 \logSH(S)^\omega\simeq
 \left(\operatorname*{colim}_i
                    \logSH(S_i)^\omega\right)^{\idem}.
\end{equation}
Mapping spectra between objects coming from compact objects at a
fixed stage are the filtered colimits of their stagewise mapping
spectra: this holds in the filtered colimit of small stable categories,
and idempotent completion is fully faithful on the objects already
present. This proves the compact-target case of
\eqref{eq:continuity}.

For completeness, reduce the arbitrary target to this compact
case. A compactly generated stable category is the category of
filtered colimits of its compact objects, also called their
\(\Ind\)-completion \cite[\S2.4]{BGT}. Choose
\(F_0\simeq\operatorname*{colim}_a F_{0a}\), with a filtered
diagram of compact objects. Pullback preserves colimits.
Using compactness of \(C\) and the compact case just proved,
\begin{align*}
 \Map(C,F)
 &\simeq\operatorname*{colim}_a\Map(C,p_0^*F_{0a})\\
 &\simeq\operatorname*{colim}_a\operatorname*{colim}_i
       \Map(C_i,p_{i0}^*F_{0a})\\
 &\simeq\operatorname*{colim}_i\Map(C_i,F_i).
\end{align*}
The last equality uses compactness of \(C_i\), and the
interchange of two colimits. This proves
\eqref{eq:continuity} for mapping spectra.
\end{proof}

\begin{corollary}\label{cor:cats-continuity}
For \(S=\lim_i S_i\) as above, the comparison functors
\[
 \colim_i\mathcal N_{S_i}\longrightarrow\mathcal N_S,\qquad
 \colim_i\mathcal L_{S_i}\longrightarrow\mathcal L_S
\]
are equivalences. The colimits are in the category of
presentable categories and colimit-preserving functors,
with pullback as the transition functors.
\end{corollary}
\begin{proof}
The ordinary smooth sites and their finite Nisnevich
excision relations descend to finite stages by the same
arguments as in Lemma~\ref{lem:finite-data}, with no log
structures. In the effective log case, use the quotient
on the compact presheaves in the proof of
Proposition~\ref{prop:continuity} before inverting \(T_S\).
In the ordinary case the quotient uses only the empty
scheme and Nisnevich excision. In either case, the finite
presentations identify the compact categories at \(S\)
with the idempotent completions of their filtered colimits.
Completion under filtered colimits identifies small
idempotent-complete stable categories and exact functors
with compactly generated stable categories and
colimit-preserving functors preserving compact objects.
Their colimits here are also colimits among presentable
categories and colimit-preserving functors: such functors
out of each completion are specified by exact functors
out of its compact category. This proves the
asserted equivalences.
\end{proof}

\section{Park's toric construction}\label{sec:toric}

Write \(\ParkSch\) for separated Noetherian schemes of finite
Krull dimension, Park's base class in \cite{ParkI}. Over these
bases one may use the separated smooth schemes and smooth pairs
as sites, since they form Zariski bases stable under products with
the projective pairs. Park's presentation of the effective
logarithmic category agrees with \(\mathcal L_S\)
\cite[Theorem~2.11]{ParkI}\cite[Corollary~3.4.7]{BPOh}.

For a Nisnevich sheaf \(H\) of spectra on the ordinary smooth
site, the logarithm \(\Log H\) is the extension of
\cite[Definition~5.25]{ParkI}. The following formula suffices
for the comparisons in this paper.

We describe the indexing schemes and maps for \(X\times\bx_S^r\).
A fan in a lattice is a finite collection of rational polyhedral cones closed
under taking faces, with intersections that are common faces. A cone is smooth if the shortest nonzero
lattice vectors along its rays extend to a lattice basis. A subdivision refines a fan without changing
its support. The fan of \((\Pp^1)^{r+m}\) has cones spanned by subsets of
\(\{e_1,-e_1,\ldots,e_{r+m},-e_{r+m}\}\) containing at most one of \(e_i,-e_i\) in each coordinate. The
cone spanned by the positive generators describes the affine chart with coordinates taking values in
\(\A^1\).

Fix integers \(r,m\geq 0\). Let \(I_{r,m}\) be the category of \(r\)-standard fans, whose objects are smooth subdivisions of
this fan satisfying the following two conditions.
\begin{enumerate}
\item If a ray has positive \(i\)-th coordinate for some \(1\leq i\leq r\), its primitive generator is \(e_i\) for that \(i\).
\item The cone generated by \(e_{r+1},\ldots,e_{r+m}\) belongs to the subdivision. For \(m=0\) this is the zero cone.
\end{enumerate}
There is a morphism from a finer subdivision to a coarser one when it refines it. These are the fans of
the standard blow-ups of Park's monoid version of \(\bx^r\times\bx^m\), by
\cite[Proposition~4.4 and Example~4.5]{ParkI}. For each cone, the integral linear forms nonnegative on
it form an additive monoid. In Park's monoid scheme these monoids, with an absorbing element adjoined,
replace the affine coordinate rings and are glued using the localizations for common faces. The
realization over \(\Z\) takes the monoid algebra of the original lattice points as the coordinate ring.
Gluing the resulting affine schemes gives the toric scheme \(Y_\Sigma\) over \(\Z\). A subdivision in
\(I_{r,m}\) is smooth and complete, meaning that its support is the whole ambient real vector space, so
\(Y_\Sigma\) is smooth proper of relative dimension \(r+m\) over \(\Z\). Write
\(Y_{\Sigma,S}=Y_\Sigma\times_{\Spec\Z}S\). Each \(Y_\Sigma\) has its toric map to
\((\Pp^1_{\Z})^{r+m}\), and hence a structure map to \(\Spec\Z\).

For clarity, the simplicial variable in the logarithm comes from the category \(\Delta\) of finite
nonempty ordered sets \([m]=\{0,\ldots,m\}\) and order-preserving maps. For \(\alpha\colon[m']\to[m]\),
the map on the last interval factors from \(\bx^{m'}\) to \(\bx^m\) is given as follows. In
output coordinate \(j\), let \(i\) be the least index with \(\alpha(i)\geq j\), or \(m'+1\) if no such
index exists. Use the input coordinate \(i\) when \(1\leq i\leq m'\), the constant \(1\) when \(i=0\),
and the constant \(0\) when \(i=m'+1\). The first \(r\) factors are fixed. This is the cosimplicial
interval convention in \cite[p.~88]{MorelVoevodsky} used in \cite[Construction~3.7]{ParkI}.

To obtain maps between the colimits over subdivisions, one uses maps of the underlying
compactifications lifting these coordinate maps. For a specified target subdivision, Park shows that
one can find such a lift after refining the source, and that lifts agree after common refinement. More
precisely, Park forms a category of pairs \((m,\Sigma)\) over \(\Delta\), with arrows induced by lifts
in the monoid schemes, and then formally inverts the refinement arrows at fixed \(m\). The left Kan
extension along the opposite projection is the functor on \(\Delta^\op\) obtained by taking colimits,
and its value at \([m]\) is the colimit over \(I_{r,m}^\op\). The common-refinement result ensures this
formula as a diagram, including its compositions \cite[Propositions~4.13--4.14 and
Definition~4.15]{ParkI}. This describes how the simplicial structure is obtained even when a coordinate
map does not lift between a given pair of subdivisions. The realizations over \(\Z\) of the
monoid-scheme lifts are ordinary morphisms between the displayed toric schemes. Their products with
\(X\) give the lifts used for \(X\times\bx_S^r\).

With these maps understood, Park's base change formula gives a natural equivalence
\begin{equation}\label{eq:park}
 \Log H(X\times\bx_S^r)
 \simeq
 \left|[m]\longmapsto
           \colim_{\Sigma\in I_{r,m}^\op} H(X\times_S Y_{\Sigma,S})\right|.
\end{equation}
Here \(|-|\) is geometric realization, the colimit of a simplicial diagram. The formula follows by
using the map to the monoid model of \(\bx^r\), which pulls back its log structure to the boundary
log structure on \(X\times\bx_S^r\). These toric subdivisions
suffice to compute the colimit among the standard blow-ups of the product over \(S\)
\cite[Proposition~5.14]{ParkI}. Thus \eqref{eq:park} is the formula in
\cite[Equation~(5.10)]{ParkI}. The canonical comparison from \(H(X)\) is induced by the
projections \(X\times_S Y_{\Sigma,S}\to X\), compatibly with all simplicial maps.

We write \(C_m^{(r)}(X,H)\) for the term inside the realization.
The natural augmentation from \(H(X)\) is an equivalence after
realization for every \(X\in\Sm/S\) and \(r\geq0\) precisely
when \(H\) is logarithmic in Park's sense
\cite[Definition~5.26]{ParkI}. The same terms also form a
cubical diagram. The ordered cube category \(\Box\) has objects
\(\{0,1\}^m\) and maps generated by inserting \(0,1\) and by
deleting coordinates. The filtration and the endpoint homotopy
in Section~\ref{sec:spectral} use these faces and degeneracies.
The first \(r\) coordinates in the
toric charts are fixed. The cubical face maps are given by
restriction to \(0\) or \(1\) in a varying coordinate,
while a degeneracy is pullback along a projection. A zero
face is the toric orbit closure corresponding to the ray
\(e_i\), with its fan obtained from the cones containing that
ray by passage to the quotient lattice. The one face is
given by the fan on the hyperplane where that coordinate
of the lattice is zero. The induced diagrams and their
maps up to common refinement are the ones of
\cite[Example~4.22 and Remark~7.8]{ParkI}.

\subsection{Cubical and simplicial realizations}\label{sec:cube-comparison}

We compare the two realizations for sheaves of spectra. First let
\(\mathcal E\) be a small category with finite products and an interval \(J\)
equipped with multiplication and endpoints \(0,1\)
satisfying the zero and unit identities. Let \(G\)
be a presheaf of spectra on \(\mathcal E\). The
presheaf
\[
             \operatorname{Cub}_{J}G(U)=
              \colim_{\Box^{\mathrm{op}}}G(U\times J^\bullet)
\]
is the localization of \(G\) at the projections
\(U\times J\to U\), as a presheaf of spectra. Here
is a verification of both parts of that statement.
The maps induced by multiplying the parameter \(J\)
in \(U\times J\) by an extra last cubical coordinate
give a homotopy from the identity to the map induced
by setting that parameter equal to zero after
cubical realization. For this, append a last
coordinate in the cube category. Each of its endpoints gives a
natural transformation from the appended diagram to the original
one, and hence the same map to its realization by the colimit
cocone. The multiplication defines a map to the appended diagram
because it commutes with the old cubical operations. This makes
\(\operatorname{Cub}_{J}G\) local.

For the local-equivalence assertion, set
\(G_n=G(-\times J^n)\). Pullback from the projection
and from the all-zero section exhibit \(G\) as a
retract of \(G_n\). The other composite is homotopic
to the identity after localization. Indeed multiplying
the \(n\) coordinates by a variable map into \(J\)
gives a map \(G_n\otimes \mathsf y(J)\to G_n\) with
the two desired endpoints. The tensor denotes the
pointwise action of a presheaf of spaces on a
presheaf of spectra, and \(\mathsf y(J)\) is the
represented presheaf of spaces. This is a valid
homotopy after localization because
\(K\otimes\mathsf y(J)\to K\) is a local equivalence
for every presheaf of spectra \(K\). It suffices
to check this on suspension spectra of represented
presheaves and their shifts, where it is the defining
projection relation. The maps \(G\to G_n\) are
therefore local equivalences in all degrees and are
compatible with the cubical structure. Passing to
the colimit gives the assertion, since the cube
category has contractible nerve. The simplex-based
singular construction for the same interval also
computes this localization in spectra
\cite[Construction~3.7]{ParkI}. Both identifications
are the ones induced by the unit of localization.
They give a natural equivalence of the two
realizations, compatible with the maps from \(G\).

Apply this observation to the diagrams defining the logarithm. For fixed \(X/S\) and an ordinary Zariski
sheaf of spectra \(H\), the rule
\[
       V\longmapsto H(X\times_S V_S)
\]
is a Zariski sheaf of spectra on the category of
smooth monoid schemes used in \cite[Section~4]{ParkI}.
The scheme \(V_S\) in this notation is the realization
over \(S\) of the ordinary monoid scheme \(V\).
Its restrictions to toric open covers are Zariski
covers after realization over \(S\), which explains
the sheaf assertion. The toric terms at \(\bx^r\)
are those of the logarithm and its cubical version
for this sheaf by the standard subdivision formula
\cite[Example~4.5 and Remark~7.8]{ParkI}.

Here the comparison of the two realizations is local
in the monoid-scheme parameter. We give the details
to distinguish it from an application of the
comparison stated for complexes in
\cite[Proposition~4.30]{ParkI}. On an affine
monoid-scheme chart of the form \(\A^P\times U\),
\(P\) is a monoid isomorphic to
\(\mathbb N^a\times\mathbb Z^b\), and \(\A^P\) has its
toric log structure. Its scheme realization over a
ring \(R\) is \(\Spec R[P]\) with the chart from \(P\).
The factor \(U\) is an ordinary
affine space or a product of an affine space with
a torus. Define the presheaf \(G_U\) on smooth log
monoid schemes by evaluating \(H\) on the product of
\(X\), \(U_S\), and the underlying scheme of the
argument. The termwise calculation of
\cite[Lemmas~4.25, 4.27, and~4.28]{ParkI} identifies
the standard subdivision colimit on this chart
with
\[
              (L_{\mathrm{Adm}}G_U)
                            (\A^P\times\bx^m).
\]
In this formula \(L_{\mathrm{Adm}}\) is presheaf
localization by admissible blow-ups, computed by
the filtered colimit over such blow-ups. The
termwise calculation works for spectra. It uses
cofinality of the toric modification systems
under products with the ordinary affine or
torus factor and equality of standard and
admissible blow-ups on the charts \(\A^P\).
The maps in the two diagrams are the interval
maps after admissible localization, by
\cite[Proposition~4.13 and Definitions~4.15 and~4.20]{ParkI}.
The interval in the cartesian category of smooth
log monoid schemes with admissible blow-ups
inverted is \(\bx\), with multiplication
induced by the toric resolution of ordinary
multiplication \cite[Proposition~3.10]{ParkI}.
The comparison for presheaves of spectra established
above therefore identifies the two realizations
on each affine chart.

These identifications are natural for restrictions
of the toric charts. In the cover of \(\bx^r\)
by products of its two affine charts, intersections
replace factors by ordinary tori. Such a factor
can be included either in \(U\) or in \(\A^P\)
as a factor with trivial log structure. The two
termwise formulas are identified by product with
that torus, since the admissible refinement fan
is unchanged by adding torus factors. The interval
comparison is compatible with these identifications
by naturality of its localization unit. The affine
identifications thus glue by finite Zariski descent.
For both diagrams this descent follows from the
same cofinality under restriction to opens, as
in \cite[Lemma~4.24 and Proposition~4.26]{ParkI},
and commutation of colimits of spectra with finite
limits. For a finite open cover these limits can be taken
by successive two-open unions and their intersections.
This proves the natural equivalence
between the cubical realization and the simplicial
toric realization for spectra, respecting the
augmentation from \(H(X)\).

\section{Restriction of logarithmic spectra}\label{sec:spectral}

Giving an ordinary smooth scheme the trivial log structure
induces a functor \(\lambda_\sharp:\mathcal N_S\to\mathcal L_S\)
on the categories in Section~\ref{sec:spectra-definitions}.
It sends the ordinary represented spectrum of \(Y\) to
\(M(Y,\triv)\). Its right adjoint \(\lambda^*\) evaluates
a logarithmic spectrum on ordinary smooth schemes. These
are the adjoints of \cite[Construction~4.0.16]{BPOh}.

\begin{theorem}[Answer to {\cite[Question~4.0.17]{BPOh}}]\label{thm:q417-full-scaling}
For every quasi-compact, quasi-separated ordinary scheme \(S\), the
right adjoint \(\lambda^*:\mathcal L_S\to\mathcal N_S\) is fully
faithful. In particular it induces equivalences on mapping spectra.
\end{theorem}

\subsection{Square scaling on toric schemes}

A lattice \(N\) determines the split torus with character
lattice \(\operatorname{Hom}(N,\mathbb Z)\). A cocharacter
\(v\in N\) is a homomorphism \(\Gm\to T_N\). It acts on
a toric scheme by translation in its dense torus, extended to
the whole scheme. Here \(V_B\) is the base change to the ordinary
base \(B\) of the split toric scheme over \(\mathbb Z\) defined by
the fan.

\begin{lemma}\label{lem:q417-toric-scaling}
Let \(V\) be the smooth complete toric scheme associated
to a fan in \(N\), with a union \(D\) of invariant boundary components
marked. Let \(B\) be an ordinary qcqs base and let \(a\) be a unit on \(B\).
For \(v\in N\), let \(u\) be the endomorphism of \(M(V_B,D_B)\)
induced by the torus action of \(v(a^2)\). Then \(1-u\) is nilpotent
in the homotopy category of \(\mathcal L_B\), meaning that a
positive power of \(1-u\) is zero.
\end{lemma}

\begin{proof}

For \(V=(\Pp^1)^n\) with empty marking, translation by \(v(a^2)\)
is a product of the square scalings \((a^{v_i})^2\).
The conclusion follows from Lemma~\ref{lem:unmarked-scaling}.

Suppose \(A\to D\to C\) is a cofiber sequence with a
compatible endomorphism \(w\) on its terms. If \(w^e=0\)
on \(C\) and \(w^f=0\) on \(A\), then \(w^{e+f}=0\)
on \(D\). Indeed \(w^e:D\to D\) factors through \(A\)
by exactness of the mapping sequence, and postcomposing
with \(w^f\) kills the factorization. Rotation of the
cofiber sequence shows that nilpotence on any two terms
implies nilpotence on the third.

Use toric weak factorization to connect the fan of \(V\)
to the product fan, on the common lattice \(N\)
\cite[Theorem~4.7.1]{WlodarczykWeak}.
Each step is a blow-up or blow-down in a smooth invariant
center \(Z\) of codimension at least two. Its exceptional
divisor \(E\) and its center \(Z\) are smooth complete
toric schemes of smaller dimension, component by component if the
center has several components. Their torus actions
are obtained from the original torus by the homomorphisms
to the dense tori of the orbit closures. In particular
\(v(a^2)\) acts through the value at \(a^2\) of the induced
cocharacters. All these constructions are made on fans,
so the corresponding ordinary schemes, their blow-up
squares, and the actions are defined over \(\mathbb Z\)
and can be base changed to \(B\).

By induction on dimension, \(1-u\) is nilpotent on the
motives of \(E\) and \(Z\). The ordinary blow-up square
is cocartesian in \(\mathcal L_B\) by
\cite[Theorem~3.3.7]{BPOh}. For a blow-up \(W\to V'\)
it gives an equivariant sequence
\[
 M(E_B)\longrightarrow M(W_B)\oplus M(Z_B)
                         \longrightarrow M(V'_B).
\]
Thus nilpotence on \(M(W_B)\) implies nilpotence on
\(M(V'_B)\). Conversely, nilpotence on \(M(V'_B)\)
implies nilpotence on the middle term and hence on
its invariant summand \(M(W_B)\). The product case
starts the factorization, and dimension zero starts
the dimension induction. This proves the claim when the
marking is empty.

For a nonempty marking, remove one marked component \(Z\) and
write \(D=D'+Z\). The logarithmic Gysin sequence for the strict
immersion of \(Z\), with the boundary induced by \(D'\), identifies
the cofiber of \(M(V_B,D_B)\to M(V_B,D'_B)\) with the logarithmic
Thom motive of the normal line bundle \(N_{Z_B/V_B}\).
Indeed, blowing up a Cartier divisor gives the original scheme,
and adding the exceptional divisor adds the marking \(Z_B\).
The Gysin sequence is natural for automorphisms of the pair
\cite[Theorem~3.2.19 and Remark~3.2.20]{BPOh}. The projective model
of the Thom motive is the cofiber of the infinity inclusion
\[
 M\bigl(\Pp(N_{Z_B/V_B})\bigr)\longrightarrow
 M\bigl(\Pp(N_{Z_B/V_B}\oplus\mathcal O)\bigr),
\]
where both schemes carry the markings pulled back from the
induced boundary on \(Z_B\) \cite[Proposition~3.2.13]{BPOh}.
The projective bundles are smooth complete split toric schemes.
The equivariant line bundle gives the torus action on their
dense tori. More explicitly, over the dense torus of \(Z_B\)
the normal line is equivariantly trivial after choosing a
splitting of the quotient of split tori. The kernel acts on the
line by a character, and the resulting action on the dense
projective-bundle torus is by translation through a homomorphism
of tori. Thus the induced automorphisms are again given by
cocharacters evaluated at \(a^2\). Two distinct invariant divisors
of a smooth toric scheme either do not meet or meet in a single
orbit closure. Hence the induced boundary on \(Z_B\), and its
pullbacks to the projective bundles, have at most as many
components as \(D'\). Induction on the number of components and
the cofiber-sequence argument prove the claim.
\end{proof}

\subsection{The filtration and its realized quotients}

Take \(S\in\ParkSch\), \(X\in\Sm/S\), \(r\geq0\), and
\(F\in\mathcal L_S\). Write \(C_\bullet\)
for the cubical diagram of \(C_m^{(r)}(X,\lambda^*F)\). By cubical
realization we mean the colimit over \(\Box^{\op}\). A mark means
that the indicated divisors are recorded by the log structure.

For a model \(V=Y_{\Sigma,S}\) in degree \(m\) and a subset
\(J\subseteq [r+m]=\{1,\ldots,r+m\}\), let \(D_J(V)\) consist of the
toric divisors lying over the infinity hyperplanes in the
coordinates indexed by \(J\), each with multiplicity one.
Use the log structure given by the local equations of those
divisors. In particular the construction uses these equations
also when \(S\) is nonreduced. Set
\[
 A_m(J)=\colim_\Sigma F\bigl(X\times_S(V,D_J(V))\bigr),
 \qquad
 (L_k C)_m=\colim_{|J|\le k} A_m(J),
\]
where the second colimit uses inclusions of subsets. Adding
divisors gives a morphism of pairs from the more marked pair
to the less marked pair, so evaluation on \(F\) gives maps in
the direction of those inclusions.

A cubical face deletes its coordinate from the set of marks.
A degeneracy inserts an unmarked coordinate. On the toric models
the face maps are restrictions at \(0\) or \(1\), and the
degeneracy maps are projections after refinement. The support of
the pullback of a marked divisor is contained in the corresponding
marked boundary in the source. For the deleted coordinate the
face lies entirely off the infinity hyperplane. The toric chart
formulas pull the marked equations back to monomials in the source
boundary equations times units, which supplies the maps of log
structures. Thus these operations define the cubical diagrams
\(L_k C\). The first stage is \(L_0 C=C\), and
\(\colim_k L_kC\) is the diagram with all coordinates marked.

For \(k>0\), put \(B=\cofib(L_{k-1}C\to L_k C)\) in cubical
spectra. In each degree,
\begin{equation}\label{eq:q417-marking-quotient}
 B_m\simeq\bigoplus_{|J|=k} B_m(J),
 \qquad
 B_m(J)=\cofib\left(\colim_{I\subsetneq J}A_m(I)\longrightarrow A_m(J)\right).
\end{equation}

The following description specifies the functoriality, including
the zero maps. For each cube let \(J\) range over the subsets
of size at most \(k\), and let \(J\mapsto uJ\) denote transport
of labels under a cubical operation \(u\). The sets of labels
with their subset inclusions form a diagram of posets over
the opposite of the cube category. Its category of pairs
\((m,J)\) has arrows over \(u:m\to m'\) allowed when
\(uJ\subseteq J'\). Let
\[
 P_m(J)=\colim_{I\subseteq J,\ |I|<k} A_m(I),\qquad
 Q_m(J)=\cofib(P_m(J)\to A_m(J)).
\]
The map \(P\to A\) is natural in this category of pairs.
In fact \(P\) is the left Kan extension of \(A\) from the
subsets of size less than \(k\). To compute it at \((m,J)\),
an arrow \((m',I')\to(m,J)\) with \(|I'|<k\) factors through
\((m,uI')\) with \(uI'\subseteq J\) and \(|uI'|<k\).
Thus the colimit can be taken in the fiber over \(m\).
It follows that \(Q_m(J)\) vanishes for \(|J|<k\)
and that its colimit over the fiber is \(B_m\).

For clarity, restrict \(Q\) to the category with the same
objects \((m,J)\) and only arrows over \(u\) with target
\((m',uJ)\). This is the category of pairs for the
diagram of \emph{sets} of labels. The inclusion of the
objects with \(|J|=k\) is a sieve, meaning that there
are no arrows from a smaller subset to a subset of
size \(k\). Thus \(Q\) on this category is the right
Kan extension by zero of its values and maps with
\(|J|=|uJ|=k\). At a smaller subset the indexing category
for this right Kan extension is empty, so its value is
the zero spectrum. This proves that any natural comparison
on the retained labels extends canonically to the diagrams
with the zero maps for deleted marks. Taking the direct
sum over the discrete set of labels in each degree gives
\eqref{eq:q417-marking-quotient} naturally in the cube.
Indeed, in the original poset the only nonzero values are
the maximal subsets of size \(k\), and the colimit from
those objects alone is the same.

An \(r\)-standard fan is recovered from its restriction to
\(x_a\leq0\) for \(1\leq a\leq r\) by adjoining \(e_a\) to cones
in \(x_a=0\), in every combination of such indices
\cite[Proposition~3.1.3]{ParkII}. In particular it contains the full
positive orthant.

\begin{lemma}\label{lem:q417-scaling-refinement}
Let \(V\) be the toric scheme of an \(r\)-standard fan \(\Sigma\)
in dimension \(n\geq1\), and fix \(1\leq i\leq n\).
Given finitely many positive integers \(j\), there is an \(r\)-standard
refinement \(W\to V\times\Pp^1_t\), unchanged over \(t\ne0\), on
whose finite \(t\)-chart each map
\[
 \phi_j:(z_1,\ldots,z_n,t)\longmapsto
        (z_1,\ldots,t^{2j}z_i,\ldots,z_n)
\]
extends to \(V\). The strict \(t=0\) face maps into the
orbit closure corresponding to \(e_i\). On the \(t=1\) face the
map is the identity.
\end{lemma}

\begin{proof}
On lattices the maps are \(g_j(x,s)=x+2js e_i\), where \(s\)
is the \(t\)-coordinate. Restrict the product fan to \(x_a\leq0\)
for \(a\leq r\). For \(s\geq0\), intersect its cones with inverse
images of cones of \(\Sigma\) under the finitely many \(g_j\).
Along \(s=0\) the original fan is unchanged, since \(g_j\) is
the identity there. The cone spanned by
\(e_{r+1},\ldots,e_n,e_t\) is unchanged as well, since its image
lies in the full positive orthant of \(\Sigma\). Use smooth toric
subdivision relative to these unchanged smooth cones
\cite[Theorem~11.1.9]{CLSToric}. For \(s\leq0\) keep the product
fan in the region \(x_a\leq0\) for \(a\leq r\). The two fans glue
along \(s=0\). Adjoin the fixed positive rays by the reconstruction
rule. The result is smooth and refines the original product fan
because the original product satisfies the same rule.

For the fan-map property after adjoining the positive rays, let
\(\tau\) be a cone in \(s\geq0\) before adjoining any \(e_a\).
If \(e_a\) is to be adjoined, \(x_a=0\) on \(\tau\) and
\(g_j(\tau)\) has nonnegative \(a\)-th coordinate. Its smallest
containing target cone either contains \(e_a\) or lies in
\(x_a=0\). In the second case the reconstruction rule allows
\(e_a\) to be adjoined. Since \(g_j(e_a)=e_a\), this gives the
required containment for the enlarged cone. The \(t=0\) assertion
uses \(g_j(e_t)=2j e_i\). The \(t=1\) assertion follows from the
unchanged product over \(t\ne0\) and the formula for \(\phi_j\).
\end{proof}

\begin{proposition}\label{prop:q417-all-cubical}
The augmentation \(F(X)\to |C_\bullet^{(r)}(X,\lambda^*F)|\)
is an equivalence for every \(r\geq0\).
\end{proposition}

\begin{proof}
Fix \(k>0\) and use \(B\) from \eqref{eq:q417-marking-quotient}.
Within the diagram \(m\mapsto B_{m+1}\) with a new last varying
coordinate \(t\), take the summands with \(t\) unmarked. Denote
this diagram, with only the old cubical operations, by \(B^{\mathrm{app}}\).
The complement of these summands is preserved by the old
operations as well, since these operations either keep \(t\)
or induce the zero map on the summand. The corresponding diagrams
formed by restricting the toric models to \(t\ne\infty\), \(t\ne0\),
and their intersection are denoted \(B^{\mathrm{app}}_0\),
\(B^{\mathrm{app}}_1\), and \(B^{\mathrm{app}}_*\). Zariski descent
gives a pullback square on the two opens for the evaluation of
\(F\) at each marking. Finite colimits, filtered colimits, and
cubical realization are exact in spectra. Applying these
operations to the squares gives
\begin{equation}\label{eq:q417-appended-gluing}
 |B^{\mathrm{app}}|\simeq
 |B^{\mathrm{app}}_0|\times_{|B^{\mathrm{app}}_*|}
 |B^{\mathrm{app}}_1|.
\end{equation}

On the summand \(B_m(J)\), choose \(i=\min J\) in the coordinate
order, including the fixed coordinates. Lemma~\ref{lem:q417-scaling-refinement}
gives a pullback \(b_j:B\to B^{\mathrm{app}}_0\) for every \(j\ge1\).
On the finite \(t\)-chart a ray \((x,s)\) has \(s\geq0\).
If its image under \(g_j\) has negative \(i\)-th coordinate,
then \(x_i<0\). The map leaves the other coordinates unchanged.
Thus a source divisor mapped into a marked infinity component lies
over a marked infinity hyperplane of the source. The toric formulas
pull the boundary equations back to monomials in the source boundary
equations times units. This gives maps of pairs for all
\(I\subseteq J\) with \(t\) unmarked. An old face
deleting a mark gives the zero map on both quotient diagrams by
\eqref{eq:q417-marking-quotient}. On every surviving face, and on
every degeneracy, the coordinate rule carries \(i\) to the same
chosen coordinate, so it commutes with \(\phi_j\).

For the compatibility after taking refinements, one may regard
each system of models as a pro-object, that is, a cofiltered
diagram in the category of pairs. For a fixed list \(J\), one
uses the cofiltered diagram of the whole finite diagram of
\(I\)-markings. A morphism to any specified target model is
represented on a sufficiently fine source model. Two
representatives of the same coordinate formula agree after
a common refinement. These comparisons can be made over \(\Z\)
before base change and product with \(X\). The formulas agree on the
common affine locus, which is schematically dense, and the targets
are separated. In toric charts the log maps are determined by the
pullback of the marked monomials. For the zero face one uses the
strict transform of that affine face, with its orbit-closure
fan. Its coordinate formulas likewise agree on a schematically
dense affine locus. Thus compositions on the retained labels
are equal as pro-morphisms, before evaluating the presheaf and
taking filtered colimits. The construction of
\eqref{eq:q417-marking-quotient} also handles compositions that
delete a mark. This yields the asserted maps of cubical spectra
and their endpoint and chart restrictions.

At \(t=0\) the image of \(\phi_j\) misses the \(i\)-th infinity
divisor. With the evaluation arrows written from the \(V\)-diagram
to the \(t=0\) diagram, the arrow at \(I\) factors through
evaluation on the face with the labels \(I\setminus\{i\}\).
This factorization is compatible with inclusions of subsets
of \(J\). The total cofiber of such a diagram is the cofiber
of the map from the colimit over proper subsets to the value
at \(J\). It can be computed by taking cofibers successively
along the coordinates. The intermediate diagram has identity
maps in the \(i\)-direction, so its total cofiber is zero. The \(t=0\)
restriction of \(b_j\) is therefore zero. At \(t=1\), the
restriction of \(b_j\) is the identity.

Write \(\widehat B\) for \(B\) with \(X\) replaced by
\(X\times\Gm\), and let \(t\) be the coordinate on this extra
factor. The overlap diagram \(B^{\mathrm{app}}_*\) identifies
with \(\widehat B\) as follows. For an appended
\(r\)-standard fan \(\Lambda\), restricting to \(t\in\Gm\)
retains precisely its cones in the hyperplane \(s=0\), where
\(s\) is the last lattice coordinate. These cones form a smooth
\(r\)-standard fan in the old lattice. The model on this open
is the product of its toric scheme with \({\Gm}_t\), and each
submarking is pulled back from the old factor since \(t\)
is unmarked. Every old fan \(\Sigma\) occurs on the overlap
of the product model \(\Sigma\times\operatorname{Fan}(\Pp^1)\).
An appended fan and any such product fan have a common
\(r\)-standard refinement. For each old fan, the appended
fans whose restriction refines it therefore form a nonempty
filtered category when ordered by pullback. This proves
cofinality of restriction to \(s=0\) on opposite posets. The colimits on
the overlap therefore agree for each submarking. Restriction
to the overlap commutes with the old face and projection
maps, so the identifications give
\(B^{\mathrm{app}}_*\simeq\widehat B\) as cubical diagrams.

On each summand, let \(u^*\) be pullback by scaling
\(z_i\) by \(t^2\). This is compatible with the cubical maps
for the same reason as \(b_j\). The map \(1-u^*\) acts locally
nilpotently on each homotopy group of each \(\widehat B_m\).
In fact, at a single toric model
Lemma~\ref{lem:q417-toric-scaling} applies to every marking
with the same coordinate scaling. One can use the base
\(X\times\Gm\), since evaluating \(F\) on smooth pairs over
this base, viewed over \(S\), preserves the local relations.
Nilpotence passes to the finite cofiber in
\eqref{eq:q417-marking-quotient}. Indeed, the colimit over
the proper subsets can be built by adjoining the subsets in
order of cardinality, with each step a pushout along the
colimits over their own proper subsets. This uses finitely
many finite sums and cofiber sequences, to which the
nilpotence argument applies. Each homotopy class in the
filtered colimit comes from a single stage.
The telescope of \(1-u^*\), meaning the colimit of consecutive
copies of the diagram under this map, is consequently zero in
each cubical degree. Realization commutes with that colimit.
It follows that \(1-u^*\) acts locally nilpotently on
\(\pi_q|\widehat B|\) for every integer \(q\).

Let \(x\in\pi_q|B|\) and let \(\widehat x\) be its pullback
to \(\pi_q|\widehat B|\). For some \(d\ge1\) one has
\((1-u^*)^d\widehat x=0\), so
\[
  \widehat x
   =\sum_{j=1}^d(-1)^{j+1}\binom dj (u^*)^j\widehat x .
\]
Use \(\sum_{j=1}^d(-1)^{j+1}\binom dj b_j(x)\) on
\(|B^{\mathrm{app}}_0|\) and the pullback of \(x\) by the
projection on \(|B^{\mathrm{app}}_1|\). The product
identification on the overlap takes the restriction of
\(b_j(x)\) to \((u^*)^j\widehat x\) and the restriction of the
projection class to \(\widehat x\). The displayed identity
therefore gives equality of the two classes on
\(|B^{\mathrm{app}}_*|\simeq|\widehat B|\). The homotopy exact sequence of
\eqref{eq:q417-appended-gluing} gives a class in
\(\pi_q|B^{\mathrm{app}}|\) with these restrictions. Its
\(t=0\) restriction is zero and its \(t=1\) restriction is \(x\).

The two endpoint maps \(|B^{\mathrm{app}}|\to |B|\) are
homotopic. To see this, include \(B^{\mathrm{app}}\) in the
diagram \(m\mapsto B_{m+1}\) obtained by appending the last
coordinate. The endpoints give natural transformations from
the appended diagram to \(B\). After passing to colimits each
agrees with the colimit cocone on the appended diagram.
Consequently \(x=0\). This holds for every \(q\), so \(|B|=0\).

We have shown that \(|L_{k-1}C|\to |L_kC|\) is an equivalence
for each \(k>0\). The filtration is increasing and every term
in degree \(m\) has reached the fully marked diagram by
\(k=r+m\). Taking the increasing colimit and realizing therefore
identifies \(|C|\) with the fully marked realization. Its
terms are naturally equivalent to \(F(X)\). This is the
degreewise assertion in the proof of
\cite[Proposition~8.3]{ParkI}, applied to \(X\times\bx^r\).
On affine log charts \cite[Proposition~8.2]{ParkI} provides
cofinal refinements by admissible smooth blow-ups and one uses
interval invariance. The degreewise calculation glues by
Zariski descent, as in the proof of
\cite[Proposition~5.24]{ParkI}. The equivalences come from
projection to \(X\), so they identify the augmentation with
the constant diagram at \(F(X)\). The cube category has
contractible nerve since the zero-dimensional cube is terminal.
This proves the proposition.
\end{proof}

\begin{proof}[Proof of Theorem~\ref{thm:q417-full-scaling}]
For \(S\in\ParkSch\) and \(F\in\mathcal L_S\), put \(H=\lambda^*F\).
Proposition~\ref{prop:q417-all-cubical} and
Section~\ref{sec:cube-comparison} show that
\[
 H(X)\xrightarrow{\sim}(\Log H)(X\times\bx^r)
 \qquad (X\in\Sm/S,\ r\geq0).
\]
This is the logarithmic condition of \cite[Definition~5.26]{ParkI}.
For a sheaf satisfying it, Park identifies \(\Log H\) with
\(\lambda_\sharp H\), compatibly with the map from the sheaf on
pairs given by \(H\) of the underlying scheme
\cite[Corollary~8.9]{ParkI}. The reconstruction equivalence
\(\Log\lambda^*F\simeq F\) is induced by pullback along the maps
\((Y,D)\to(Y,\varnothing)\) \cite[Theorem~8.5]{ParkI}. The composite
\(\lambda_\sharp H\simeq\Log H\simeq F\) is therefore the counit
\(\lambda_\sharp\lambda^*F\to F\). Since the counit is an equivalence,
\(\lambda^*\) is fully faithful over \(S\).

For any qcqs \(S\), open restrictions to quasi-compact affine
opens jointly detect equivalences by strict Nisnevich descent.
The pullback of a sheaf along such an open immersion is
restriction on the sites. It commutes with \(\lambda^*\)
by evaluation and with \(\lambda_\sharp\) by its values
on represented spectra. Thus one can check the counit
over affine bases. Express an affine base as an inverse
limit of the spectra of its finitely generated
\(\Z\)-subalgebras, which belong to \(\ParkSch\).
By Corollary~\ref{cor:cats-continuity}, \(\lambda_\sharp\)
is the colimit of the left adjoints at the stages, since
they agree on ordinary represented spectra.
Taking right adjoints identifies \(\lambda^*\) with the
limit of the right adjoints. Limits of presentable categories
with right adjoints are computed in categories, and the
mapping spaces of the limit are the limits of mapping spaces.
Stagewise full faithfulness therefore gives full faithfulness
over \(S\). Since \(\lambda^*\) is exact, it also gives
equivalences on mapping spectra.
\end{proof}

\subsection{The essential image on ordinary schemes}
\label{sec:q418-essential-image}

Binda, Park, and \O stv\ae r ask which Nisnevich sheaves of spectra
on smooth schemes are restrictions of effective logarithmic spectra
\cite[Remark~4.0.18]{BPOh}. They seek a test on a sheaf \(H\) that
uses only ordinary schemes and the maps induced by their morphisms. By
Theorem~\ref{thm:q417-full-scaling}, the essential image is already a
full subcategory of \(\mathcal N_S\). A criterion for membership in
this subcategory therefore describes the effective logarithmic
category entirely in terms of ordinary sheaves.

Fix \(r,m\geq0\), put \(n=r+m\), and set
\(Q=(\Pp^1_{\Z})^n\), with coordinates \(z_1,\ldots,z_n\) on the
standard affine open. Use smooth toric modifications \(p:V\to Q\).
Here these are maps defined by smooth subdivisions of the product
fan, so \(V\) is a smooth split toric scheme over \(\Z\) and \(p\)
is proper, equivariant, and an isomorphism on the dense torus.
Impose the following two conditions on \(p\).
\begin{enumerate}
\item The map is an isomorphism over \((\A^1_{\Z})^n\).
\item For each \(a\leq r\), put \(U_a=\{z_a\ne\infty\}\) and
\(V_{a,1}=p^{-1}(\{z_a=1\})\). Under the identification
\(U_a=\A^1_{\Z}\times(\Pp^1_{\Z})^{n-1}\) with \(z_a\) the first
coordinate, the restriction of \(p\) is isomorphic to the product
of \(\mathrm{id}_{\A^1_{\Z}}\) and
\(V_{a,1}\to(\Pp^1_{\Z})^{n-1}\).
\end{enumerate}
These modifications form a category with arrows over \(Q\) from a
refinement to the coarser modification. Evaluate \(H(X\times_S V_S)\)
on each model and take the colimit along the resulting pullbacks.
The first \(r\) coordinates stay fixed. In a varying coordinate,
restriction to the strict transform of \(z_i=0\) or \(z_i=1\)
gives the faces. These strict transforms are the scheme-theoretic
closures of the corresponding coordinate sections on
\((\A^1_{\Z})^n\). The degeneracies use pullback along the coordinate
projections, which lift to morphisms after refinement. The construction
in Section~\ref{sec:toric} gives compatibility after common
refinement.

To identify this diagram with \(C_\bullet^{(r)}(X,H)\), consider the
half-space \(x_a\geq0\) in the lattice of the fan. The open
\(p^{-1}(U_a)\) is the toric scheme of the cones in that half-space.
The preimage of \(\{z_a\ne 0,\infty\}\) is the toric scheme of
the cones in \(x_a=0\). It is the product of \(\Gm\) in coordinate
\(a\) with the toric scheme of this fan in the remaining lattice.
In particular, that latter scheme is \(V_{a,1}\). A product
isomorphism in (2) agrees with the canonical identification on the
dense torus because it is over \(U_a\). It is therefore torus
equivariant: the two action maps to the separated target agree
on the product of the torus with the dense torus of the source.
Consequently (2) says exactly that the fan on the half-space
is obtained by adjoining \(e_a\) to the cones in \(x_a=0\).

For the converse from the ray condition, suppose the only ray with
positive \(a\)-th coordinate is generated by \(e_a\). For each cone
\(\tau\) in \(x_a=0\), choose a point in its relative interior and
approach it by small positive displacements along \(e_a\).
Since the fan is finite and complete, one cone contains infinitely
many of these points, so it contains their limit and hence all of
\(\tau\). That cone has \(e_a\) as a ray. Smooth cones are
simplicial, so \(\operatorname{Cone}(\tau,e_a)\) is a face of it.
Each cone meeting \(x_a>0\) is spanned by \(e_a\) and a face in
\(x_a=0\). This proves the product description.
Condition (1) preserves the full positive orthant. Under (2) this
is equivalent to preserving \(\operatorname{Cone}(e_{r+1},\ldots,e_n)\),
since one can adjoin the first \(r\) positive rays successively.
Thus these modifications give exactly the \(r\)-standard subdivisions.
The strict transforms and projections induce their cubical maps
from Section~\ref{sec:toric}
\cite[Examples~4.5 and~4.22]{ParkI}\cite[Proposition~3.1.3]{ParkII}.

\begin{corollary}[Answer to {\cite[Remark~4.0.18]{BPOh}}]\label{cor:q418-essential-image}
Let \(S\) be a quasi-compact quasi-separated ordinary scheme.
A sheaf \(H\in\mathcal N_S\) is in the essential image of
\(\lambda^*\) if and only if the maps induced by projection
\begin{equation}\label{eq:q418-ordinary-test}
 H(X)\longrightarrow \bigl|C_\bullet^{(r)}(X,H)\bigr|
\end{equation}
are equivalences for every \(X\in\Sm/S\) and \(r\geq0\).
The realization is computed by the ordinary scheme diagrams above.
\end{corollary}

\begin{proof}
First assume \(S\in\ParkSch\). For \(H=\lambda^*F\),
Proposition~\ref{prop:q417-all-cubical} gives
\eqref{eq:q418-ordinary-test}. Conversely,
Section~\ref{sec:cube-comparison} identifies this map with
\(H(X)\to(\Log H)(X\times\bx^r)\). Thus the tests say that
\(H\) is logarithmic in the sense of
\cite[Definition~5.26]{ParkI}. By \cite[Proposition~5.28]{ParkI},
\(\Log H\) belongs to \(\mathcal L_S\), and the tests for \(r=0\)
give \(H\simeq\lambda^*(\Log H)\). Park identifies \(\Log H\) with
\(\lambda_\sharp H\) in \cite[Corollary~8.9]{ParkI}.

Now take \(S=\Spec A=\varprojlim_i S_i\) with the affine stages
used in the proof of Theorem~\ref{thm:q417-full-scaling}. Under the
limit description by right adjoints, the projection of \(H\) to
stage \(i\) is \(H_i(Y)=H(Y\times_{S_i}S)\) for \(Y\in\Sm/S_i\).
This follows by adjunction on represented spectra. The
modifications and coordinate maps are defined over \(\Z\), so
the diagram for \(H_i\) at \(Y\) is the diagram for \(H\) at
\(Y\times_{S_i}S\). Every smooth test scheme over \(S\) descends
to a stage by finite presentation. The tests over \(S\) are
therefore equivalent to the tests at all stages.

The transition right adjoints commute with restriction and hence
preserve the stagewise essential images. Each image is full,
so all compatibility maps and their higher coherences between
objects in the images lie in those subcategories. Thus the
essential image of the restriction on the limit is precisely
the objects whose projections belong to the stagewise essential
images. This proves the criterion for affine \(S\).

For general \(S\), take a finite affine open cover. Any smooth
scheme over one of these opens is also smooth and finitely
presented over \(S\), and the two diagrams used to test it agree
by base change. Conversely, evaluations on \(X\times_S V_S\)
have Zariski descent in \(X\). The colimits of spectra defining
the test preserve finite limits, so the tests on the cover imply
the test over \(S\) by Zariski excision. Membership in the essential
image is also detected on this cover. By full faithfulness,
membership is equivalent to the unit
\(H\to\lambda^*\lambda_\sharp H\) being an equivalence.
This unit commutes with open restriction, as in the proof of
Theorem~\ref{thm:q417-full-scaling}. The open restrictions jointly
detect equivalences. The affine criterion now gives the result.
\end{proof}

\subsection[The stable comparison]{The stable comparison\texorpdfstring{ \cite[Construction~4.0.19]{BPOh}}{}}
\label{sec:ms-log-comparison}

Let \(S\) be a quasi-compact, quasi-separated scheme. Following
\cite[Construction~4.0.19]{BPOh}, let \(\mathcal E_S\) be the
full subcategory of presheaves of spectra on \(\Sm/S\) satisfying
Zariski descent and elementary blow-up excision. A closed immersion
\(Z\hookrightarrow X\) in \(\Sm/S\) is elementary if Zariski locally
on \(X\) it is the zero-section inclusion of \(Z\) into
\(\A^n_Z\), with a possible disjoint smooth summand in the ambient
scheme. Excision means that the presheaf is zero on the empty scheme
and sends the square with vertices \(X\), \(Z\), the blow-up of \(X\)
along \(Z\), and its exceptional divisor to a pullback square
\cite[Definition~4.26]{KMSramification}. Denote by \(m_S(X)\)
the represented suspension spectrum localized into \(\mathcal E_S\).
This category has the symmetric monoidal structure induced by
products over \(S\), so \(m_S(S)\) is its tensor unit.

Put \(T=\cofib(m_S(S)\xrightarrow{1}m_S(\Pp^1_S))\), where the arrow
comes from the \(1\)-section. The comparison uses Annala and Iwasa's
category \(\mathrm{MS}_S=\mathcal E_S[T^{-1}]\)
\cite[Construction~4.0.19]{BPOh}. Here formal inversion
means the universal presentable symmetric monoidal category
under \(\mathcal E_S\) in which the image of \(T\) has a tensor
inverse. It is modeled by symmetric \(T\)-spectra
\cite[Proposition~1.3.14]{AIMotivic}. Likewise, put
\(T_{\log}=\cofib(M(S)\xrightarrow{1}M(\Pp^1_S))\) in
\(\mathcal L_S\), using the projective line with trivial log structure.
Then \(\log\SH(S)=\mathcal L_S[T_{\log}^{-1}]\)
\cite[Proposition~3.3.6 and Corollary~3.4.7]{BPOh}.

To construct the comparison, give each ordinary smooth scheme the
trivial log structure. The resulting functor carries ordinary
Zariski distinguished squares to strict Nisnevich distinguished
squares. Ordinary smooth blow-up excision in \(\mathcal L_S\)
\cite[Theorem~3.3.7]{BPOh} gives the elementary blow-up relations as
well. Thus the assignment \(m_S(X)\mapsto M(X)\) extends to a
colimit-preserving symmetric monoidal functor
\(F_0:\mathcal E_S\to\mathcal L_S\). It sends \(T\) to \(T_{\log}\).
Formal inversion gives the adjunction of
\cite[Construction~4.0.19]{BPOh}
\[
  \lambda_\sharp:\mathrm{MS}_S
      \;\rightleftarrows\;\log\SH(S):\lambda^* .
\]
Binda, Park, and \O stv\ae r expect the right adjoint
\(\lambda^*\) to be fully faithful for every such \(S\).

\begin{corollary}\label{cor:ms-log-comparison}
For every quasi-compact, quasi-separated scheme \(S\), the
right adjoint \(\lambda^*:\log\SH(S)\to\mathrm{MS}_S\) in
\cite[Construction~4.0.19]{BPOh} is fully faithful.
\end{corollary}

\begin{proof}
Write \(G_0\) for the right adjoint of \(F_0\). By adjunction
on represented spectra, \(G_0E\) is the presheaf
\(X\mapsto E(X,\varnothing)\). It satisfies elementary blow-up
excision because mapping into \(E\) takes the ordinary smooth
blow-up pushout squares in \(\mathcal L_S\) to pullback squares
of spectra. The presheaf is also a Nisnevich sheaf. Thus \(G_0\)
and the restriction in Theorem~\ref{thm:q417-full-scaling} have
the same underlying presheaf. Both \(\mathcal E_S\) and the
category of ordinary Nisnevich sheaves of spectra are full
subcategories of the presheaf category. The theorem therefore
implies that \(G_0\) is fully faithful, so \(F_0\) is a reflective
localization, meaning a left adjoint with fully faithful right adjoint.

Apply \cite[Lemma~1.5.4]{AIMotivic} to the localization
\(F_0:\mathcal E_S\to\mathcal L_S\), using \(\mathcal E_S\)
as the acting symmetric monoidal category. The action on
\(\mathcal L_S\) is tensor product after \(F_0\), because
\(F_0\) is symmetric monoidal. The lemma gives a
localization on the categories of symmetric spectra for
\(T\) and \(T_{\log}\) respectively. These categories
are the formal tensor inversions in the displayed
adjunction by \cite[Proposition~1.3.14]{AIMotivic}.
Consequently its right adjoint is fully faithful.
\end{proof}

Over normal bases, Koizumi, Miyazaki, and Saito also factor the
stable comparison through their category of motivic spectra with
modulus. They prove full faithfulness for the right adjoint from
logarithmic spectra to spectra with modulus and construct the
adjunction between \(\mathrm{MS}_S\) and spectra with modulus
\cite[Corollary~4.25 and Theorem~4.27]{KMSramification}.

\section{Logarithmicity of algebraic cobordism}\label{sec:cobordism}

The ordinary stable motivic category \(\SH(S)\) is obtained from
\(\mathcal N_S\) by imposing invariance under product with \(\A^1_S\)
and by inverting the pointed projective line, as in
\cite{MorelVoevodsky,Voevodsky}.

We write \(\wedge\) for the smash product, \(\one_S\) for its unit, and \(M_S(X)=\Sigma^\infty_+X\) for
the suspension spectrum of a smooth \(S\)-scheme with a disjoint base point. For integers \(a,b\), the
notation \(\Sigma^{a,b}\) means suspension by \(a-b\) copies of the simplicial circle and \(b\) copies
of \((\Gm,1)\), using inverse suspensions for negative exponents. Set \(T=\Sigma^{2,1}\one_S\), the
suspension spectrum of the pointed projective line. We write \(T^q\wedge E=\Sigma^{2q,q}E\) for any
integer \(q\). For \(E\in\SH(S)\), its sheaf of spectra in Tate degree \(q\) is
\[
 E^{(q)}(X)=\Map_{\SH(S)}\bigl(M_S(X),T^q\wedge E\bigr).
\]
For a stable motivic spectrum, these are the sheaves in the
\(\Pp^1\)-spectrum levels for \(q\geq0\) in \cite[Definition~5.29]{ParkI}.
The suspensions with \(q<0\) will also satisfy our comparison.

Voevodsky's algebraic cobordism spectrum \(\MGL_S\) is constructed from Thom spaces of the tautological
vector bundles on Grassmannians \cite[Section~6.3]{Voevodsky}\cite[Section~3]{HoyoisCobordism}. For a vector bundle
\(V\), its Thom space \(\Th(V)\) is \(V/(V\setminus 0_V)\), where \(0_V\) is the zero section. In the
construction, the spectrum of the Thom space of a rank \(d\) tautological bundle is desuspended by
\(\Sigma^{2d,d}\). The colimit uses the inclusions of Grassmannians and stabilization by direct sum
with a trivial line. These constructions give the canonical pullback equivalence
\(f^*\MGL_S\simeq\MGL_{S'}\) for every \(f\colon S'\to S\) \cite[Section~7]{HoyoisCobordism}.

For \(S\in\ParkSch\), use the logarithmicity comparison
from Section~\ref{sec:toric} at every Tate degree \(q\geq0\).
Park asks whether it holds for \(\MGL_S\) over every regular
\(S\in\ParkSch\) \cite[Conjecture~7.16]{ParkI}.

\begin{theorem}[Logarithmicity of cobordism {\cite[Conjecture~7.16]{ParkI}}]\label{thm:park-mgl}
For every \(S\in\ParkSch\), every integer \(q\), every \(r\geq 0\), and every \(X\in\Sm/S\), the canonical comparison
\[
 \MGL_S^{(q)}(X)\longrightarrow
 \Log\bigl(\MGL_S^{(q)}\bigr)(X\times\bx_S^r)
\]
is an equivalence. In particular, the conjecture holds with integral
coefficients over all regular bases in \(\ParkSch\)
\cite[Conjecture~7.16]{ParkI}.
\end{theorem}

\subsection{Representing the comparison by a cofiber}

An object \(A\) of a symmetric monoidal category is dualizable if there is an object \(A^\vee\) with
evaluation \(A^\vee\wedge A\to\one\) and coevaluation \(\one\to A\wedge A^\vee\) satisfying the two
triangle identities. In particular, for dualizable \(A\),
\(\Map(B\wedge A,E)\simeq\Map(B,A^\vee\wedge E)\).
Motivic Atiyah duality says that \(M_S(Y)\) is dualizable for a smooth proper \(p\colon Y\to S\) and gives the formula
\begin{equation}\label{eq:duality}
 M_S(Y)^\vee\simeq p_\#\bigl(S^{-T_p}\bigr).
\end{equation}
Here \(p_\#\) is the left adjoint of pullback along \(p\), sending \(M_Y(U)\) to \(M_S(U)\) for \(U\)
smooth over \(Y\). The bundle \(T_p\) is the relative tangent bundle. For a vector bundle \(V\) over
\(Y\), \(S^V\) is the suspension spectrum in \(\SH(Y)\) of its Thom space over \(Y\) and \(S^{-V}\) its
smash inverse. Formula~\eqref{eq:duality} is smooth proper duality in the form established in
\cite[Sections~2--3]{HoyoisDuality}.

Define a simplicial object and two associated objects in \(\SH(S)\) by
\begin{align}
 B_{r,m,S}
   &=\colim_{\Sigma\in I_{r,m}^\op} M_S(Y_{\Sigma,S})^\vee,
       \label{eq:degree}\\
 B_{r,S}&=|B_{r,\bullet,S}|,\qquad
 D_{r,S}=\cofib\bigl(\one_S\longrightarrow B_{r,S}\bigr).
       \label{eq:D}
\end{align}
The simplicial maps of \eqref{eq:degree} use the same indexing construction as \eqref{eq:park}, applied
to the contravariant functor \(Y\mapsto M_S(Y_S)^\vee\) on the smooth proper toric schemes and the
lifts described above. Equivalently, one takes the left Kan extension of this functor along the
opposite projection to \(\Delta^\op\). The map from the unit is induced by the duals of the structure
maps \(M_S(Y_{\Sigma,S})\to\one_S\). These maps are compatible with every lift, so the unit maps are
compatible with the simplicial structure.

\begin{proposition}\label{prop:comparison}
For \(E\in\SH(S)\), \(q\in\Z\), and \(X\in\Sm/S\), the cofiber of the canonical logarithm comparison on
\(E^{(q)}\) is naturally
\[
 \Map_{\SH(S)}\bigl(M_S(X),T^q\wedge E\wedge D_{r,S}\bigr).
\]
Moreover, \(f^*D_{r,S}\simeq D_{r,S'}\) for every morphism \(f\colon S'\to S\) of base schemes. In
particular, for fixed \(r\), the comparisons for all \(q\geq 0\) and all \(X\) are equivalences if and
only if \(E\wedge D_{r,S}\simeq0\).
\end{proposition}

\begin{proof}
The spectrum \(M_S(X)\) is compact, that is, mapping out of it commutes with filtered colimits. The
suspension spectra of smooth schemes and their shifts \(\Sigma^{a,b}\) are compact generators of
\(\SH(S)\) \cite[Section~2.1]{HoyoisCobordism}\cite[Section~4.1]{BachmannHoyois}. The mapping spectrum functor
\(\Map_{\SH(S)}(M_S(X),-)\) is exact and preserves all coproducts, since a coproduct is a filtered
colimit of finite sums. It therefore preserves all colimits.

For a given \(Y=Y_{\Sigma,S}\), the equality \(M_S(X\times_S Y)\simeq M_S(X)\wedge M_S(Y)\) and duality
give a natural equivalence
\[
 E^{(q)}(X\times_S Y)\simeq
 \Map_{\SH(S)}\bigl(M_S(X),T^q\wedge E\wedge M_S(Y)^\vee\bigr).
\]
The equivalence is contravariantly natural in \(Y\), including for maps that lie over a morphism in
\(\Delta\). Taking the colimits of \eqref{eq:park} gives
\[
 \Log E^{(q)}(X\times\bx_S^r)
 \simeq\Map_{\SH(S)}\bigl(M_S(X),T^q\wedge E\wedge B_{r,S}\bigr).
\]
Naturality identifies the projection-induced comparison with the map induced by \(\one_S\to B_{r,S}\).
Taking its cofiber proves the formula. Also, \(f^*\) is symmetric monoidal and preserves colimits. It
sends \(M_S(Y_{\Sigma,S})\) to \(M_{S'}(Y_{\Sigma,S'})\) and preserves their duals and all the
structure maps. This proves the pullback assertion.

If the comparisons are equivalences for \(q\geq0\), the formula implies that all mapping spectra from
\(M_S(X)\wedge T^{-q}\) to \(E\wedge D_{r,S}\) vanish. The objects \(M_S(X)\wedge T^{-q}\) for
\(q\geq0\), together with ordinary suspensions and desuspensions, generate \(\SH(S)\) by
\(\Pp^1\)-stabilization. To see this from the bigraded generators above, given \(\Sigma^{a,b}M_S(X)\)
choose \(q\geq0\) with \(q+b\geq0\). Express \(T^{q+b}\wedge M_S(X)\) by taking the iterated cofibers
of the sections at infinity in a product of projective lines with \(X\), then apply \(T^{-q}\) and the
ordinary shift \(a-2b\). Thus \(E\wedge D_{r,S}\simeq0\). Conversely, its vanishing gives the
comparison for every integer \(q\).
\end{proof}

\subsection{The calculation over fields}

Let \(k\) be a perfect field, and let \(\HZ_k\) be the spectrum that represents integral motivic
cohomology. The proof of \cite[Theorem~B.4]{BachmannHoyois} identifies the spectrum of
\cite{Spitzweck} over a field with the one used in \cite[Section~4]{HoyoisCobordism}. For \(q\geq0\), the
sheaf \(\HZ_k^{(q)}\) is the spectrum associated to the complex
\(R\Gamma_{\mathrm{mot}}(-,\Z(q))[2q]\). In this notation \(\Z(q)\) is the weight \(q\) motivic
complex, whose hypercohomology gives the motivic cohomology groups, and \([2q]\) is the indicated
cohomological shift. For \(q<0\) the sheaf \(\HZ_k^{(q)}\) vanishes, since motivic cohomology of smooth
schemes over a field vanishes in negative weights
\cite[Corollary~4.26(1)]{HoyoisCobordism}.

Park proves that \(R\Gamma_{\mathrm{mot}}(-,\Z(q))\) is logarithmic for \(q\geq0\)
over every perfect field \cite[Theorem~7.9]{ParkI}. The calculation of Chow groups
used there \cite[Theorem~7.7]{ParkI} is proved in \cite[Theorem~1.1]{ParkII}.
Shifts commute with the colimits in \eqref{eq:park}.
Thus all the sheaves \(\HZ_k^{(q)}\) are logarithmic. By Proposition~\ref{prop:comparison},
\begin{equation}\label{eq:HZvanishes}
 D_{r,k}\wedge\HZ_k\simeq0.
\end{equation}
We use the homotopy \(t\)-structure of \(\SH(k)\) \cite[Section~2.1]{HoyoisCobordism}. Its nonnegative part
is the subcategory generated under colimits and extensions by \(\Sigma^{a,b}M_k(U)\) for \(U\) smooth
and \(a-b\geq0\). We call an object \(c\)-connective if it belongs to the \(c\)-fold ordinary
suspension of this subcategory. The subcategories of \(c\)-connective objects are closed under colimits
and extensions. Suspension by \(\Sigma^{a,b}\) changes the lower bound by \(a-b\). The smash product of
a \(c\)-connective and a \(d\)-connective object is \((c+d)\)-connective. This follows on the
generators from the product of smooth schemes and then for all objects from preservation of colimits
and extensions in both variables.

\begin{lemma}\label{lem:connectivity}
For every field \(k\) and integer \(r\geq0\), the spectrum \(D_{r,k}\wedge\MGL_k\) is \((-r)\)-connective.
\end{lemma}

\begin{proof}
Let \(p\colon Y\to\Spec k\) be smooth proper of pure dimension \(d\). Choose a finite Zariski open
cover of \(Y\) trivializing \(T_p\). For a nonempty finite intersection \(U\) of opens in this cover,
let \(p_U\colon U\to\Spec k\) be the structure map. Over \(U\), the inverse tangent Thom sphere is
\(\Sigma^{-2d,-d}\one_U\). Applying \((p_U)_\#\) gives \(\Sigma^{-2d,-d}M_k(U)\) by the projection
formula, so it is \((-d)\)-connective. The dual \(p_\#S^{-T_p}\) in \eqref{eq:duality} is the colimit
of these terms along the Čech diagram of the cover. Indeed, descent first writes \(\one_Y\) as the
colimit of the spectra of the intersections over \(Y\). Smashing that expression with \(S^{-T_p}\) and
applying \(p_\#\) gives the claimed diagram by the projection formula for the open immersions. The dual
is therefore \((-d)\)-connective.

It follows for the degree \(m\) terms that \(B_{r,m,k}\) is \((-r-m)\)-connective. We check directly
that their realization is \((-r)\)-connective. The \(m\)-th latching object \(L_mB\) is the colimit of
the lower-degree terms prescribed by the proper degeneracies into degree \(m\), namely the surjections
\([m]\to[j]\) with \(j<m\). It is \((-r-m+1)\)-connective. Including the \(m\)-simplices in the
realization adds a cofiber
\[
 \Sigma^{m,0}\cofib\bigl(L_mB\longrightarrow B_{r,m,k}\bigr).
\]
This is the usual skeletal filtration of the realization, obtained by attaching each \(m\)-simplex
along its boundary and its degenerate simplices. The displayed cofiber is \((-r)\)-connective, also for
\(m=0\) with \(L_0B=0\). Induction on the skeleta and passage to their colimit prove the bound for
\(B_{r,k}\). Since \(\one_k\) is connective, \eqref{eq:D} gives the same bound for \(D_{r,k}\).
Finally, \(\MGL_k\) is connective \cite[Corollary~3.9]{HoyoisCobordism}. The smash-product bound gives
the result.
\end{proof}

\begin{proposition}\label{prop:field}
If \(k\) is a perfect field, then \(D_{r,k}\wedge\MGL_k\simeq0\) for every \(r\geq0\).
\end{proposition}

\begin{proof}
Put \(D=D_{r,k}\), \(R=\MGL_k\), and \(N=D\wedge R\). By \eqref{eq:HZvanishes},
\(D\wedge\HZ_k\simeq0\), and by Lemma~\ref{lem:connectivity}, \(N\) is bounded below. Hoyois's
vanishing lemma \cite[Lemma~7.10]{HoyoisCobordism} says that if an object \(F\) of \(\SH(k)\) satisfies
\(F\wedge\HZ_k\simeq0\), then every map from \(F\) to a bounded-below weak \(\MGL_k\)-module is zero in
the homotopy category. A weak module is a module in the homotopy category. In particular, its action
satisfies the unit identity there. The free module \(N\) has such an action by multiplication on the
\(R\)-factor. Apply the lemma to \(F=D\) and the unit-induced map \(D\to D\wedge R\). Smashing with
\(R\) and composing with multiplication gives
\[
 D\wedge R\longrightarrow D\wedge R\wedge R
       \xrightarrow{\ \mathrm{id}_D\wedge\mu_R\ }D\wedge R.
\]
The composite is the identity by the unit axiom, and it is zero because the first map is zero in the
homotopy category. Thus \(D\wedge R\simeq0\).
\end{proof}

\begin{proof}[Second proof of Proposition~\ref{prop:field}]
Put \(R=\MGL_k\). The Lazard ring is the universal coefficient ring
for one-dimensional commutative formal group laws.
Choose the homogeneous polynomial generators \(t_i\) in degree \(i\) used in the Hopkins--Morel
isomorphism. The orientation of \(R=\MGL_k\) gives the corresponding classes in bidegree \((2i,i)\).
Write \(Q_0=R\) and define \(Q_i\) in the category of \(R\)-modules by the cofiber sequence
\[
 \Sigma^{2i,i}Q_{i-1}\xrightarrow{t_i}Q_{i-1}\longrightarrow Q_i
 \qquad(i\geq1).
\]
Here multiplication is the map obtained from the action of the indicated class. The integral
Hopkins--Morel isomorphism identifies the underlying spectrum of the quotient \(Q_\infty=\colim_i Q_i\)
with \(\HZ_k\) \cite[Theorem~7.2]{BachmannHopkins}.

Let \(N_i=D_{r,k}\wedge Q_i\), so that \(N_0=D_{r,k}\wedge R\). By Lemma~\ref{lem:connectivity},
\(N_0\) is \((-r)\)-connective. The cofiber sequences, after smashing with \(D_{r,k}\), show by
induction that all \(N_i\) have this same lower bound. The functor \(\tau_{\leq -r}\) taking homotopy
truncation in degrees at most \(-r\) is left adjoint to the inclusion of the right orthogonal to the
\((1-r)\)-connective objects. It therefore preserves colimits into that subcategory. The first term in
the cofiber sequence for \(N_{i-1}\to N_i\) is \(i-r\)-connective. Mapping the sequence into any object
truncated in degrees at most \(-r\) gives an equivalence on mapping spaces from \(N_i\) and \(N_{i-1}\)
by orthogonality. Hence \(N_{i-1}\to N_i\) induces an equivalence on \(\tau_{\leq -r}\).
By \eqref{eq:HZvanishes},
\[
 \colim_i N_i\simeq D_{r,k}\wedge\HZ_k\simeq0.
\]
Thus \(\tau_{\leq -r}N_0\simeq0\) and \(N_0\) is \((1-r)\)-connective. Repeat the argument with
this improved bound. For every integer \(b\), \(N_0\) is \(b\)-connective. The homotopy \(t\)-structure
over a field is left complete, meaning that the canonical map from an object to the inverse limit
of its truncations \(\tau_{\leq n}\) is an equivalence \cite[Corollary~2.4]{HoyoisCobordism}.
All these truncations of \(N_0\) are zero, so \(D_{r,k}\wedge\MGL_k\simeq0\).
\end{proof}

\subsection{Passage to the base and logarithmicity for modules}

\begin{proof}[Proof of Theorem~\ref{thm:park-mgl}]
Work first over \(S_0=\Spec\Z\). For a point \(s\) of \(S_0\), let \(\kappa(s)\) denote its residue
field and \(i_s\colon\Spec\kappa(s)\to S_0\) the canonical morphism. The two pullback formulas for
\(D_{r,S_0}\) and \(\MGL_{S_0}\) give
\[
 i_s^*(D_{r,S_0}\wedge\MGL_{S_0})
   \simeq D_{r,\kappa(s)}\wedge\MGL_{\kappa(s)} .
\]
The residue fields of \(S_0\) are \(\mathbb Q\) and the prime fields \(\mathbb F_p\). They are perfect,
so Proposition~\ref{prop:field} makes every displayed pullback zero. Residue field pullbacks detect
zero objects of \(\SH(S_0)\) by \cite[Proposition~B.3]{BachmannHoyois}. Consequently
\(D_{r,S_0}\wedge\MGL_{S_0}\simeq0\). Pulling back along \(S\to S_0\) gives
\(D_{r,S}\wedge\MGL_S\simeq0\) for every \(S\in\ParkSch\). Proposition~\ref{prop:comparison} now
yields the canonical equivalence at all the stated \(X,r,q\).
\end{proof}

The passage through \(\Spec\Z\) allows the field calculation to use
perfect fields throughout. In particular the theorem applies
integrally to singular bases, bases with imperfect residue fields,
and bases of mixed characteristic.

\begin{corollary}
Every \(\MGL_S\)-module in \(\SH(S)\), for \(S\in\ParkSch\),
has logarithmic sheaves in every Tate degree.
\end{corollary}
\begin{proof}
For such a module \(E\), its unit and action make \(E\wedge D_{r,S}\)
a retract of \(E\wedge\MGL_S\wedge D_{r,S}\). The latter
vanishes by the preceding proof, and the formula for the
comparison cofiber in Proposition~\ref{prop:comparison}
gives the result.
\end{proof}

Park's conjecture about modules concerns
\(\mathrm{MS}_{\Nis}(S)\), the Nisnevich version of the
projective-line stable category with elementary blow-up
excision before imposing affine-line invariance
\cite[\S5 and Conjecture~7.17]{ParkI}. The corollary concerns
the ordinary affine-line invariant category.

\section{Algebraic \texorpdfstring{\(K\)}{K}-theory on qcqs schemes}\label{sec:k-theory}

Take \(S\in\Sch\) with the trivial log structure. Write
\(K(S)\) for the nonconnective algebraic \(K\)-theory spectrum of
perfect complexes on \(S\). A perfect complex is a complex locally
quasi-isomorphic to a bounded complex of vector bundles of finite rank.
Nonconnective \(K\)-theory includes the negative \(K\)-groups.

Use the logarithmic stable motivic category \(\logSH(S)\) of
Section~\ref{sec:spectra-definitions}. If
\(f_S:S\to\Spec\Z\) is the structural morphism, put
\begin{equation}\label{eq:absolute}
 E_S=f_S^*\mathrm{logKGL}_{\Z},\qquad
 \log K(S)=\Map_{\logSH(S)}(\one_S,E_S).
\end{equation}
Here \(\mathrm{logKGL}_{\Z}\) is the Bott-periodic spectrum
constructed from algebraic \(K\)-theory over the regular base
\(\Spec\Z\). Its projective-line structure maps are induced by
the Bott class, the reduced \(K\)-theory class of the tautological
line bundle on the projective line. See
\cite[Definition~6.5.6]{BPOh}. Formula \eqref{eq:absolute} is the
instance of \cite[Definition~6.5.9]{BPOh} used for ordinary schemes.
The predicted equivalence in \cite[Conjecture~6.5.10]{BPOh} is
\begin{equation}\label{eq:conjecture}
 \log K(S)\simeq K(S).
\end{equation}
\begin{theorem}[Comparison with algebraic \(K\)-theory {\cite[Conjecture~6.5.10]{BPOh}}]\label{thm:log-k}
For every qcqs scheme \(S\), there is a natural equivalence of
spectra \(K(S)\simeq\log K(S)\), with \(\log K\) defined by
\eqref{eq:absolute}. Thus \eqref{eq:conjecture} follows from
Park's representability theorems and finite-presentation continuity.
\end{theorem}

For \(S\in\ParkSch\), Park constructs a logarithmic extension
\(\Log K\) of \(X\mapsto K(X)\) by the construction in
Section~\ref{sec:toric}, and its Bott-periodic spectrum
\(\Log\KGL_S\). We use the following results from \cite{ParkI}.
\begin{enumerate}
\item \(\Log K(X)=K(X)\) for ordinary smooth
\(S\)-schemes \(X\) by \cite[Theorem~7.10 and Definition~5.26]{ParkI}.
\item If \(S\) is regular, then
\(\Log\KGL_S\simeq\mathrm{logKGL}_S\) by \cite[Theorem~8.11]{ParkI}.
\item For a morphism \(g:S'\to S\) in this class,
\(g^*\Log\KGL_S\simeq
\Log\KGL_{S'}\) by \cite[Proposition~8.12]{ParkI}.
\item For ordinary smooth \(X/S\),
\[
 \Map_{\logSH(S)}(\Sigma^\infty_{\Pp^1}X_+,
       \Log\KGL_S)
 \simeq\Log K(X)
\]
by \cite[Theorem~8.13]{ParkI}.
\end{enumerate}
We use these comparisons with their compatibility under pullback.
The Chow group calculation \cite[Theorem~7.7]{ParkI} used in
\cite[Theorem~7.10]{ParkI} is proved in \cite[Theorem~1.1]{ParkII}.

\begin{proposition}\label{prop:noetherian}
Theorem~\ref{thm:log-k} holds for separated Noetherian schemes of
finite Krull dimension, naturally for their morphisms.
\end{proposition}
\begin{proof}
Apply \cite[Theorem~8.11]{ParkI} over \(\Spec\Z\) and
\cite[Proposition~8.12]{ParkI} along \(f_S\) to obtain
\[
 E_S=f_S^*\mathrm{logKGL}_{\Z}
   \simeq f_S^*\Log\KGL_{\Z}
   \simeq\Log\KGL_S.
\]
The identity \(S\to S\) is smooth. Taking \(X=S\) in
\cite[Theorem~8.13]{ParkI} and using \cite[Theorem~7.10]{ParkI} gives
\[
 \Map_{\logSH(S)}(\one_S,E_S)
   \simeq\Log K(S)\simeq K(S).
\]
The comparisons commute with pullback.
\end{proof}

\subsection{The comparison for arbitrary qcqs schemes}

We first make the reconstruction of natural transformations explicit.
A contravariant functor \(F\) from qcqs schemes to spectra is
\emph{affine-continuous} if every filtered colimit of rings
\(A=\operatorname*{colim}_i A_i\) induces an equivalence
\[
 \operatorname*{colim}_i F(\Spec A_i)\simeq F(\Spec A).
\]
A Zariski sheaf satisfies descent for open covers. For qcqs schemes
this can be checked by the empty-object
condition and Mayer--Vietoris for covers by two quasi-compact open
subschemes. Write \(\mathrm{Aff}^{\mathrm{fp}}_{\Z}\) for
the category of affine spectra of finitely presented \(\Z\)-algebras.

\begin{lemma}[Reconstruction]\label{lem:reconstruction}
Restriction to \(\mathrm{Aff}^{\mathrm{fp}}_{\Z}\) is fully
faithful on the category of affine-continuous Zariski sheaves of
spectra on qcqs schemes. In particular, a natural comparison between
two such sheaves on these finite-presentation affines extends
uniquely, up to a contractible space of choices, to a natural
comparison on all qcqs schemes.
\end{lemma}
\begin{proof}
For a ring \(A\), let \(I_A\) be the category of maps
\(B\to A\) with \(B\) a finitely presented \(\Z\)-algebra.
This category is filtered. Tensor products give a common target for
two objects, and a finite set of relations gives a common target
equalizing two parallel arrows. Its colimit of rings \(B\) is
\(A\): finitely many elements are represented by a polynomial ring,
and each finite set of relations can be imposed at a later object.
Affine continuity therefore gives the functorial formula
\begin{equation}\label{eq:kan-formula}
 F(\Spec A)\simeq
 \operatorname*{colim}_{(B\to A)\in I_A}F(\Spec B).
\end{equation}
This is the left Kan extension formula on rings: its universal
property identifies natural transformations between the functors on
all rings with natural transformations between their restrictions to
finitely presented rings. The statement concerns the spaces of
natural transformations, so it retains their homotopies and higher
compatibilities.

Affine open subschemes form a basis for the Zariski topology.
Consequently restriction to affines is fully faithful on Zariski
sheaves. Descent reconstructs a sheaf on each scheme from this basis,
and compatibility with a morphism is checked after refining its
source and target by affine open subschemes. Combining this with
\eqref{eq:kan-formula} proves the assertion.
\end{proof}

\begin{proof}[Proof of Theorem~\ref{thm:log-k}]
First suppose \(S=\Spec A\). Write
\(A=\operatorname*{colim}_i A_i\), where the \(A_i\) are
finitely presented \(\Z\)-algebras. Include \(A_0=\Z\)
as the initial stage. No injectivity of the ring maps is required.
Every \(S_i=\Spec A_i\) is separated, Noetherian, and of
finite Krull dimension. Proposition~\ref{prop:noetherian}
gives compatible equivalences \(K(S_i)\simeq\log K(S_i)\).
Apply Proposition~\ref{prop:continuity} with
\(C_0=\one_{\Z}\) and
\(F_0=\mathrm{logKGL}_{\Z}\). We obtain
\begin{align}\label{eq:affine-chain}
 \log K(\Spec A)
 &\simeq\operatorname*{colim}_i\log K(\Spec A_i)\notag\\
 &\simeq\operatorname*{colim}_i K(A_i)
 \simeq K(A).
\end{align}
The last equivalence is continuity of nonconnective \(K\)-theory.
Perfect complexes on a filtered colimit of rings descend to a
stage, as do their maps and retract data. Nonconnective
\(K\)-theory preserves filtered colimits of the resulting small
stable categories \cite[\S9]{BGT}. This also follows from
\cite[\S3.20 and \S6]{TT}.

We now pass to any qcqs \(S\). Both sides of
\eqref{eq:conjecture} satisfy Zariski descent. For ordinary
nonconnective \(K\)-theory this is \cite[Theorem~8.1]{TT}.
For \(\log K\), if \(S=U\cup V\) is a cover by
quasi-compact open subschemes, the defining sheaf relations give
the cofiber sequence in \(\logSH(S)\)
\[
 M_S(U\cap V)\longrightarrow M_S(U)\oplus M_S(V)
                 \longrightarrow\one_S.
\]
Mapping into \(E_S\) gives a pullback square of spectra.
For an open immersion \(j:W\to S\), the smooth adjunction
identifies \(\Map(M_S(W),E_S)\) with
\(\Map(\one_W,j^*E_S)=\log K(W)\).
This proves the required Mayer--Vietoris property and hence
Zariski descent. Affine continuity follows from
Proposition~\ref{prop:continuity} for \(\log K\) and from continuity
of nonconnective \(K\)-theory for \(K\). Lemma~\ref{lem:reconstruction}
therefore extends the natural equivalence of
Proposition~\ref{prop:noetherian} on finite-presentation affines
to a natural equivalence on all qcqs schemes. The calculation
\eqref{eq:affine-chain} is its affine formula. This proves
Theorem~\ref{thm:log-k} as a statement about functors, with all
compatibilities under pullback.
\end{proof}

\subsection{Smooth schemes, coefficients, support, and completion}

Every presentable stable category has an action of spectra. For
\(F\in\logSH(S)\) and a spectrum \(R\), denote the resulting
object by \(F\otimes R\). The product of two ordinary spectra
is denoted by \(\wedge\). Here \(R\) need not be a ring
spectrum. These conventions specify exactly what is meant by
coefficients in this section.

\begin{lemma}\label{lem:spectral-coefficients}
If \(G\) is compact in a presentable stable category, then
for every object \(F\) and every spectrum \(R\) the natural
map
\[
 \Map(G,F)\wedge R\longrightarrow\Map(G,F\otimes R)
\]
is an equivalence.
\end{lemma}
\begin{proof}
For the sphere spectrum this is the identity. Mapping spectra are
exact in their target, so the assertion is preserved by shifts,
finite cofibers, and retracts in \(R\). It therefore holds for
finite spectra, meaning retracts of finite cell spectra built from
the sphere. Every spectrum is a filtered colimit of finite spectra.
The left side preserves these colimits, and the right side does so
by compactness of \(G\). This proves the assertion for every
\(R\).
\end{proof}

\begin{corollary}\label{cor:smooth-representability}
For an ordinary smooth finitely presented morphism \(p:X\to S\)
between qcqs schemes and every spectrum \(R\), there is a natural
equivalence
\begin{equation}\label{eq:smooth-coefficients}
 \Map_{\logSH(S)}(M_S(X),E_S\otimes R)
       \simeq K(X)\wedge R.
\end{equation}
In particular, the pullback of \(\mathrm{logKGL}_{\Z}\)
represents ordinary nonconnective \(K\)-theory on the ordinary
smooth site of every qcqs base.
\end{corollary}
\begin{proof}
The smooth adjunction \(p_\sharp\dashv p^*\) identifies
\(p_\sharp\one_X=M_S(X)\). The functor \(p_\sharp\) is induced by
composition of smooth pairs with \(p\), which respects the defining relations.
See \cite[Proposition~2.3.9 and \S2.4]{BPOh}.
Since \(p^*E_S=E_X\), it gives
\[
 \Map_{\logSH(S)}(M_S(X),E_S)
 \simeq\Map_{\logSH(X)}(\one_X,E_X)
 \simeq K(X).
\]
The last equivalence is Theorem~\ref{thm:log-k}.
Now \(M_S(X)\) is compact by Lemma~\ref{lem:compact}, so
Lemma~\ref{lem:spectral-coefficients} gives
\eqref{eq:smooth-coefficients}.
\end{proof}

For an integer \(n\geq2\), let \(\mathbb S/n\) be the
cofiber of multiplication by \(n\) on the sphere spectrum.
For any spectrum or object of a stable category, the notation
\(F/n\) means \(\cofib(n:F\to F)\). Thus
\(F/n=F\otimes\mathbb S/n\) in a category acted on by
spectra. Taking \(R=\mathbb S/n\) in
\eqref{eq:smooth-coefficients} gives the comparison with these
finite coefficients, even when \(n\) is not invertible on
\(S\). Taking \(R\) to be the rationalized sphere gives the
rational comparison. If \(R=H\Lambda\) is the
Eilenberg--MacLane spectrum of a ring \(\Lambda\), the formula
instead gives \(K(X)\wedge H\Lambda\). The last convention is
different from Moore-spectrum coefficients and is stated separately.

Let \(Z\) be a closed subset of a qcqs scheme \(S\), and
suppose \(U=S\setminus Z\) is quasi-compact. Define
\[
 Q_{S,U}=\cofib\bigl(M_S(U)\longrightarrow\one_S\bigr).
\]
This is compact. Let \(\mathrm{Perf}_Z(S)\) be the stable
category of perfect complexes on \(S\) whose restriction to
\(U\) is zero, and put \(K_Z(S)=K(\mathrm{Perf}_Z(S))\).
These definitions use the closed subset, rather than a choice of
scheme structure on it.

\begin{corollary}\label{cor:support}
With this notation, for every spectrum \(R\) there is a natural
equivalence
\[
 \Map_{\logSH(S)}(Q_{S,U},E_S\otimes R)
       \simeq K_Z(S)\wedge R.
\]
The comparison identifies the logarithmic fiber sequence for
\(U\subset S\) with the nonconnective \(K\)-theory localization
sequence.
\end{corollary}
\begin{proof}
Mapping out of the defining cofiber gives
\[
 \Map(Q_{S,U},E_S)
  \simeq\operatorname{fib}\bigl(\log K(S)\to\log K(U)\bigr).
\]
Naturality in Theorem~\ref{thm:log-k} identifies this with
\(\operatorname{fib}(K(S)\to K(U))\), which is \(K_Z(S)\)
by \cite[Theorem~7.4]{TT}. Tensoring spectra with \(R\)
preserves fiber sequences, and the same conclusion follows either
from this observation or from Lemma~\ref{lem:spectral-coefficients}
applied to the compact object \(Q_{S,U}\).
\end{proof}

The \emph{derived \(n\)-adic completion} of an object \(F\)
in a stable category is, by definition,
\(F^\wedge_n=\varprojlim_r F/n^r\), with the natural
reduction maps. This definition is made on spectra or objects of
the category, before taking homotopy groups.

\begin{corollary}\label{cor:completion}
For every qcqs scheme \(S\), every ordinary smooth finitely
presented \(X/S\), and every integer \(n\geq2\),
\[
 \Map_{\logSH(S)}(M_S(X),(E_S)^\wedge_n)
              \simeq K(X)^\wedge_n.
\]
\end{corollary}
\begin{proof}
Mapping into a target preserves inverse limits. At every finite
stage, Corollary~\ref{cor:smooth-representability} identifies the
mapping spectrum with \(K(X)/n^r\), compatibly with the
reduction maps. Taking their inverse limit gives the assertion.
No interchange of homotopy groups with inverse limits is required.
\end{proof}

The comparison applies to the scheme with trivial log structure.
It gives a useful calculation for the marked interval as well.
Logarithmic invariance gives \(\log K(\bx_S)=\log K(S)=K(S)\),
whereas the trivial locus of \(\bx_S\) has \(K\)-theory
\(K(\A^1_S)\). These may differ on a singular base, as discussed
in \cite[Remark~6.5.11]{BPOh}.

For a concrete singular example, let \(A=k[\epsilon]/(\epsilon^2)\)
for any field \(k\). The unit \(1+\epsilon t\) in
\(A[t]\), with inverse \(1-\epsilon t\), gives a nonzero
class in the kernel of
\(K_1(A[t])\to K_1(A)\) induced by \(t\mapsto0\).
Indeed, the determinant of this class is the nonconstant unit
\(1+\epsilon t\), so the class is nonzero and cannot come from
\(K_1(A)\). This verifies concretely why the logarithmic interval
must not be replaced by an ordinary affine line in the argument.

\end{document}